\documentclass{ajmlatex}
\usepackage{latexsym}

\relax 
\relax 
\relax 
\relax 
\newtheorem{conj}[theorem]{Conjecture}
\newtheorem{prop}[theorem]{Proposition}
\newtheorem{importnota}[theorem]{Important Notation}
\newtheorem{prblm}[theorem]{Problem}
\newtheorem{notation}[theorem]{Notation}
\newtheorem{defin}[theorem]{Definition}
\newtheorem{caution}[theorem]{Caution}
\newtheorem{remark}[theorem]{Remark}

\newtheorem{construction}[theorem]{Construction}

\newtheorem{example}[theorem]{Example}
\newtheorem{conclusion}[theorem]{Conclusion}
\newtheorem{triviality}[theorem]{Triviality}
\newtheorem{proto}[theorem]{Prototype Quasifibration}
\newtheorem{cauex}[theorem]{Cautionary Example}
\newtheorem{subth}{ }[theorem]
\newtheorem{ssubth}{ }[subth]

\input amssym.def
\input amssym.tex

\newcommand\nc{\newcommand}

\nc\BA{\begin{array}{c}}
\nc\EA{\end{array}}

\newcount\blopone 

\newcount\fone
\newcount\ftwo
\newcount\fthr
\newcount\eone
\newcount\etwo
\newcount\ethr

\newcount\gone
\newcount\gtwo
\newcount\gthr
\newcount\gfou
\newcount\gfiv
\newcount\gsix
\newcount\gsev

\newcount\xone
\newcount\xtwo
\newcount\xthr
\newcount\xfou
\newcount\xfiv
\newcount\xsix
\newcount\xsev

\newcount\hone
\newcount\htwo
\newcount\hthr
\newcount\hfou
\newcount\hfiv
\newcount\hsix
\newcount\hsev

\newcount\yone
\newcount\ytwo
\newcount\ythr
\newcount\yfou
\newcount\yfiv
\newcount\ysix
\newcount\ysev

\nc\nctr{\newcount}

\newcount\varone
\newcount\vartwo
\newcount\varthr
\newcount\varfou

\nctr{\ione}
\nctr{\itwo}
\nctr{\ithr}
\nctr{\ifou}
\nctr{\ifiv}
\nctr{\isix}
\nctr{\isev}
\nctr{\ieig}
\nctr{\inin}

\nctr{\Ione}
\nctr{\Itwo}
\nctr{\Ithr}
\nctr{\Ifou}
\nctr{\Ifiv}
\nctr{\Isix}
\nctr{\Isev}
\nctr{\Ieig}
\nctr{\Inin}

\nctr{\tione}
\nctr{\titwo}
\nctr{\tithr}
\nctr{\tifou}
\nctr{\tifiv}
\nctr{\tisix}
\nctr{\tisev}
\nctr{\tieig}
\nctr{\tinin}

\nctr{\tIone}
\nctr{\tItwo}
\nctr{\tIthr}
\nctr{\tIfou}
\nctr{\tIfiv}
\nctr{\tIsix}
\nctr{\tIsev}
\nctr{\tIeig}
\nctr{\tInin}

\nctr{\jone}
\nctr{\jtwo}
\nctr{\jthr}
\nctr{\jfou}
\nctr{\jfiv}
\nctr{\jsix}
\nctr{\jsev}
\nctr{\jeig}
\nctr{\jnin}

\nctr{\Jone}
\nctr{\Jtwo}
\nctr{\Jthr}
\nctr{\Jfou}
\nctr{\Jfiv}
\nctr{\Jsix}
\nctr{\Jsev}
\nctr{\Jeig}
\nctr{\Jnin}

\nctr{\tjone}
\nctr{\tjtwo}
\nctr{\tjthr}
\nctr{\tjfou}
\nctr{\tjfiv}
\nctr{\tjsix}
\nctr{\tjsev}
\nctr{\tjeig}
\nctr{\tjnin}

\nctr{\tJone}
\nctr{\tJtwo}
\nctr{\tJthr}
\nctr{\tJfou}
\nctr{\tJfiv}
\nctr{\tJsix}
\nctr{\tJsev}
\nctr{\tJeig}
\nctr{\tJnin}

\nc\siv[9]{\ione =#1 \multiply \ione by 10 \tione =0 \advance \tione by \ione
\itwo =#2 \multiply \itwo by 10
 \titwo =\tione \advance \titwo by \itwo
\ithr =#3 \multiply \ithr by 10
 \tithr =\titwo \advance \tithr by \ithr
\ifou =#4 \multiply \ifou by 10
 \tifou =\tithr \advance \tifou by \ifou
\ifiv =#5 \multiply \ifiv by 10
 \tifiv =\tifou \advance \tifiv by \ifiv
\isix =#6 \multiply \isix by 10
 \tisix =\tifiv \advance \tisix by \isix
\isev =#7 \multiply \isev by 10
 \tisev =\tisix \advance \tisev by \isev
\ieig =#8 \multiply \ieig by 10
 \tieig =\tisev \advance \tieig by \ieig
\inin =#9 \multiply \inin by 10
 \tinin =\tieig \advance \tinin by \inin 
}

\nc\sIv[9]{\Ione =#1 \multiply \Ione by 10 
\tIone =\tinin \advance \tIone by \Ione
\Itwo =#2 \multiply \Itwo by 10
 \tItwo =\tIone \advance \tItwo by \Itwo
\Ithr =#3 \multiply \Ithr by 10
 \tIthr =\tItwo \advance \tIthr by \Ithr
\Ifou =#4 \multiply \Ifou by 10
 \tIfou =\tIthr \advance \tIfou by \Ifou
\Ifiv =#5 \multiply \Ifiv by 10
 \tIfiv =\tIfou \advance \tIfiv by \Ifiv
\Isix =#6 \multiply \Isix by 10
 \tIsix =\tIfiv \advance \tIsix by \Isix
\Isev =#7 \multiply \Isev by 10
 \tIsev =\tIsix \advance \tIsev by \Isev
\Ieig =#8 \multiply \Ieig by 10
 \tIeig =\tIsev \advance \tIeig by \Ieig
\Inin =#9 \multiply \Inin by 10
 \tInin =\tIeig \advance \tInin by \Inin 
}

\nc\sjv[9]{\jone =#1 \multiply \jone by 10 \tjone =0 \advance \tjone by \jone
\jtwo =#2 \multiply \jtwo by 10
 \tjtwo =\tjone \advance \tjtwo by \jtwo
\jthr =#3 \multiply \jthr by 10
 \tjthr =\tjtwo \advance \tjthr by \jthr
\jfou =#4 \multiply \jfou by 10
 \tjfou =\tjthr \advance \tjfou by \jfou
\jfiv =#5 \multiply \jfiv by 10
 \tjfiv =\tjfou \advance \tjfiv by \jfiv
\jsix =#6 \multiply \jsix by 10
 \tjsix =\tjfiv \advance \tjsix by \jsix
\jsev =#7 \multiply \jsev by 10
 \tjsev =\tjsix \advance \tjsev by \jsev
\jeig =#8 \multiply \jeig by 10
 \tjeig =\tjsev \advance \tjeig by \jeig
\jnin =#9 \multiply \jnin by 10
 \tjnin =\tjeig \advance \tjnin by \jnin 
}

\nc\sJv[9]{\Jone =#1 \multiply \Jone by 10 
\tJone =\tjnin \advance \tJone by \Jone
\Jtwo =#2 \multiply \Jtwo by 10
 \tJtwo =\tJone \advance \tJtwo by \Jtwo
\Jthr =#3 \multiply \Jthr by 10
 \tJthr =\tJtwo \advance \tJthr by \Jthr
\Jfou =#4 \multiply \Jfou by 10
 \tJfou =\tJthr \advance \tJfou by \Jfou
\Jfiv =#5 \multiply \Jfiv by 10
 \tJfiv =\tJfou \advance \tJfiv by \Jfiv
\Jsix =#6 \multiply \Jsix by 10
 \tJsix =\tJfiv \advance \tJsix by \Jsix
\Jsev =#7 \multiply \Jsev by 10
 \tJsev =\tJsix \advance \tJsev by \Jsev
\Jeig =#8 \multiply \Jeig by 10
 \tJeig =\tJsev \advance \tJeig by \Jeig
\Jnin =#9 \multiply \Jnin by 10
 \tJnin =\tJeig \advance \tJnin by \Jnin 
}

\nc\icoo[1]{
\if#11\varone=0\varthr=\ione
\else\if#12\varone=\tione\varthr=\itwo
\else\if#13\varone=\titwo\varthr=\ithr
\else\if#14\varone=\tithr\varthr=\ifou
\else\if#15\varone=\tifou\varthr=\ifiv
\else\if#16\varone=\tifiv\varthr=\isix
\else\if#17\varone=\tisix\varthr=\isev
\else\if#18\varone=\tisev\varthr=\ieig
\else\if#19\varone=\tieig\varthr=\inin
\fi\fi\fi\fi\fi\fi\fi\fi\fi
}

\nctr\nten
\nctr\nele
\nctr\ntwe
\nctr\nthi
\nctr\nfou
\nctr\nfif
\nctr\nsix
\nctr\nsev
\nctr\neig

\nc\setnvar{\nten=10\nele=11\ntwe=12\nthi=13\nfou=14\nfif=15\nsix=16\nsev=17\neig=18}

\nc\Icoo[1]{\setnvar
\if#11\varone=0\varthr=\ione
\else\if#12\varone=\tione\varthr=\itwo
\else\if#13\varone=\titwo\varthr=\ithr
\else\if#14\varone=\tithr\varthr=\ifou
\else\if#15\varone=\tifou\varthr=\ifiv
\else\if#16\varone=\tifiv\varthr=\isix
\else\if#17\varone=\tisix\varthr=\isev
\else\if#18\varone=\tisev\varthr=\ieig
\else\if#19\varone=\tieig\varthr=\inin
\else\if#1\nten\varone=\tinin\varthr=\Ione
\else\if#1\nele\varone=\tIone\varthr=\Itwo
\else\if#1\ntwe\varone=\tItwo\varthr=\Ithr
\else\if#1\nthi\varone=\tIthr\varthr=\Ifou
\else\if#1\nfou\varone=\tIfou\varthr=\Ifiv
\else\if#1\nfif\varone=\tIfiv\varthr=\Isix
\else\if#1\nsix\varone=\tIsix\varthr=\Isev
\else\if#1\nsev\varone=\tIsev\varthr=\Ieig
\else\if#1\neig\varone=\tIeig\varthr=\Inin
\fi\fi\fi\fi\fi\fi\fi\fi\fi\fi\fi\fi\fi\fi\fi\fi\fi\fi
}

\nc\jcoo[1]{
\if#11\vartwo=0\varfou=\jone
\else\if#12\vartwo=\tjone\varfou=\jtwo
\else\if#13\vartwo=\tjtwo\varfou=\jthr
\else\if#14\vartwo=\tjthr\varfou=\jfou
\else\if#15\vartwo=\tjfou\varfou=\jfiv
\else\if#16\vartwo=\tjfiv\varfou=\jsix
\else\if#17\vartwo=\tjsix\varfou=\jsev
\else\if#18\vartwo=\tjsev\varfou=\jeig
\else\if#19\vartwo=\tjeig\varfou=\jnin
\fi\fi\fi\fi\fi\fi\fi\fi\fi
}

\nc\Jcoo[1]{\setnvar
\if#1\vartwo=0\varfou=\jone
\else\if#12\vartwo=\tjone\varfou=\jtwo
\else\if#13\vartwo=\tjtwo\varfou=\jthr
\else\if#14\vartwo=\tjthr\varfou=\jfou
\else\if#15\vartwo=\tjfou\varfou=\jfiv
\else\if#16\vartwo=\tjfiv\varfou=\jsix
\else\if#17\vartwo=\tjsix\varfou=\jsev
\else\if#18\vartwo=\tjsev\varfou=\jeig
\else\if#19\vartwo=\tjeig\varfou=\jnin
\else\if#1\nten\vartwo=\tjnin\varfou=\Jone
\else\if#1\nele\vartwo=\tJone\varfou=\Jtwo
\else\if#1\ntwe\vartwo=\tJtwo\varfou=\Jthr
\else\if#1\nthi\vartwo=\tJthr\varfou=\Jfou
\else\if#1\nfou\vartwo=\tJfou\varfou=\Jfiv
\else\if#1\nfif\vartwo=\tJfiv\varfou=\Jsix
\else\if#1\nsix\vartwo=\tJsix\varfou=\Jsev
\else\if#1\nsev\vartwo=\tJsev\varfou=\Jeig
\else\if#1\neig\vartwo=\tJeig\varfou=\Jnin
\fi\fi\fi\fi\fi\fi\fi\fi\fi\fi\fi\fi\fi\fi\fi\fi\fi\fi
}

\newsavebox\boxone
\newsavebox\boxtwo
\newsavebox\boxthr
\newsavebox\boxfou
\newsavebox\boxfiv
\newsavebox\boxsix
\newsavebox\boxsev
\newsavebox\boxeig
\newsavebox\boxnin
\newsavebox\boxten
\newsavebox\boxele
\newsavebox\boxtwe
\newsavebox\boxthi
\newsavebox\boxfot
\newsavebox\boxfit

\nc\uvariablearrowpos[3]{\eone =#1
\etwo =\eone \advance \etwo by -2 
\ethr =\eone \advance \ethr by 2
\gone =#2 \fone =#3
\ftwo =\fone \divide \ftwo by 2
\fthr =\fone \advance \fthr by -2
\gtwo =\gone \advance \gtwo by 2
\gthr =\gone \advance \gthr by 4
\gfou =\gone \multiply \gfou by 2
\advance \gfou by \fone
\gfiv =\gfou \divide \gfiv by 2 
\advance \gfou by -5
\divide  \gfou by 2
\gsev =\gone \advance \gsev by \fone
\gsix =\gone \advance \gsix by \fthr
}

\nc\rvariablearrowpos[3]{\eone =#2
\etwo =\eone \advance \etwo by -2 
\ethr =\eone \advance \ethr by 2
\gone =#1 \fone =#3
\ftwo =\fone \divide \ftwo by 2
\fthr =\fone \advance \fthr by -2
\gtwo =\gone \advance \gtwo by 2
\gthr =\gone \advance \gthr by 4
\gfou =\gone \multiply \gfou by 2
\advance \gfou by \fone
\gfiv =\gfou \divide \gfiv by 2 
\advance \gfou by -5
\divide  \gfou by 2
\gsev =\gone \advance \gsev by \fone
\gsix =\gone \advance \gsix by \fthr
}

\nc\vume[3]{\uvariablearrowpos{#1}{#2}{#3}
\put(\eone ,\gtwo){\line(0,1){\fthr}}
\put(\etwo ,\gtwo){\oval(4,4)[br]}
\put(\etwo ,\gsev){\oval(4,4)[br]}
\put(\etwo ,\gsix){\oval(4,4)[br]}
\put(\ethr ,\gtwo){\oval(4,4)[bl]}
\put(\ethr ,\gsev){\oval(4,4)[bl]}
\put(\ethr ,\gsix){\oval(4,4)[bl]}
}

\nc\vum[3]{\uvariablearrowpos{#1}{#2}{#3}
\put(\eone ,\gtwo){\line(0,1){\fthr}}
\put(\etwo ,\gtwo){\oval(4,4)[br]}
\put(\etwo ,\gsev){\oval(4,4)[br]}
\put(\ethr ,\gtwo){\oval(4,4)[bl]}
\put(\ethr ,\gsev){\oval(4,4)[bl]}
}

\nc\vue[3]{\uvariablearrowpos{#1}{#2}{#3}
\put(\eone ,\gone){\line(0,1){\fone}}
\put(\etwo ,\gsev){\oval(4,4)[br]}
\put(\etwo ,\gsix){\oval(4,4)[br]}
\put(\ethr ,\gsev){\oval(4,4)[bl]}
\put(\ethr ,\gsix){\oval(4,4)[bl]}
}

\nc\vua[3]{\uvariablearrowpos{#1}{#2}{#3}
\put(\eone ,\gone){\line(0,1){\fone}}
\put(\etwo ,\gsev){\oval(4,4)[br]}
\put(\ethr ,\gsev){\oval(4,4)[bl]}
}

\nc\vuhe[3]{\uvariablearrowpos{#1}{#2}{#3}
\put(\eone ,\gone){\line(0,1){\fone}}
\put(\etwo ,\gfou){\line(0,1){6}}
\put(\ethr ,\gfou){\line(0,1){6}}
\put(\etwo ,\gfiv){\line(1,0){4}}
\put(\etwo ,\gsev){\oval(4,4)[br]}
\put(\etwo ,\gsix){\oval(4,4)[br]}
\put(\ethr ,\gsev){\oval(4,4)[bl]}
\put(\ethr ,\gsix){\oval(4,4)[bl]}
}

\nc\vuh[3]{\uvariablearrowpos{#1}{#2}{#3}
\put(\eone ,\gone){\line(0,1){\fone}}
\put(\etwo ,\gfou){\line(0,1){6}}
\put(\ethr ,\gfou){\line(0,1){6}}
\put(\etwo ,\gfiv){\line(1,0){4}}
\put(\etwo ,\gsev){\oval(4,4)[br]}
\put(\ethr ,\gsev){\oval(4,4)[bl]}
}

\nc\vuhc[3]{\uvariablearrowpos{#1}{#2}{#3}
\put(\eone ,\gtwo){\line(0,1){\fthr}}
\put(\etwo ,\gfou){\line(0,1){6}}
\put(\ethr ,\gfou){\line(0,1){6}}
\put(\etwo ,\gfiv){\line(1,0){4}}
\put(\etwo ,\gsev){\oval(4,4)[br]}
\put(\ethr ,\gsev){\oval(4,4)[bl]}
\put(\eone ,\gone){\oval(4,4)[t]}
\put(\eone ,\gthr){\oval(4,4)[b]}
}

\nc\vuc[3]{\uvariablearrowpos{#1}{#2}{#3}
\put(\eone ,\gtwo){\line(0,1){\fthr}}
\put(\etwo ,\gsev){\oval(4,4)[br]}
\put(\ethr ,\gsev){\oval(4,4)[bl]}
\put(\eone ,\gone){\oval(4,4)[t]}
\put(\eone ,\gthr){\oval(4,4)[b]}
}

\nc\vrme[3]{\rvariablearrowpos{#1}{#2}{#3}
\put(\gtwo ,\eone){\line(1,0){\fthr}}
\put(\gtwo ,\ethr){\oval(4,4)[bl]}
\put(\gsev ,\ethr){\oval(4,4)[bl]}
\put(\gsix ,\ethr){\oval(4,4)[bl]}
\put(\gtwo ,\etwo){\oval(4,4)[tl]}
\put(\gsev ,\etwo){\oval(4,4)[tl]}
\put(\gsix ,\etwo){\oval(4,4)[tl]}
}

\nc\vrm[3]{\rvariablearrowpos{#1}{#2}{#3}
\put(\gtwo ,\eone){\line(1,0){\fthr}}
\put(\gtwo ,\ethr){\oval(4,4)[bl]}
\put(\gsev ,\ethr){\oval(4,4)[bl]}
\put(\gtwo ,\etwo){\oval(4,4)[tl]}
\put(\gsev ,\etwo){\oval(4,4)[tl]}
}

\nc\vre[3]{\rvariablearrowpos{#1}{#2}{#3}
\put(\gone ,\eone){\line(1,0){\fone}}
\put(\gsev ,\ethr){\oval(4,4)[bl]}
\put(\gsix ,\ethr){\oval(4,4)[bl]}
\put(\gsev ,\etwo){\oval(4,4)[tl]}
\put(\gsix ,\etwo){\oval(4,4)[tl]}
}

\nc\vra[3]{\rvariablearrowpos{#1}{#2}{#3}
\put(\gone ,\eone){\line(1,0){\fone}}
\put(\gsev ,\ethr){\oval(4,4)[bl]}
\put(\gsev ,\etwo){\oval(4,4)[tl]}
}

\nc\vrhe[3]{\rvariablearrowpos{#1}{#2}{#3}
\put(\gone ,\eone){\line(1,0){\fone}}
\put(\gfou ,\ethr){\line(1,0){6}}
\put(\gfou ,\etwo){\line(1,0){6}}
\put(\gfiv ,\etwo){\line(0,1){4}}
\put(\gsev ,\ethr){\oval(4,4)[bl]}
\put(\gsix ,\ethr){\oval(4,4)[bl]}
\put(\gsev ,\etwo){\oval(4,4)[tl]}
\put(\gsix ,\etwo){\oval(4,4)[tl]}
}

\nc\vrh[3]{\rvariablearrowpos{#1}{#2}{#3}
\put(\gone ,\eone){\line(1,0){\fone}}
\put(\gfou ,\ethr){\line(1,0){6}}
\put(\gfou ,\etwo){\line(1,0){6}}
\put(\gfiv ,\etwo){\line(0,1){4}}
\put(\gsev ,\ethr){\oval(4,4)[bl]}
\put(\gsev ,\etwo){\oval(4,4)[tl]}
}

\nc\vrhc[3]{\rvariablearrowpos{#1}{#2}{#3}
\put(\gtwo ,\eone){\line(1,0){\fthr}}
\put(\gfou ,\ethr){\line(1,0){6}}
\put(\gfou ,\etwo){\line(1,0){6}}
\put(\gfiv ,\etwo){\line(0,1){4}}
\put(\gsev ,\ethr){\oval(4,4)[bl]}
\put(\gsev ,\etwo){\oval(4,4)[tl]}
\put(\gone ,\eone){\oval(4,4)[r]}
\put(\gthr ,\eone){\oval(4,4)[l]}
}

\nc\vrc[3]{\rvariablearrowpos{#1}{#2}{#3}
\put(\gtwo ,\eone){\line(1,0){\fthr}}
\put(\gsev ,\ethr){\oval(4,4)[bl]}
\put(\gsev ,\etwo){\oval(4,4)[tl]}
\put(\gone ,\eone){\oval(4,4)[r]}
\put(\gthr ,\eone){\oval(4,4)[l]}
}

\nc\reqarr[3]{\eone=#1
\fone=#2
\gone=#3
\gtwo=\gone
\divide\gtwo by 2
\advance\eone by \gtwo
\Req{\eone}{\fone}{\gone}
}
\nc\veqarr[3]{\eone=#1
\fone=#2
\gone=#3
\gtwo=\gone
\divide\gtwo by 2
\advance\fone by \gtwo
\Veq{\eone}{\fone}{\gone}
}

\nc\vuba[3]{
\vua{#1}{#2}{#3}
\put(\eone,\gone){\makebox(0,0)[b]{$\bullet$}}
}

\nc\vrba[3]{
\vra{#1}{#2}{#3}
\put(\gone,\eone){\makebox(0,0)[l]{$\bullet$}}
}

\nc\setra[4]{
\if#4a\vra{#1}{#2}{#3}
\else\if#4m\vrm{#1}{#2}{#3}
\else\if#4e\vre{#1}{#2}{#3}
\else\if#4q\vrme{#1}{#2}{#3}
\else\if#4h\vrh{#1}{#2}{#3}
\else\if#4k\vrhe{#1}{#2}{#3}
\else\if#4H\vrhc{#1}{#2}{#3}
\else\if#4c\vrc{#1}{#2}{#3}
\fi\fi\fi\fi\fi\fi\fi\fi
\if#4M\vrM{#1}{#2}{#3}\fi
\if#4E\vrE{#1}{#2}{#3}\fi
\if#4I\vrI{#1}{#2}{#3}\fi
\if#4i\vri{#1}{#2}{#3}\fi
\if#4T\vraT{#1}{#2}{#3}\fi
\if#4t\vrat{#1}{#2}{#3}\fi
\if#4s\vrMs{#1}{#2}{#3}\fi
\if#4Q\vraMe{#1}{#2}{#3}\fi
\if#4K\vraIe{#1}{#2}{#3}\fi
\if#4=\reqarr{#1}{#2}{#3}\fi
\if#4b\vrba{#1}{#2}{#3}\fi
\if#4B\vrBa{#1}{#2}{#3}\fi
\if#4l\vla{#1}{#2}{#3}\fi
}

\nc\setua[4]{
\if#4a\vua{#1}{#2}{#3}
\else\if#4m\vum{#1}{#2}{#3}
\else\if#4e\vue{#1}{#2}{#3}
\else\if#4q\vume{#1}{#2}{#3}
\else\if#4h\vuh{#1}{#2}{#3}
\else\if#4k\vuhe{#1}{#2}{#3}
\else\if#4H\vuhc{#1}{#2}{#3}
\else\if#4c\vuc{#1}{#2}{#3}
\fi\fi\fi\fi\fi\fi\fi\fi
\if#4M\vuM{#1}{#2}{#3}\fi
\if#4E\vuE{#1}{#2}{#3}\fi
\if#4I\vuI{#1}{#2}{#3}\fi
\if#4i\vui{#1}{#2}{#3}\fi
\if#4T\vuaT{#1}{#2}{#3}\fi
\if#4t\vuat{#1}{#2}{#3}\fi
\if#4s\vuMs{#1}{#2}{#3}\fi
\if#4Q\vuaMe{#1}{#2}{#3}\fi
\if#4K\vuaIe{#1}{#2}{#3}\fi
\if#4=\veqarr{#1}{#2}{#3}\fi
\if#4b\vuba{#1}{#2}{#3}\fi
\if#4B\vuBa{#1}{#2}{#3}\fi
\if#4d\vda{#1}{#2}{#3}\fi
}

\nc\ct{{\cal T}}
\nc\cs{{\cal S}}
\nc\ca{{\cal A}}

\nc\barr{\begin{center}
$\begin{array}{ccccccccccccccccccccccccccccccccccccccccccccccccccccccccccccc}}
\nc\earr{\end{array}$
\end{center}
}

\nc\BC{\vspace*{1mm}\begin{center}}
\nc\EC{\end{center}\vspace*{1mm}}

\nc\posar[4]{\icoo{#1}\jcoo{#2}
\if#4t
\fone=\vartwo \advance \fone by \varfou
\advance \varthr by -4
\advance \varone by 2
\setra\varone\fone\varthr{#3}
\else\if#4n
\fone=\vartwo \advance \fone by \varfou
\advance \varthr by -4
\setra\varone\fone\varthr{#3}
\else\if#4l\advance \vartwo by 2
\advance\varfou by -4
\setua\varone\vartwo\varfou{#3}
\else\if#4w
\advance\varfou by -4
\setua\varone\vartwo\varfou{#3}
\else\if#4b\advance \varone by 2
\advance \varthr by -4
\setra\varone\vartwo\varthr{#3}
\else\if#4s\advance \varone by 4
\advance \varthr by -4
\setra\varone\vartwo\varthr{#3}
\else\if#4r\eone=\varone \advance \eone by \varthr
\fone=\vartwo \advance \fone by 2
\advance \varfou by -4
\setua\eone\fone\varfou{#3}
\else\if#4e\eone=\varone \advance \eone by \varthr
\fone=\vartwo \advance \fone by 4
\advance \varfou by -4
\setua\eone\fone\varfou{#3}
\fi\fi\fi\fi\fi\fi\fi\fi
}

\nc\postr[4]{\icoo{#1}\jcoo{#2}
\if#4r
\advance \varone by \varthr
\put(\varone,\vartwo){\makebox(0,\varfou){$#3$}}
\else\if#4l
\put(\varone,\vartwo){\makebox(0,\varfou){$#3$}}
\else\if#4b
\put(\varone,\vartwo){\makebox(\varthr,0){$#3$}}
\else\if#4t
\advance\vartwo by \varfou
\put(\varone,\vartwo){\makebox(\varthr,0){$#3$}}
\else\if#4e
\advance\varone by \varthr
\put(\varone,\vartwo){\makebox(0,0){$#3$}}
\else\if#4s
\put(\varone,\vartwo){\makebox(0,0){$#3$}}
\else\if#4w
\advance\vartwo by \varfou
\put(\varone,\vartwo){\makebox(0,0){$#3$}}
\else\if#4n
\advance\varone by \varthr
\advance\vartwo by \varfou
\put(\varone,\vartwo){\makebox(0,0){$#3$}}
\else
\put(\varone,\vartwo){\makebox(\varthr,\varfou){$#3$}}
\fi\fi\fi\fi\fi\fi\fi\fi
}

\nc\bvu[4]{ \icoo{#1} \jcoo{#2}
\hone =#3
\htwo =\varone \advance \htwo by \hone
\hthr =\varone \advance \hthr by \varthr
\advance \hthr by -\hone
\hfou =\vartwo \advance \hfou by 2
\hfiv =\varfou \advance \hfiv by -4
\setua\htwo\hfou\hfiv{#4}
\setua\hthr\hfou\hfiv{#4}
}

\nc\fbvu[5]{ \icoo{#1} \jcoo{#2}\fone=\vartwo\advance
\fone by \varfou
\jcoo{#3}\advance \vartwo by -\fone
\varfou=\vartwo \vartwo=\fone
\hone =#4
\htwo =\varone \advance \htwo by \hone
\hthr =\varone \advance \hthr by \varthr
\advance \hthr by -\hone
\hfou =\vartwo \advance \hfou by 2
\hfiv =\varfou \advance \hfiv by -4
\setua\htwo\hfou\hfiv{#5}
\setua\hthr\hfou\hfiv{#5}
}

\nc\fbu[4]{
\fbvu{#1}{#2}{#3}7{#4}
}

\nc\bu[3]{\bvu{#1}{#2}{7}{#3}}

\nc\bvr[4]{ \icoo{#1} \jcoo{#2}
\hone =#3
\htwo =\vartwo \advance \htwo by \hone
\hthr =\vartwo \advance \hthr by \varfou
\advance \hthr by -\hone
\hfou =\varone \advance \hfou by 2
\hfiv =\varthr \advance \hfiv by -4
\setra\hfou\htwo\hfiv{#4}
\setra\hfou\hthr\hfiv{#4}
}

\nc\fbvr[5]{\icoo{#2} \fone=\varone
\icoo{#1} \advance\varone by \varthr
\varthr=\fone \advance\varthr by -\varone \jcoo{#3}
\hone =#4
\htwo =\vartwo \advance \htwo by \hone
\hthr =\vartwo \advance \hthr by \varfou
\advance \hthr by -\hone
\hfou =\varone \advance \hfou by 2
\hfiv =\varthr \advance \hfiv by -4
\setra\hfou\htwo\hfiv{#5}
\setra\hfou\hthr\hfiv{#5}
}

\nc\br[3]{\bvr{#1}{#2}{7}{#3}}

\nc\fbr[4]{\fbvr{#1}{#2}{#3}{7}{#4}}

\nc\barb[2]{\hbox{\begin{picture}(21,21)
\setua25{16}{#1}
\setra52{16}{#2}
\end{picture}}}

\nc\bxs[7]{\icoo{#1}\jcoo{#2}
\put(\varone ,\vartwo){\framebox(\varthr ,\varfou){$\ct^{\barb{#3}{#4}}_{[{#5},{#6}]}{#7}$}}
}

\nc\grbxs[7]{\icoo{#1}\jcoo{#2}
\put(\varone ,\vartwo){\framebox(\varthr ,\varfou){$Gr^{\barb{#3}{#4}}_{[{#5},{#6}]}{#7}$}}
}

\nc\bx[6]{\icoo{#1}\jcoo{#2}
\put(\varone ,\vartwo){\framebox(\varthr ,\varfou){$\ct^{\barb{#3}{#4}}_{[{#5},{#6}]}$}}
}

\nc\grbx[6]{\icoo{#1}\jcoo{#2}
\put(\varone ,\vartwo){\framebox(\varthr ,\varfou){$Gr^{\barb{#3}{#4}}_{[{#5},{#6}]}$}}
}

\nc\fbx[5]{\icoo{#1}\jcoo{#2}
\put(\varone ,\vartwo){\framebox(\varthr ,\varfou){$#3^{\barb{#4}{#5}}$}}
}

\nc\fbxs[6]{\icoo{#1}\jcoo{#2}
\put(\varone ,\vartwo){\framebox(\varthr ,\varfou){$#3^{\barb{#4}{#5}}{#6}$}}
}

\nc\Fbx[3]{\icoo{#1}\jcoo{#2}
\put(\varone ,\vartwo){\framebox(\varthr ,\varfou){$#3$}}
}

\nc\hbx[8]{\icoo{#1}\jcoo{#2}
\put(\varone ,\vartwo){\framebox(\varthr ,\varfou){$\ct^{\barb{#3}{#4}}_{[{#5},{#6}]}{#7}$}}
\hone=\varone \advance \hone by 5
\htwo=\varone \advance \htwo by \varthr
\advance \htwo by -5
\hthr =\vartwo \advance \hthr by 5
\hfou =\vartwo \advance \hfou by \varfou
\advance \hfou by -5
\hfiv =\varthr \advance \hfiv by -10
\hsix =\varfou \advance \hsix by -10
\if#8r\put(\hone ,\hthr){\vector(1,0){\hfiv}}
\else\if#8l\put(\htwo ,\hthr){\vector(-1,0){\hfiv}}
\else\if#8u\put(\hone ,\hthr){\vector(0,1){\hsix}}
\else\if#8d\put(\hone ,\hfou){\vector(0,-1){\hsix}}
\fi\fi\fi\fi
}

\nc\Fhbx[4]{\icoo{#1}\jcoo{#2}
\put(\varone ,\vartwo){\framebox(\varthr ,\varfou){$#3$}}
\hone=\varone \advance \hone by 5
\htwo=\varone \advance \htwo by \varthr
\advance \htwo by -5
\hthr =\vartwo \advance \hthr by 5
\hfou =\vartwo \advance \hfou by \varfou
\advance \hfou by -5
\hfiv =\varthr \advance \hfiv by -10
\hsix =\varfou \advance \hsix by -10
\if#4r\put(\hone ,\hthr){\vector(1,0){\hfiv}}
\else\if#4l\put(\htwo ,\hthr){\vector(-1,0){\hfiv}}
\else\if#4u\put(\hone ,\hthr){\vector(0,1){\hsix}}
\else\if#4d\put(\hone ,\hfou){\vector(0,-1){\hsix}}
\fi\fi\fi\fi
}

\nc\fhbx[7]{\icoo{#1}\jcoo{#2}
\put(\varone ,\vartwo){\framebox(\varthr ,\varfou){${#3}^{\barb{#4}{#5}}{#6}$}}
\hone=\varone \advance \hone by 5
\htwo=\varone \advance \htwo by \varthr
\advance \htwo by -5
\hthr =\vartwo \advance \hthr by 5
\hfou =\vartwo \advance \hfou by \varfou
\advance \hfou by -5
\hfiv =\varthr \advance \hfiv by -10
\hsix =\varfou \advance \hsix by -10
\if#7r\put(\hone ,\hthr){\vector(1,0){\hfiv}}
\else\if#7l\put(\htwo ,\hthr){\vector(-1,0){\hfiv}}
\else\if#7u\put(\hone ,\hthr){\vector(0,1){\hsix}}
\else\if#7d\put(\hone ,\hfou){\vector(0,-1){\hsix}}
\fi\fi\fi\fi
}

\nc\grhbx[8]{\icoo{#1}\jcoo{#2}
\put(\varone ,\vartwo){\framebox(\varthr ,\varfou){$Gr^{\barb{#3}{#4}}_{[{#5},{#6}]}{#7}$}}
\hone=\varone \advance \hone by 5
\htwo=\varone \advance \htwo by \varthr
\advance \htwo by -5
\hthr =\vartwo \advance \hthr by 5
\hfou =\vartwo \advance \hfou by \varfou
\advance \hfou by -5
\hfiv =\varthr \advance \hfiv by -10
\hsix =\varfou \advance \hsix by -10
\if#8r\put(\hone ,\hthr){\vector(1,0){\hfiv}}
\else\if#8l\put(\htwo ,\hthr){\vector(-1,0){\hfiv}}
\else\if#8u\put(\hone ,\hthr){\vector(0,1){\hsix}}
\else\if#8d\put(\hone ,\hfou){\vector(0,-1){\hsix}}
\fi\fi\fi\fi
}

\nc\htpl[6]{\icoo{#1}\jcoo{#2}
\fone=\vartwo \advance \fone by \varfou
\eone=\varone \advance \eone by 10
\ftwo=\fone \advance \ftwo by -5
\fthr=\fone \advance \fthr by -15
\hthr=\vartwo \advance \hthr by 20
\hone=\varthr \multiply \hone by 2
\htwo=\varfou \multiply \htwo by 2
\advance \htwo by -50
\put(\eone ,\ftwo){\oval(\hone ,\htwo)[br]}
\put(\eone ,\hthr){\vector(-1,0){5}}
\put(\varone ,\fone){\line(1,0){\varthr}}
\advance\fone by -2
\put(\varone ,\fone){\makebox(\varthr,0)[t]{${#3}^{\arrb{#4}}_{#5}{#6}$}}
}

\nc\htpd[6]{\icoo{#1}\jcoo{#2}
\fone=\varone \advance \fone by \varthr
\eone=\vartwo \advance \eone by 10
\ftwo=\fone \advance \ftwo by -5
\fthr=\fone \advance \fthr by -15
\hthr=\varone \advance \hthr by 20
\hone=\varfou \multiply \hone by 2
\htwo=\varthr \multiply \htwo by 2
\advance \htwo by -50
\put(\ftwo ,\eone){\oval(\htwo ,\hone)[tl]}
\put(\hthr ,\eone){\vector(0,-1){5}}
\put(\fone ,\vartwo){\line(0,1){\varfou}}
\advance\fone by -2
\put(\fone ,\vartwo){\makebox(0 ,\varfou)[r]{${#3}^{\arub{#4}}_{#5}{#6}$}}
}

\nc\htpu[6]{\icoo{#1}\jcoo{#2}
\eone=\varone \advance \eone by \varthr
\advance \eone by -20
\etwo=\varone \advance \etwo by 5
\fone=\vartwo \advance \fone by \varfou
\advance \fone by -10
\hone=\varfou \multiply \hone by 2
\htwo=\varthr \multiply \htwo by 2
\advance \htwo by -50
\put(\etwo ,\fone){\oval(\htwo ,\hone)[br]}
\put(\eone ,\fone){\vector(0,1){5}}
\put(\varone ,\vartwo){\line(0,1){\varfou}}
\advance\varone by 2
\put(\varone ,\vartwo){\makebox(0,\varfou)[l]{${#3}^{\arub{#4}}_{#5}{#6}$}}
}

\nc\htpr[6]{\icoo{#1}\jcoo{#2}
\eone=\vartwo \advance \eone by \varfou
\advance \eone by -20
\etwo=\vartwo \advance \etwo by 5
\fone=\varone \advance \fone by \varthr
\advance \fone by -10
\hone=\varthr \multiply \hone by 2
\htwo=\varfou \multiply \htwo by 2
\advance \htwo by -50
\put(\fone ,\etwo){\oval(\hone ,\htwo)[tl]}
\put(\fone ,\eone){\vector(1,0){5}}
\put(\varone ,\vartwo){\line(1,0){\varthr}}
\advance\vartwo by 7
\put(\varone ,\vartwo){\makebox(\varthr ,0)[b]{${#3}^{\arrb{#4}}_{#5}{#6}$}}
}

\nc\htp[6]{
\if#6u\htpu{#1}{#2}{\ct}{#3}{[{#4},{#5}]}{}
\else\if#6d\htpd{#1}{#2}{\ct}{#3}{[{#4},{#5}]}{}
\else\if#6r\htpr{#1}{#2}{\ct}{#3}{[{#4},{#5}]}{}
\else\if#6l\htpl{#1}{#2}{\ct}{#3}{[{#4},{#5}]}{}
\fi\fi\fi\fi
}

\nc\grhtp[6]{
\if#6u\htpu{#1}{#2}{Gr}{#3}{[{#4},{#5}]}{}
\else\if#6d\htpd{#1}{#2}{Gr}{#3}{[{#4},{#5}]}{}
\else\if#6r\htpr{#1}{#2}{Gr}{#3}{[{#4},{#5}]}{}
\else\if#6l\htpl{#1}{#2}{Gr}{#3}{[{#4},{#5}]}{}
\fi\fi\fi\fi
}

\nc\coosq[4]{\icoo{#2}\eone=\varone
\advance \eone by \varthr
\icoo{#1}\advance \eone by -\varone
\varthr=\eone
\jcoo{#4}\fone=\vartwo
\advance \fone by \varfou
\jcoo{#3}\advance \fone by -\vartwo
\varfou=\fone
}

\nc\fhtpl[4]{\coosq{#1}{#2}{#3}{#4}
\fone=\vartwo \advance \fone by \varfou
\eone=\varone \advance \eone by 10
\ftwo=\fone \advance \ftwo by -5
\fthr=\fone \advance \fthr by -15
\hthr=\vartwo \advance \hthr by 20
\hone=\varthr \multiply \hone by 2
\htwo=\varfou \multiply \htwo by 2
\advance \htwo by -50
\put(\eone ,\ftwo){\oval(\hone ,\htwo)[br]}
\put(\eone ,\hthr){\vector(-1,0){5}}
}

\nc\fhtpd[4]{\coosq{#1}{#2}{#3}{#4}
\fone=\varone \advance \fone by \varthr
\eone=\vartwo \advance \eone by 10
\ftwo=\fone \advance \ftwo by -5
\fthr=\fone \advance \fthr by -15
\hthr=\varone \advance \hthr by 20
\hone=\varfou \multiply \hone by 2
\htwo=\varthr \multiply \htwo by 2
\advance \htwo by -50
\put(\ftwo ,\eone){\oval(\htwo ,\hone)[tl]}
\put(\hthr ,\eone){\vector(0,-1){5}}
}

\nc\fhtpu[4]{\coosq{#1}{#2}{#3}{#4}
\eone=\varone \advance \eone by \varthr
\advance \eone by -20
\etwo=\varone \advance \etwo by 5
\fone=\vartwo \advance \fone by \varfou
\advance \fone by -10
\hone=\varfou \multiply \hone by 2
\htwo=\varthr \multiply \htwo by 2
\advance \htwo by -50
\put(\etwo ,\fone){\oval(\htwo ,\hone)[br]}
\put(\eone ,\fone){\vector(0,1){5}}
}

\nc\fhtpr[4]{\coosq{#1}{#2}{#3}{#4}
\eone=\vartwo \advance \eone by \varfou
\advance \eone by -20
\etwo=\vartwo \advance \etwo by 5
\fone=\varone \advance \fone by \varthr
\advance \fone by -10
\hone=\varthr \multiply \hone by 2
\htwo=\varfou \multiply \htwo by 2
\advance \htwo by -50
\put(\fone ,\etwo){\oval(\hone ,\htwo)[tl]}
\put(\fone ,\eone){\vector(1,0){5}}
}

\nc\diff[5]{\coosq{#1}{#2}{#3}{#4}
\fone=\vartwo \advance \fone by \varfou
\advance \fone by 5
\eone=\varthr \divide \eone by 2
\advance \eone by \varone
\advance \eone by -10
\etwo=\varthr \advance \etwo by \varone
\advance \varthr by 20
\put(\eone,\fone){\oval(\varthr,30)[t]}
\put(\etwo,\fone){\line(0,-1){3}}
\advance \varone by -5
\multiply \varfou by 2
\advance \varfou by 10
\put(\varone,\fone){\oval(30,\varfou)[bl]}
\put(\varone,\vartwo){\vector(1,0){3}}
\advance \varone by -20
\advance \fone by 20
\put(\varone,\fone){\makebox(0,0){$#5$}}
}

\nc\xdiff[7]{\coosq{#1}{#2}{#3}{#4}
\fone=\vartwo \advance \fone by \varfou
\advance \fone by 5
\eone=\varthr \divide \eone by 2
\advance \eone by \varone
\advance \eone by -10
\etwo=\varthr \advance \etwo by \varone
\advance \varthr by 20
\put(\eone,\fone){\oval(\varthr,30)[t]}
\put(\etwo,\fone){\line(0,-1){#7}}
\advance \varone by -5
\multiply \varfou by 2
\advance \varfou by 10
\put(\varone,\fone){\oval(30,\varfou)[bl]}
\put(\varone,\vartwo){\vector(1,0){#6}}
\advance \varone by -20
\advance \fone by 20
\put(\varone,\fone){\makebox(0,0){$#5$}}
}

\nc\vdiff[7]{\icoo#2
\xone=\varone
\icoo#1
\advance\xone by -\varone
\advance\xone by 3
\jcoo#6
\yone=\vartwo
\advance\yone by \varfou
\jcoo#5
\advance\vartwo by \varfou
\advance\yone by -\vartwo
\advance\yone by 3
\xdiff{#1}{#3}{#4}{#6}{#7}{\xone}{\yone}
}

\nc\xbx[2]{\icoo{#1}\jcoo{#2}
\eone=\varone\advance \eone by -5
\etwo=\varthr\advance \etwo by 10
\fone=\vartwo\advance \fone by -5
\ftwo=\vartwo \advance \ftwo by \varfou
\advance \ftwo by 5
\put(\eone ,\fone){\line(1,1){\etwo}}
\put(\eone ,\ftwo){\line(1,-1){\etwo}}
\put(\varone,\vartwo){\framebox(\varthr,\varfou){}}
}

\nc\trui[3]{\icoo{#1}\jcoo{#2}
\eone=\varone \advance \eone by \varthr
\fone=\varthr \divide \fone by 3
\etwo=\eone \advance \etwo by -\fone
\advance  \fone by \vartwo
\put(\varone,\vartwo){\line(1,1){\varthr}}
\put(\varone,\vartwo){\line(1,0){\varthr}}
\put(\eone,\vartwo){\line(0,1){\varthr}}
\put(\etwo,\fone){\makebox(0,0){$#3$}}
}

\nc\xtrui[2]{\icoo{#1}\jcoo{#2}
\eone=\varone \advance \eone by \varthr
\fone=\varthr \divide \fone by 3
\etwo=\eone \advance \etwo by -\fone
\advance  \fone by \vartwo
\hone =\varthr \divide \hone by 2
\htwo =\hone \advance \htwo by \varone
\advance \htwo by -5
\hthr =\vartwo \advance \hthr by -5
\put(\varone,\vartwo){\line(1,1){\varthr}}
\put(\varone,\vartwo){\line(1,0){\varthr}}
\put(\eone,\vartwo){\line(0,1){\varthr}}
\advance \hone by 10
\put(\htwo,\hthr){\line(1,1){\hone}}
\advance \htwo by \hone
\put(\htwo,\hthr){\line(-1,1){\hone}}
}

\nc\truj[3]{\icoo{#1}\jcoo{#2}
\eone=\varone \advance \eone by \varthr
\ftwo=\vartwo \advance \ftwo by \varfou
\fone=\varfou \divide \fone by 3
\etwo=\eone \advance \etwo by -\fone
\advance \fone by \vartwo
\put(\eone,\ftwo){\line(-1,-1){\varfou}}
\put(\eone,\ftwo){\line(0,-1){\varfou}}
\put(\eone,\vartwo){\line(-1,0){\varfou}}
\put(\etwo,\fone){\makebox(0,0){$#3$}}
}

\nc\xtruj[2]{\icoo{#1}\jcoo{#2}
\eone=\varone \advance \eone by \varthr
\ftwo=\vartwo \advance \ftwo by \varfou
\fone=\varfou \divide \fone by 3
\etwo=\eone \advance \etwo by -\fone
\advance \fone by \vartwo
\put(\eone,\ftwo){\line(-1,-1){\varfou}}
\put(\eone,\ftwo){\line(0,-1){\varfou}}
\put(\eone,\vartwo){\line(-1,0){\varfou}}
\hone=\varfou\divide \hone by  2
\advance \ftwo by -\hone
\advance \hone by 10
\advance \ftwo by 5
\advance \eone by 5
\put(\eone,\ftwo){\line(-1,-1){\hone}}
\advance\vartwo by -5
\put(\eone,\vartwo){\line(-1,1){\hone}}
}

\nc\trdi[3]{\icoo{#1}\jcoo{#2}
\eone=\vartwo \advance \eone by \varfou
\fone=\varfou \divide \fone by 3
\etwo=\eone \advance \etwo by -\fone
\advance  \fone by \varone
\put(\varone,\vartwo){\line(1,1){\varfou}}
\put(\varone,\vartwo){\line(0,1){\varfou}}
\put(\varone,\eone){\line(1,0){\varfou}}
\put(\fone,\etwo){\makebox(0,0){$#3$}}
}

\nc\xtrdi[2]{\icoo{#1}\jcoo{#2}
\eone=\vartwo \advance \eone by \varfou
\fone=\varfou \divide \fone by 3
\etwo=\eone \advance \etwo by -\fone
\advance  \fone by \varone
\hone =\varfou \divide \hone by 2
\htwo =\hone \advance \htwo by \vartwo
\advance \htwo by -5
\hthr =\varone \advance \hthr by -5
\put(\varone,\vartwo){\line(1,1){\varfou}}
\put(\varone,\vartwo){\line(0,1){\varfou}}
\put(\varone,\eone){\line(1,0){\varfou}}
\advance \hone by 10
\put(\hthr,\htwo){\line(1,1){\hone}}
\advance \htwo by \hone
\put(\hthr,\htwo){\line(1,-1){\hone}}
}

\nc\trdj[3]{\icoo{#1}\jcoo{#2}
\fone=\vartwo \advance \fone by \varfou
\put(\varone,\vartwo){\line(1,1){\varfou}}
\put(\varone,\fone){\line(1,0){\varfou}}
\put(\varone,\vartwo){\line(0,1){\varfou}}
\yone=\varfou \divide \yone by 3
\advance \varone by \yone
\advance \fone by -\yone
\put(\varone,\fone){\makebox(0,0){$#3$}}
}

\nc\xtrdj[2]{\icoo{#1}\jcoo{#2}
\eone=\vartwo \advance \eone by \varfou
\ftwo=\varone \advance \ftwo by \varthr
\fone=\varfou \divide \fone by 3
\etwo=\eone \advance \etwo by -\fone
\advance \fone by \varone
\put(\ftwo,\eone){\line(-1,-1){\varfou}}
\put(\ftwo,\eone){\line(-1,0){\varfou}}
\put(\varone,\eone){\line(0,-1){\varfou}}
\hone=\varfou\divide \hone by  2
\advance \ftwo by -\hone
\advance \hone by 10
\advance \ftwo by 5
\advance \eone by 5
\put(\ftwo,\eone){\line(-1,-1){\hone}}
\advance\varone by -5
\put(\varone,\eone){\line(1,-1){\hone}}
}

\nc\Fli[4]{\icoo{#1}\jcoo{#2}
\if#4r\advance \varone by \varthr
\put(\varone,\vartwo){\line(0,1){\varfou}}
\advance\varone by -2
\put(\varone,\vartwo){\makebox(0,\varfou)[r]{$#3$}}
\else\if#4l
\put(\varone,\vartwo){\line(0,1){\varfou}}
\advance \varone by 2
\put(\varone,\vartwo){\makebox(0,\varfou)[l]{$#3$}}
\else\if#4b
\put(\varone,\vartwo){\line(1,0){\varthr}}
\advance \vartwo by 2
\put(\varone,\vartwo){\makebox(\varthr,20)[b]{$#3$}}
\else\if#4t
\advance \vartwo by \varfou
\put(\varone,\vartwo){\line(1,0){\varthr}}
\advance \vartwo by -2
\put(\varone,\vartwo){\makebox(\varthr,0)[t]{$#3$}}
\fi\fi\fi\fi
}

\nc\arub[1]{\hbox{\begin{picture}(4,21)
\setua25{16}{#1}
\end{picture}}}

\nc\arrb[1]{\hbox{\begin{picture}(21,4)
\setra22{16}{#1}
\end{picture}}}

\nc\fli[6]{\icoo{#1}\jcoo{#2}
\if#6r\advance \varone by \varthr
\put(\varone,\vartwo){\line(0,1){\varfou}}
\advance\varone by -2
\put(\varone,\vartwo){\makebox(0,\varfou)[r]{${#3}^{\arub{#4}}{#5}$}}
\else\if#6l
\put(\varone,\vartwo){\line(0,1){\varfou}}
\advance \varone by 2
\put(\varone,\vartwo){\makebox(0,\varfou)[l]{${#3}^{\arub{#4}}{#5}$}}
\else\if#6b
\put(\varone,\vartwo){\line(1,0){\varthr}}
\advance \vartwo by 4
\put(\varone,\vartwo){\makebox(\varthr,0)[b]{${#3}^{\arrb{#4}}{#5}$}}
\else\if#6t
\advance \vartwo by \varfou
\put(\varone,\vartwo){\line(1,0){\varthr}}
\advance \vartwo by -2
\put(\varone,\vartwo){\makebox(\varthr,0)[t]{${#3}^{\arrb{#4}}{#5}$}}
\fi\fi\fi\fi
}

\nc\li[7]{\icoo{#1}\jcoo{#2}
\if#7r\advance \varone by \varthr
\put(\varone,\vartwo){\line(0,1){\varfou}}
\advance\varone by -2
\put(\varone,\vartwo){\makebox(0,\varfou)[r]{${\ct}^{\arub{#3}}_{[#4,#5]}{#6}$}}
\else\if#7l
\put(\varone,\vartwo){\line(0,1){\varfou}}
\advance \varone by 2
\put(\varone,\vartwo){\makebox(0,\varfou)[l]{${\ct}^{\arub{#3}}_{[#4,#5]}{#6}$}}
\else\if#7b
\put(\varone,\vartwo){\line(1,0){\varthr}}
\advance \vartwo by 6
\put(\varone,\vartwo){\makebox(\varthr,0)[b]{${\ct}^{\arrb{#3}}_{[#4,#5]}{#6}$}}
\else\if#7t
\advance \vartwo by \varfou
\put(\varone,\vartwo){\line(1,0){\varthr}}
\advance \vartwo by -2
\put(\varone,\vartwo){\makebox(\varthr,0)[t]{${\ct}^{\arrb{#3}}_{[#4,#5]}{#6}$}}
\fi\fi\fi\fi
}

\nc\grli[7]{\icoo{#1}\jcoo{#2}
\if#7r\advance \varone by \varthr
\put(\varone,\vartwo){\line(0,1){\varfou}}
\advance\varone by -2
\put(\varone,\vartwo){\makebox(0,\varfou)[r]{${Gr}^{\arub{#3}}_{[#4,#5]}{#6}$}}
\else\if#7l
\put(\varone,\vartwo){\line(0,1){\varfou}}
\advance \varone by 2
\put(\varone,\vartwo){\makebox(0,\varfou)[l]{${Gr}^{\arub{#3}}_{[#4,#5]}{#6}$}}
\else\if#7b
\put(\varone,\vartwo){\line(1,0){\varthr}}
\advance \vartwo by 6
\put(\varone,\vartwo){\makebox(\varthr,0)[b]{${Gr}^{\arrb{#3}}_{[#4,#5]}{#6}$}}
\else\if#7t
\advance \vartwo by \varfou
\put(\varone,\vartwo){\line(1,0){\varthr}}
\advance \vartwo by -2
\put(\varone,\vartwo){\makebox(\varthr,0)[t]{${Gr}^{\arrb{#3}}_{[#4,#5]}{#6}$}}
\fi\fi\fi\fi
}

\nc\xli[3]{\icoo{#1}\jcoo{#2}
\if#3r\advance \varone by \varthr
\put(\varone,\vartwo){\line(0,1){\varfou}}
\yone=\varfou \divide \yone by 2
\advance \yone by \vartwo
\advance \yone by -20
\advance \varone by -5
\put(\varone,\yone){\line(1,4){10}}
\advance \varone by 10
\put(\varone,\yone){\line(-1,4){10}}
\else\if#3l
\put(\varone,\vartwo){\line(0,1){\varfou}}
\yone=\varfou \divide \yone by 2
\advance \yone by \vartwo
\advance \yone by -20
\advance \varone by -5
\put(\varone,\yone){\line(1,4){10}}
\advance \varone by 10
\put(\varone,\yone){\line(-1,4){10}}
\else\if#3b
\put(\varone,\vartwo){\line(1,0){\varthr}}
\xone=\varthr \divide \xone by 2
\advance \xone by \varone
\advance \xone by -20
\advance \vartwo by -5
\put(\xone,\vartwo){\line(4,1){40}}
\advance \vartwo by 10
\put(\xone,\vartwo){\line(4,-1){40}}
\else\if#3t
\advance \vartwo by \varfou
\put(\varone,\vartwo){\line(1,0){\varthr}}
\xone=\varthr \divide \xone by 2
\advance \xone by \varone
\advance \xone by -20
\advance \vartwo by -5
\put(\xone,\vartwo){\line(4,1){40}}
\advance \vartwo by 10
\put(\xone,\vartwo){\line(4,-1){40}}
\fi\fi\fi\fi
}

\nc\warn[3]{\icoo{#1} \jcoo{#2}
\xone=\varthr \divide \xone by 2
\yone=\varfou \divide \yone by 2
\advance \varone by \xone
\advance \vartwo by \yone
\if#3r \advance \varone by \xone
\else\if#3l \advance \varone by -\xone
\else\if#3t \advance \vartwo by \yone
\else\if#3b \advance \vartwo by -\yone
\fi\fi\fi\fi
\put(\varone,\vartwo){\circle{3}}
\put(\varone,\vartwo){\circle{7}}
\put(\varone,\vartwo){\circle{11}}
\put(\varone,\vartwo){\circle{15}}
\put(\varone,\vartwo){\circle{19}}
}

\nc\uSbx[9]{\icoo{#1}\jcoo{#2}
\eone=\varthr \divide \eone by 2
\xthr=\varone \advance \xthr by \eone
\advance \eone by -35
\divide \eone by 2
\fone=\varfou \divide \fone by 2
\ythr=\vartwo \advance \ythr by \fone
\advance \fone by -35
\divide \fone by 2
\xone=\varone \advance \xone by \eone
\xtwo=\xone \advance \xtwo by \eone
\xsix=\varone \advance \xsix by \varthr 
\advance \xsix by -\eone
\xfiv =\xsix \advance \xfiv by -\eone
\xfou =\xfiv \advance \xfou by -20
\yone=\vartwo \advance \yone by \fone
\ytwo=\yone \advance \ytwo by \fone
\ysix=\vartwo \advance \ysix by \varfou 
\advance \ysix by -\fone
\yfiv =\ysix \advance \yfiv by -\fone
\yfou =\yfiv \advance \yfou by -20
\advance\ythr by 3
\setra{\xtwo}{\ysix}{20}{#4}
\setra{\xfou}{\ysix}{20}{#4}
\setua{\xsix}{\ytwo}{20}{#3}
\setua{\xsix}{\yfou}{20}{#3}
\put(\xsix,\ythr){\makebox(0,0){$\vdots$}}
\put(\xthr,\ysix){\makebox(0,0){$\cdots$}}
\let\J=#7 \let\I=#8
\put(\xsix,\yone){\makebox(0,0){${#5}$}}
\let\J=#6 \let\I=#9
\put(\xone,\ysix){\makebox(0,0){${#5}$}}
\let\J=#7 \let\I=#9
\put(\xsix,\ysix){\makebox(0,0){${#5}$}}
}

\nc\Sbx[9]{\icoo{#1}\jcoo{#2}
\eone=\varthr \divide \eone by 2
\xthr=\varone \advance \xthr by \eone
\advance \eone by -35
\divide \eone by 2
\fone=\varfou \divide \fone by 2
\ythr=\vartwo \advance \ythr by \fone
\advance \fone by -35
\divide \fone by 2
\xone=\varone \advance \xone by \eone
\xtwo=\xone \advance \xtwo by \eone
\xsix=\varone \advance \xsix by \varthr 
\advance \xsix by -\eone
\xfiv =\xsix \advance \xfiv by -\eone
\xfou =\xfiv \advance \xfou by -20
\yone=\vartwo \advance \yone by \fone
\ytwo=\yone \advance \ytwo by \fone
\ysix=\vartwo \advance \ysix by \varfou 
\advance \ysix by -\fone
\yfiv =\ysix \advance \yfiv by -\fone
\yfou =\yfiv \advance \yfou by -20
\advance\ythr by 3
\setra{\xtwo}{\yone}{20}{#4}
\setra{\xtwo}{\ysix}{20}{#4}
\setra{\xfou}{\yone}{20}{#4}
\setra{\xfou}{\ysix}{20}{#4}
\setua{\xone}{\ytwo}{20}{#3}
\setua{\xsix}{\ytwo}{20}{#3}
\setua{\xone}{\yfou}{20}{#3}
\setua{\xsix}{\yfou}{20}{#3}
\put(\xone,\ythr){\makebox(0,0){$\vdots$}}
\put(\xsix,\ythr){\makebox(0,0){$\vdots$}}
\put(\xthr,\yone){\makebox(0,0){$\cdots$}}
\put(\xthr,\ysix){\makebox(0,0){$\cdots$}}
\let\J=#6 \let\I=#8
\put(\xone,\yone){\makebox(0,0){${#5}$}}
\let\J=#7 \let\I=#8
\put(\xsix,\yone){\makebox(0,0){${#5}$}}
\let\J=#6 \let\I=#9
\put(\xone,\ysix){\makebox(0,0){${#5}$}}
\let\J=#7 \let\I=#9
\put(\xsix,\ysix){\makebox(0,0){${#5}$}}
}

\nc\Sbxsmall[9]{\icoo{#1}\jcoo{#2}
\eone=\varthr \divide \eone by 2
\xthr=\varone \advance \xthr by \eone
\advance \eone by -35
\divide \eone by 2
\fone=\varfou \divide \fone by 2
\ythr=\vartwo \advance \ythr by \fone
\advance \fone by -35
\divide \fone by 2
\xone=\varone \advance \xone by \eone
\xtwo=\xone \advance \xtwo by \eone
\xsix=\varone \advance \xsix by \varthr 
\advance \xsix by -\eone
\xfiv =\xsix \advance \xfiv by -\eone
\xfou =\xfiv \advance \xfou by -20
\yone=\vartwo \advance \yone by \fone
\ytwo=\yone \advance \ytwo by \fone
\ysix=\vartwo \advance \ysix by \varfou 
\advance \ysix by -\fone
\yfiv =\ysix \advance \yfiv by -\fone
\yfou =\yfiv \advance \yfou by -20
\advance\ythr by 3
\advance\xtwo by 2
\setra{\xtwo}{\yone}{18}{#4}
\setra{\xtwo}{\ysix}{18}{#4}
\setra{\xfou}{\yone}{18}{#4}
\setra{\xfou}{\ysix}{18}{#4}
\setua{\xone}{\ytwo}{20}{#3}
\setua{\xsix}{\ytwo}{20}{#3}
\setua{\xone}{\yfou}{20}{#3}
\setua{\xsix}{\yfou}{20}{#3}
\put(\xone,\ythr){\makebox(0,0){$\vdots$}}
\put(\xsix,\ythr){\makebox(0,0){$\vdots$}}
\put(\xthr,\yone){\makebox(0,0){$\cdots$}}
\put(\xthr,\ysix){\makebox(0,0){$\cdots$}}
\let\J=#6 \let\I=#8
\put(\xone,\yone){\makebox(0,0){${#5}$}}
\let\J=#7 \let\I=#8
\put(\xsix,\yone){\makebox(0,0){${#5}$}}
\let\J=#6 \let\I=#9
\put(\xone,\ysix){\makebox(0,0){${#5}$}}
\let\J=#7 \let\I=#9
\put(\xsix,\ysix){\makebox(0,0){${#5}$}}
}

\nc\tinySbx[9]{\icoo{#1}\jcoo{#2}
\eone=\varthr \divide \eone by 2
\xthr=\varone \advance \xthr by \eone
\advance \eone by -25
\divide \eone by 2
\fone=\varfou \divide \fone by 2
\ythr=\vartwo \advance \ythr by \fone
\advance \fone by -35
\divide \fone by 2
\xone=\varone \advance \xone by \eone
\xtwo=\xone \advance \xtwo by \eone
\xsix=\varone \advance \xsix by \varthr 
\advance \xsix by -\eone
\xfiv =\xsix \advance \xfiv by -\eone
\xfou =\xfiv \advance \xfou by -15
\yone=\vartwo \advance \yone by \fone
\ytwo=\yone \advance \ytwo by \fone
\ysix=\vartwo \advance \ysix by \varfou 
\advance \ysix by -\fone
\yfiv =\ysix \advance \yfiv by -\fone
\yfou =\yfiv \advance \yfou by -20
\advance\ythr by 3
\setra{\xtwo}{\yone}{15}{#4}
\setra{\xtwo}{\ysix}{15}{#4}
\setra{\xfou}{\yone}{15}{#4}
\setra{\xfou}{\ysix}{15}{#4}
\setua{\xone}{\ytwo}{20}{#3}
\setua{\xsix}{\ytwo}{20}{#3}
\setua{\xone}{\yfou}{20}{#3}
\setua{\xsix}{\yfou}{20}{#3}
\put(\xone,\ythr){\makebox(0,0){$\vdots$}}
\put(\xsix,\ythr){\makebox(0,0){$\vdots$}}
\put(\xthr,\yone){\makebox(0,0){$\cdots$}}
\put(\xthr,\ysix){\makebox(0,0){$\cdots$}}
\let\J=#6 \let\I=#8
\put(\xone,\yone){\makebox(0,0){${#5}$}}
\let\J=#7 \let\I=#8
\put(\xsix,\yone){\makebox(0,0){${#5}$}}
\let\J=#6 \let\I=#9
\put(\xone,\ysix){\makebox(0,0){${#5}$}}
\let\J=#7 \let\I=#9
\put(\xsix,\ysix){\makebox(0,0){${#5}$}}
}

\nc\Adecone[2]{\let\aone=#1
\let\bone=#2
}

\nc\Adectwo[2]{\let\atwo=#1
\let\btwo=#2
}

\nc\Adecthr[2]{\let\athr=#1
\let\bthr=#2
}

\nc\Adecfou[2]{\let\afou=#1
\let\bfou=#2
}

\nc\Adecfiv[2]{\let\afiv=#1
\let\bfiv=#2
}

\nc\Adecsix[2]{\let\asix=#1
\let\bsix=#2
}

\nc\Adecsev[2]{\let\asev=#1
\let\bsev=#2
}

\nc\Adeceig[2]{\let\aeig=#1
\let\beig=#2
}

\nc\Adecnin[2]{\let\anin=#1
\let\bnin=#2
}

\nc\Bdecone[2]{\let\cone=#1
\let\done=#2
}

\nc\Bdectwo[2]{\let\ctwo=#1
\let\dtwo=#2
}

\nc\Bdecthr[2]{\let\cthr=#1
\let\dthr=#2
}

\nc\Bdecfou[2]{\let\cfou=#1
\let\dfou=#2
}

\nc\Bdecfiv[2]{\let\cfiv=#1
\let\dfiv=#2
}

\nc\Bdecsix[2]{\let\csix=#1
\let\dsix=#2
}

\nc\Bdecsev[2]{\let\csev=#1
\let\dsev=#2
}

\nc\Bdeceig[2]{\let\ceig=#1
\let\deig=#2
}

\nc\Bdecnin[2]{\let\cnin=#1
\let\dnin=#2
}

\nc\A[3]{
\if#11\Adecone{#2}{#3}\fi
\if#12\Adectwo{#2}{#3}\fi
\if#13\Adecthr{#2}{#3}\fi
\if#14\Adecfou{#2}{#3}\fi
\if#15\Adecfiv{#2}{#3}\fi
\if#16\Adecsix{#2}{#3}\fi
\if#17\Adecsev{#2}{#3}\fi
\if#18\Adeceig{#2}{#3}\fi
\if#19\Adecnin{#2}{#3}\fi
}

\nc\B[3]{
\if#11\Bdecone{#2}{#3}\fi
\if#12\Bdectwo{#2}{#3}\fi
\if#13\Bdecthr{#2}{#3}\fi
\if#14\Bdecfou{#2}{#3}\fi
\if#15\Bdecfiv{#2}{#3}\fi
\if#16\Bdecsix{#2}{#3}\fi
\if#17\Bdecsev{#2}{#3}\fi
\if#18\Bdeceig{#2}{#3}\fi
\if#19\Bdecnin{#2}{#3}\fi
}

\nc\Ipar[1]{\if#11 \let\varfiv=\aone
\let\varsix=\bone
\else\if#12 \let\varfiv=\atwo
\let\varsix=\btwo
\else\if#13 \let\varfiv=\athr
\let\varsix=\bthr
\else\if#14 \let\varfiv=\afou
\let\varsix=\bfou
\else\if#15 \let\varfiv=\afiv
\let\varsix=\bfiv
\else\if#16 \let\varfiv=\asix
\let\varsix=\bsix
\else\if#17 \let\varfiv=\asev
\let\varsix=\bsev
\else\if#18 \let\varfiv=\aeig
\let\varsix=\beig
\else\if#19 \let\varfiv=\anin
\let\varsix=\bnin
\fi\fi\fi\fi\fi\fi\fi\fi\fi
}

\nc\Jpar[1]{\if#11 \let\varsev=\cone
\let\vareig=\done
\else\if#12 \let\varsev=\ctwo
\let\vareig=\dtwo
\else\if#13 \let\varsev=\cthr
\let\vareig=\dthr
\else\if#14 \let\varsev=\cfou
\let\vareig=\dfou
\else\if#15 \let\varsev=\cfiv
\let\vareig=\dfiv
\else\if#16 \let\varsev=\csix
\let\vareig=\dsix
\else\if#17 \let\varsev=\csev
\let\vareig=\dsev
\else\if#18 \let\varsev=\ceig
\let\vareig=\deig
\else\if#19 \let\varsev=\cnin
\let\vareig=\dnin
\fi\fi\fi\fi\fi\fi\fi\fi\fi
}

\nc\Sbxp[5]{\Ipar{#1}\Jpar{#2}
\Sbx{#1}{#2}{#3}{#4}{#5}{\varfiv}{\varsix}{\varsev}{\vareig}
}

\nc\tinySbxp[5]{\Ipar{#1}\Jpar{#2}
\tinySbx{#1}{#2}{#3}{#4}{#5}{\varfiv}{\varsix}{\varsev}{\vareig}
}

\nc\Sbxpsmall[5]{\Ipar{#1}\Jpar{#2}
\Sbxsmall{#1}{#2}{#3}{#4}{#5}{\varfiv}{\varsix}{\varsev}{\vareig}
}

\nc\uSbxp[5]{\Ipar{#1}\Jpar{#2}
\uSbx{#1}{#2}{#3}{#4}{#5}{\varfiv}{\varsix}{\varsev}{\vareig}
}

\nc\zeroboxes[5]{
\put(\xone,\yone){\makebox(0,0){0}}
\let\J=#2 \let\I=#3
\put(\xsix,\yone){\makebox(0,0){${#5}$}}
\put(\xsix,\ysix){\makebox(0,0){0}}
}

\nc\fullboxes[5]{
\let\J=#1 \let\I=#3
\put(\xone,\yone){\makebox(0,0){${#5}$}}
\let\J=#2 \let\I=#3
\put(\xsix,\yone){\makebox(0,0){${#5}$}}
\let\J=#2 \let\I=#4
\put(\xsix,\ysix){\makebox(0,0){${#5}$}}
}

\nc\quesboxes[5]{
\let\J=#1 \let\I=#3
\put(\xone,\yone){\makebox(0,0){${#5}$}}
\put(\xsix,\yone){\makebox(0,0){$?$}}
\let\J=#2 \let\I=#4
\put(\xsix,\ysix){\makebox(0,0){${#5}$}}
}

\nc\formalStru[4]{\icoo{#1}\jcoo{#2}
\eone=\varthr \divide \eone by 2
\xthr=\varone \advance \xthr by \eone
\advance \eone by -35
\divide \eone by 2
\fone=\varfou \divide \fone by 2
\ythr=\vartwo \advance \ythr by \fone
\advance \fone by -35
\divide \fone by 2
\xone=\varone \advance \xone by \eone
\xtwo=\xone \advance \xtwo by \eone
\xsix=\varone \advance \xsix by \varthr 
\advance \xsix by -\eone
\xfiv =\xsix \advance \xfiv by -\eone
\xfou =\xfiv \advance \xfou by -20
\yone=\vartwo \advance \yone by \fone
\ytwo=\yone \advance \ytwo by \fone
\ysix=\vartwo \advance \ysix by \varfou 
\advance \ysix by -\fone
\yfiv =\ysix \advance \yfiv by -\fone
\yfou =\yfiv \advance \yfou by -20
\setra{\xtwo}{\yone}{20}{#4}
\setra{\xfou}{\yone}{20}{#4}
\setua{\xsix}{\ytwo}{20}{#3}
\setua{\xsix}{\yfou}{20}{#3}
\put(\xsix,\ythr){\makebox(0,0){$\vdots$}}
\put(\xthr,\yone){\makebox(0,0){$\cdots$}}
\advance \xtwo by 20
\advance \xfiv by -20
\advance \ytwo by 20
\advance \yfiv by -20
\put(\xthr,\ythr){\makebox(0,0){$\cdot$}}
\put(\xfiv,\yfiv){\makebox(0,0){$\cdot$}}
\put(\xtwo,\ytwo){\makebox(0,0){$\cdot$}}
}

\nc\formaltinyStru[4]{
\icoo{#1}\jcoo{#2}
\eone=\varthr \divide \eone by 2
\xthr=\varone \advance \xthr by \eone
\advance \eone by -25
\divide \eone by 2
\fone=\varfou \divide \fone by 2
\ythr=\vartwo \advance \ythr by \fone
\advance \fone by -35
\divide \fone by 2
\xone=\varone \advance \xone by \eone
\xtwo=\xone \advance \xtwo by \eone
\xsix=\varone \advance \xsix by \varthr 
\advance \xsix by -\eone
\xfiv =\xsix \advance \xfiv by -\eone
\xfou =\xfiv \advance \xfou by -15
\yone=\vartwo \advance \yone by \fone
\ytwo=\yone \advance \ytwo by \fone
\ysix=\vartwo \advance \ysix by \varfou 
\advance \ysix by -\fone
\yfiv =\ysix \advance \yfiv by -\fone
\yfou =\yfiv \advance \yfou by -20
\advance\ythr by 3
\setra{\xtwo}{\yone}{15}{#4}
\setra{\xfou}{\yone}{15}{#4}
\setua{\xsix}{\ytwo}{20}{#3}
\setua{\xsix}{\yfou}{20}{#3}
\put(\xsix,\ythr){\makebox(0,0){$\vdots$}}
\put(\xthr,\yone){\makebox(0,0){$\cdots$}}
\advance \xtwo by 20
\advance \xfiv by -20
\advance \ytwo by 20
\advance \yfiv by -20
\put(\xthr,\ythr){\makebox(0,0){$\cdot$}}
\advance\xthr by -20
\advance\ythr by -20
\put(\xthr,\ythr){\makebox(0,0){$\cdot$}}
\advance\xthr by 40
\advance\ythr by 40
\put(\xthr,\ythr){\makebox(0,0){$\cdot$}}
}

\nc\ftinyStru[9]{
\formaltinyStru{#1}{#2}{#3}{#4}
\fullboxes{#6}{#7}{#8}{#9}{#5}
}

\nc\ftinyStrup[5]{\Ipar{#1}\Jpar{#2}
\ftinyStru{#1}{#2}{#3}{#4}{#5}{\varfiv}{\varsix}{\varsev}{\vareig}
}

\nc\Stru[9]{
\formalStru{#1}{#2}{#3}{#4}
\zeroboxes{#6}{#7}{#8}{#9}{#5}
}

\nc\fStru[9]{
\formalStru{#1}{#2}{#3}{#4}
\fullboxes{#6}{#7}{#8}{#9}{#5}
}

\nc\gStru[9]{
\formalStru{#1}{#2}{#3}{#4}
\quesboxes{#6}{#7}{#8}{#9}{#5}
}

\nc\Strup[5]{\Ipar{#1}\Jpar{#2}
\Stru{#1}{#2}{#3}{#4}{#5}{\varfiv}{\varsix}{\varsev}{\vareig}
}

\nc\fStrup[5]{\Ipar{#1}\Jpar{#2}
\fStru{#1}{#2}{#3}{#4}{#5}{\varfiv}{\varsix}{\varsev}{\vareig}
}

\nc\gStrup[5]{\Ipar{#1}\Jpar{#2}
\gStru{#1}{#2}{#3}{#4}{#5}{\varfiv}{\varsix}{\varsev}{\vareig}
}

\nc\Strd[9]{\icoo{#1}\jcoo{#2}
\eone=\varthr \divide \eone by 2
\xthr=\varone \advance \xthr by \eone
\advance \eone by -35
\divide \eone by 2
\fone=\varfou \divide \fone by 2
\ythr=\vartwo \advance \ythr by \fone
\advance \fone by -35
\divide \fone by 2
\xone=\varone \advance \xone by \eone
\xtwo=\xone \advance \xtwo by \eone
\xsix=\varone \advance \xsix by \varthr 
\advance \xsix by -\eone
\xfiv =\xsix \advance \xfiv by -\eone
\xfou =\xfiv \advance \xfou by -20
\yone=\vartwo \advance \yone by \fone
\ytwo=\yone \advance \ytwo by \fone
\ysix=\vartwo \advance \ysix by \varfou 
\advance \ysix by -\fone
\yfiv =\ysix \advance \yfiv by -\fone
\yfou =\yfiv \advance \yfou by -20
\setra{\xtwo}{\ysix}{20}{#4}
\setra{\xfou}{\ysix}{20}{#4}
\setua{\xone}{\ytwo}{20}{#3}
\setua{\xone}{\yfou}{20}{#3}
\put(\xone,\ythr){\makebox(0,0){$\vdots$}}
\put(\xthr,\ysix){\makebox(0,0){$\cdots$}}
\put(\xone,\yone){\makebox(0,0){0}}
\let\J=#6 \let\I=#9
\put(\xone,\ysix){\makebox(0,0){${#5}$}}
\put(\xsix,\ysix){\makebox(0,0){0}}
\advance \xtwo by 20
\advance \xfiv by -20
\advance \ytwo by 20
\advance \yfiv by -20
\put(\xthr,\ythr){\makebox(0,0){$\cdot$}}
\put(\xfiv,\yfiv){\makebox(0,0){$\cdot$}}
\put(\xtwo,\ytwo){\makebox(0,0){$\cdot$}}
}

\nc\fStrd[9]{\icoo{#1}\jcoo{#2}
\eone=\varthr \divide \eone by 2
\xthr=\varone \advance \xthr by \eone
\advance \eone by -35
\divide \eone by 2
\fone=\varfou \divide \fone by 2
\ythr=\vartwo \advance \ythr by \fone
\advance \fone by -35
\divide \fone by 2
\xone=\varone \advance \xone by \eone
\xtwo=\xone \advance \xtwo by \eone
\xsix=\varone \advance \xsix by \varthr 
\advance \xsix by -\eone
\xfiv =\xsix \advance \xfiv by -\eone
\xfou =\xfiv \advance \xfou by -20
\yone=\vartwo \advance \yone by \fone
\ytwo=\yone \advance \ytwo by \fone
\ysix=\vartwo \advance \ysix by \varfou 
\advance \ysix by -\fone
\yfiv =\ysix \advance \yfiv by -\fone
\yfou =\yfiv \advance \yfou by -20
\setra{\xtwo}{\ysix}{20}{#4}
\setra{\xfou}{\ysix}{20}{#4}
\setua{\xone}{\ytwo}{20}{#3}
\setua{\xone}{\yfou}{20}{#3}
\put(\xone,\ythr){\makebox(0,0){$\vdots$}}
\put(\xthr,\ysix){\makebox(0,0){$\cdots$}}
\let\J=#6 \let\I=#8
\put(\xone,\yone){\makebox(0,0){${#5}$}}
\let\J=#6 \let\I=#9
\put(\xone,\ysix){\makebox(0,0){${#5}$}}
\let\J=#7 \let\I=#9
\put(\xsix,\ysix){\makebox(0,0){${#5}$}}
\advance \xtwo by 20
\advance \xfiv by -20
\advance \ytwo by 20
\advance \yfiv by -20
\put(\xthr,\ythr){\makebox(0,0){$\cdot$}}
\put(\xfiv,\yfiv){\makebox(0,0){$\cdot$}}
\put(\xtwo,\ytwo){\makebox(0,0){$\cdot$}}
}

\nc\Strdp[5]{\Ipar{#1}\Jpar{#2}
\Strd{#1}{#2}{#3}{#4}{#5}{\varfiv}{\varsix}{\varsev}{\vareig}
}

\nc\fStrdp[5]{\Ipar{#1}\Jpar{#2}
\fStrd{#1}{#2}{#3}{#4}{#5}{\varfiv}{\varsix}{\varsev}{\vareig}
}

\nc\Sbu[3]{ \icoo{#1} \jcoo{#2}
\hone=\varthr \divide \hone by 2
\advance \hone by -35
\divide \hone by 2
\htwo =\varone \advance \htwo by \hone
\hthr =\varone \advance \hthr by \varthr
\advance \hthr by -\hone
\setua\htwo\vartwo\varfou{#3}
\setua\hthr\vartwo\varfou{#3}
}

\nc\tinySbu[3]{ \icoo{#1} \jcoo{#2}
\hone=\varthr \divide \hone by 2
\advance \hone by -25
\divide \hone by 2
\htwo =\varone \advance \htwo by \hone
\hthr =\varone \advance \hthr by \varthr
\advance \hthr by -\hone
\setua\htwo\vartwo\varfou{#3}
\setua\hthr\vartwo\varfou{#3}
}

\nc\Sfbu[4]{\coosq{#1}{#1}{#2}{#3} 
\hone=\varthr \divide \hone by 2
\advance \hone by -35
\divide \hone by 2
\htwo =\varone \advance \htwo by \hone
\hthr =\varone \advance \hthr by \varthr
\advance \hthr by -\hone
\setua\htwo\vartwo\varfou{#4}
\setua\hthr\vartwo\varfou{#4}
}

\nc\Sbr[3]{ \icoo{#1} \jcoo{#2}
\hone=\varfou \divide \hone by 2
\advance \hone by -35
\divide \hone by 2
\htwo =\vartwo \advance \htwo by \hone
\hthr =\vartwo \advance \hthr by \varfou
\advance \hthr by -\hone
\setra\varone\htwo\varthr{#3}
\setra\varone\hthr\varthr{#3}
}

\nc\Sbrsmall[3]{ \icoo{#1} \jcoo{#2}
\hone=\varfou \divide \hone by 2
\advance \hone by -35
\divide \hone by 2
\htwo =\vartwo \advance \htwo by \hone
\hthr =\vartwo \advance \hthr by \varfou
\advance \hthr by -\hone
\advance\varthr by -4
\advance\varone by 2
\setra\varone\htwo\varthr{#3}
\setra\varone\hthr\varthr{#3}
}

\nc\fSrli[5]{\icoo{#1}\jcoo{#2}
\eone=\varthr \divide \eone by 2
\xthr=\varone \advance \xthr by \eone
\advance \eone by -35
\divide \eone by 2
\fone=\varfou \divide \fone by 2
\ythr=\vartwo \advance \ythr by \fone
\xone=\varone \advance \xone by \eone
\xtwo=\xone \advance \xtwo by \eone
\xsix=\varone \advance \xsix by \varthr 
\advance \xsix by -\eone
\xfiv =\xsix \advance \xfiv by -\eone
\xfou =\xfiv \advance \xfou by -20
\setra{\xtwo}{\ythr}{20}{#3}
\setra{\xfou}{\ythr}{20}{#3}
\put(\xthr,\ythr){\makebox(0,0){$\cdots$}}
\put(\xone,\ythr){\makebox(0,0){${#4}$}}
\put(\xsix,\ythr){\makebox(0,0){${#5}$}}
}

\nc\Srli[6]{\icoo{#1}\jcoo{#2}
\eone=\varthr \divide \eone by 2
\xthr=\varone \advance \xthr by \eone
\advance \eone by -35
\divide \eone by 2
\fone=\varfou \divide \fone by 2
\ythr=\vartwo \advance \ythr by \fone
\xone=\varone \advance \xone by \eone
\xtwo=\xone \advance \xtwo by \eone
\xsix=\varone \advance \xsix by \varthr 
\advance \xsix by -\eone
\xfiv =\xsix \advance \xfiv by -\eone
\xfou =\xfiv \advance \xfou by -20
\setra{\xtwo}{\ythr}{20}{#3}
\setra{\xfou}{\ythr}{20}{#3}
\put(\xthr,\ythr){\makebox(0,0){$\cdots$}}
\let\J=#5 
\put(\xone,\ythr){\makebox(0,0){${#4}$}}
\let\J=#6 
\put(\xsix,\ythr){\makebox(0,0){${#4}$}}
}

\nc\Srlip[4]{\Ipar{#1}
\Srli{#1}{#2}{#3}{#4}{\varfiv}{\varsix}
}

\nc\Suli[6]{\icoo{#1}\jcoo{#2}
\eone=\varthr \divide \eone by 2
\xthr=\varone \advance \xthr by \eone
\fone=\varfou \divide \fone by 2
\ythr=\vartwo \advance \ythr by \fone
\advance \fone by -35
\divide \fone by 2
\yone=\vartwo \advance \yone by \fone
\ytwo=\yone \advance \ytwo by \fone
\ysix=\vartwo \advance \ysix by \varfou 
\advance \ysix by -\fone
\yfiv =\ysix \advance \yfiv by -\fone
\yfou =\yfiv \advance \yfou by -20
\advance\ythr by 3
\setua{\xthr}{\ytwo}{20}{#3}
\setua{\xthr}{\yfou}{20}{#3}
\put(\xthr,\ythr){\makebox(0,0){$\vdots$}}
\let\I=#5
\put(\xthr,\yone){\makebox(0,0){${#4}$}}
\let\I=#6
\put(\xthr,\ysix){\makebox(0,0){${#4}$}}
}

\nc\Sulip[4]{\Jpar{#2}
\Suli{#1}{#2}{#3}{#4}{\varsev}{\vareig}
}

\nc\Sposar[4]{\icoo{#1}\jcoo{#2}
\xone=\varthr \divide \xone by 2
\advance \xone by \varone
\yone=\varfou \divide \yone by 2
\advance \yone by \vartwo
\if#4r
\setra\varone\yone\varthr{#3}
\else\if#4u
\setua\xone\vartwo\varfou{#3}
\fi\fi
}

\nc\fSposar[5]{
\if#5r
\horzarr{#1}{#2}{#3}
\else\if#5u
\vertarr{#1}{#2}{#3}
\fi\fi
\xone=\varthr \divide \xone by 2
\advance \xone by \varone
\yone=\varfou \divide \yone by 2
\advance \yone by \vartwo
\if#5r
\setra\varone\yone\varthr{#4}
\else\if#5u
\setua\xone\vartwo\varfou{#4}
\fi\fi
}

\nc\fmbx[4]{\coosq{#1}{#2}{#3}{#4}
\advance \varone by 2
\advance \vartwo by 2
\advance \varthr by -4
\advance \varfou by -4
\thicklines
\put(\varone,\vartwo){\framebox(\varthr,\varfou){}}
\thinlines
}

\nc\dshbx[4]{\coosq{#1}{#2}{#3}{#4}
\advance \varone by -5
\advance \vartwo by -5
\advance \varthr by 10
\advance \varfou by 10
\put(\varone,\vartwo){\dashbox{5}(\varthr,\varfou){}}
}

\nc\divli[3]{\icoo{#1}\jcoo{#2}
\xone=\varthr \divide \xone by 2
\advance \xone by \varone
\yone =\varfou \divide \yone by 2
\advance \yone by \vartwo
\advance \varthr by -4
\advance \varfou by -4
\advance \varone by 2
\advance \vartwo by 2
\if#3v 
\put(\xone,\vartwo){\line(0,1){\varfou}}
\else\if#3h
\put(\varone,\yone){\line(1,0){\varthr}}
\fi\fi
}

\nc\ovrbrc[4]{\icoo{#1}\jcoo{#2}
\advance \vartwo by \varfou
\advance\vartwo by -5
\advance \varthr by -10
\advance \varthr by #4
\put(\varone,\vartwo){\makebox(\varthr,0)[b]{$
\overbrace{\hbox{\makebox[\varthr\unitlength]{}}}^{#3}$}}
}

\nc\undrbrc[4]{\icoo{#1}\jcoo{#2}
\advance\vartwo by -5
\advance \varthr by -10
\advance \varthr by #4
\put(\varone,\vartwo){\makebox(\varthr,0)[t]{$
\underbrace{\hbox{\makebox[\varthr\unitlength]{}}}_{#3}$}}
}

\nc\plundrbrc[4]{\icoo{#1}\jcoo{#2}
\advance \varthr by #4
\put(\varone,\vartwo){\makebox(\varthr,0)[t]{$
\underbrace{\hbox{\makebox[\varthr\unitlength]{}}}_{#3}$}}
}

\nc\lftbrc[5]{\icoo{#1}\jcoo{#2}
\advance \varfou by -6
\advance \varfou by #5
\advance \varone by 10
\fone=#5
\advance \vartwo by -\fone
\advance \vartwo by 3
\put(\varone,\vartwo){\makebox(0,\varfou)[r]
{${\left( \begin{array}{c}{\scriptstyle #3 }
\\{\scriptstyle #4 } \end{array} \right) }\left\{ \begin{array}{c} 
\hbox{\begin{picture}(0,\varfou)
\put(0,0){\makebox(0,\varfou){}}
\end{picture}}
\end{array} \right. $}}
}

\nc\pllftbrc[4]{\icoo{#1}\jcoo{#2}
\advance \varfou by -6
\advance \varfou by #4
\advance \varone by 10
\fone=#4
\advance \vartwo by -\fone
\advance \vartwo by 3
\put(\varone,\vartwo){\makebox(0,\varfou)[r]
{${#3 }\left\{ \begin{array}{c} 
\hbox{\begin{picture}(0,\varfou)
\put(0,0){\makebox(0,\varfou){}}
\end{picture}}
\end{array} \right. $}}
}

\nc\lftmtrx[5]{\icoo{#1}\jcoo{#2}
\fone=#4
\advance \varone by -\fone
\advance\varthr by \fone
\fone=#5
\advance\varthr by \fone
\put(\varone,\vartwo){\makebox(\varthr,\varfou)
{$\left( \begin{array}{c} 
\hbox{\begin{picture}(\varthr,\varfou)
\put(0,0){\makebox(\varthr,\varfou){}}
\end{picture}}
\end{array} \right)^{#3} $}}
}

\nc\rghtbrc[5]{\icoo{#1}\jcoo{#2}
\advance \varfou by -6
\advance \varfou by #5
\advance \varone by \varthr
\advance \varone by -10
\fone=#5
\advance \vartwo by -\fone
\put(\varone,\vartwo){\makebox(0,\varfou)[l]{$\left. \begin{array}{c} 
\hbox{\begin{picture}(0,\varfou)
\put(0,0){\makebox(0,\varfou)}
\end{picture}}
\end{array} \right\}{\left( \begin{array}{c}{\scriptstyle #3 }
\\{\scriptstyle  #4}   \end{array} \right) } $}}
}

\nc\setwidths[5]{
\ione=#1 \multiply \ione by 10
\itwo=#2 \multiply \itwo by 10
\ithr=#3 \multiply \ithr by 10
\ifou=#4 \multiply \ifou by 10
\ifiv=#5 \multiply \ifiv by 10
}

\nc\setarrows[5]{
\jone=#1 \multiply \jone by 10
\jtwo=#2 \multiply \jtwo by 10
\jthr=#3 \multiply \jthr by 10
\jfou=#4 \multiply \jfou by 10
\jfiv=#5 \multiply \jfiv by 10
}

\nc\rowsdec[1]{
\ifnum#1=1 \varone=\ione \vartwo=\jone
\else\ifnum#1=2 \varone=\itwo \vartwo=\jtwo
\else\ifnum#1=3 \varone=\ithr \vartwo=\jthr
\else\ifnum#1=4 \varone=\ifou \vartwo=\jfou
\else\ifnum#1=5 \varone=\ifiv \vartwo=\jfiv
\fi\fi\fi\fi\fi
}

\nc\Beginsnake[2]{\tione=#1 \tjone=#2 
\advance \tione by -1
\multiply \tione by \tjone
\advance \tione by 40
\begin{center}
\begin{picture}(400,\tione)
}

\nc\Arrow[1]{\yone=#1 
\rowsdec\yone
\xone=\varone \multiply \xone by 3
\divide \xone by 2
\advance \xone by \vartwo
\advance \yone by 1
\rowsdec\yone
\xtwo=\varone \multiply \xtwo by 3
\divide \xtwo by 2
\advance \xtwo by \vartwo
\ytwo =\tjone \divide \ytwo by 2
\advance \yone by -2
\multiply \yone by \tjone
\advance \yone by 20
\advance \yone by \ytwo
\multiply \yone by -1
\advance \yone by \tione
\advance \xone by 8
\advance \xtwo by 8
\put(200,\yone){\line(1,0){\xone}}
\put(200,\yone){\line(-1,0){\xtwo}}
\advance \yone by \ytwo
\advance \xone by 200
\put(\xone,\yone){\line(-1,0){8}}
\ythr=\ytwo \divide \ythr by 2
\advance \yone by -\ythr
\put(\xone,\yone){\oval(\ytwo,\ytwo)[r]}
\advance \yone by -\ytwo
\multiply \xtwo by -1
\advance \xtwo by 200
\put(\xtwo,\yone){\oval(\ytwo,\ytwo)[l]}
\advance \yone by -\ythr
\put(\xtwo,\yone){\vector(1,0){8}}
}

\nc\Row[5]{\rowsdec{#1}
\fone=\varone \divide \fone by 2
\gthr=200
\hthr=\gthr \advance \hthr by \fone
\advance \varone by \vartwo
\gfou =\gthr \advance \gfou by \varone
\gfiv =\gfou \advance \gfiv by \varone
\gtwo =\gthr \advance \gtwo by -\varone
\gone =\gtwo \advance \gone by -\varone
\htwo =\hthr \advance \htwo by -\varone
\hone =\htwo \advance \hone by -\varone
\hfou =\hthr \advance \hfou by \varone
\yone =#1
\advance \yone by -1
\multiply \yone by \tjone
\advance \yone by 20
\multiply \yone by -1
\advance \yone by \tione
\put(\gtwo,\yone){\makebox(0,0){$#2$}}
\put(\gthr,\yone){\makebox(0,0){$#3$}}
\put(\gfou,\yone){\makebox(0,0){$#4$}}
\put(\htwo,\yone){\vector(1,0){\vartwo}}
\put(\hthr,\yone){\vector(1,0){\vartwo}}
\if#5l 
\put(\gone,\yone){\makebox(0,0){0}}
\put(\hone,\yone){\vector(1,0){\vartwo}}
\else\if#5r 
\put(\gfiv,\yone){\makebox(0,0){0}}
\put(\hfou,\yone){\vector(1,0){\vartwo}}
\fi\fi
}

\nc\Endsnake{
\end{picture}
\end{center}
}

\nc\UBP[3]{
\put(#1,#2){\makebox(0,0){\usebox{#3}}}
}

\def\UB(#1,#2)#3{
\put(#1,#2){\makebox(0,0){\usebox{#3}}}
}

\nc\SB[2]{\savebox{#1}
{\begin{picture}(\tinin,\tjnin)
{#2}
\end{picture}}}

\nc\Var[5]{\fone =#3
\divide \fone by 2
\advance \fone by #2
\put(#1,\fone){\vector(0,-1){#3}}
\eone =#1
\advance \eone by -2
\put(\eone,#2){\makebox(0,0)[r]{$#4$}}
\advance \eone by 4
\put(\eone,#2){\makebox(0,0)[l]{$#5$}}
}

\def\Va(#1,#2,#3)#4#5{\fone =#3
\divide \fone by 2
\advance \fone by #2
\put(#1,\fone){\vector(0,-1){#3}}
\eone =#1
\advance \eone by -2
\put(\eone,#2){\makebox(0,0)[r]{$#4$}}
\advance \eone by 4
\put(\eone,#2){\makebox(0,0)[l]{$#5$}}
}

\nc\Vin[5]{\fone =#3
\divide \fone by 2
\advance \fone by #2
\advance \fone by -2
\yone =#3
\advance \yone by -2
\put(#1,\fone){\vector(0,-1){\yone}}
\eone =#1
\advance \eone by -2
\put(\eone,#2){\makebox(0,0)[r]{$#4$}}
\advance \eone by 4
\put(\eone,#2){\makebox(0,0)[l]{$#5$}}
\put(\eone,\fone){\oval(4,4)[t]}
\advance \eone by 2
\put(\eone,\fone){\line(0,-1){2}}
}

\def\Vi(#1,#2,#3)#4#5{\fone =#3
\divide \fone by 2
\advance \fone by #2
\advance \fone by -2
\yone =#3
\advance \yone by -2
\put(#1,\fone){\vector(0,-1){\yone}}
\eone =#1
\advance \eone by -2
\put(\eone,#2){\makebox(0,0)[r]{$#4$}}
\advance \eone by 4
\put(\eone,#2){\makebox(0,0)[l]{$#5$}}
\put(\eone,\fone){\oval(4,4)[t]}
\advance \eone by 2
\put(\eone,\fone){\line(0,-1){2}}
}

\nc\Rin[4]{
\ythr=#2
\advance \ythr by 2
\put(#1,\ythr){\makebox(0,0)[b]{$#4$}}
\eone=#3
\divide\eone by 2
\etwo =#1 \advance \etwo by -\eone
\advance \etwo by 2
\xone =#3
\advance \xone by -2
\put(\etwo,#2){\vector(1,0){\xone}}
\yone =#2
\advance \yone by 2
\put(\etwo,\yone){\oval(4,4)[l]}
\advance \yone by 2
\put(\etwo,\yone){\line(1,0){2}}
}

\def\Ri(#1,#2,#3)#4{
\ythr=#2
\advance \ythr by 2
\put(#1,\ythr){\makebox(0,0)[b]{$#4$}}
\eone=#3
\divide\eone by 2
\etwo =#1 \advance \etwo by -\eone
\advance \etwo by 2
\xone =#3
\advance \xone by -2
\put(\etwo,#2){\vector(1,0){\xone}}
\yone =#2
\advance \yone by 2
\put(\etwo,\yone){\oval(4,4)[l]}
\advance \yone by 2
\put(\etwo,\yone){\line(1,0){2}}
}

\nc\slar[6]{\fone =#3
\multiply \fone by #5
\divide \fone by #4
\eone =#3
\divide \eone by 2
\divide \fone by 2
\yone =#2 \advance \yone by 4
\xsix =#1
\ifnum#4 < 0  \ifnum#5 < #4 
\advance \xsix by -4\fi\fi
\ifnum#4 < 0  \ifnum#5 = #4 
\advance \xsix by -3\fi\fi
\ifnum#4 < 0 \xsev=#5
\multiply \xsev by -1
\ifnum\xsev < #4 
\advance \xsix by 4\fi
\ifnum\xsev = #4 
\advance \xsix by 3\fi
\fi
\ifnum#4 > 0 \ifnum#5 > #4 
\advance \xsix by -4\fi
\ifnum#5=#4 \advance \xsix by -3\fi
\fi
\ifnum#4 > 0 \xsev=#5
\multiply \xsev by -1
\ifnum\xsev > #4 
\advance \xsix by 4\fi
\ifnum\xsev = #4 
\advance \xsix by 3\fi
\fi
\put(\xsix,\yone){\makebox(0,0)[b]{$#6$}}
\advance \yone by -4
\ifnum#4 >0
\advance \yone by -\fone
\xone =#1
\advance \xone by -\eone
\put(\xone,\yone){\vector(#4,#5){#3}}
\else\advance \yone by \fone
\xone =#1
\advance \xone by \eone
\put(\xone,\yone){\vector(#4,#5){#3}}
\fi
}

\def\sla(#1,#2,#3,#4,#5)#6{\fone =#3
\multiply \fone by #5
\divide \fone by #4
\eone =#3
\divide \eone by 2
\divide \fone by 2
\yone =#2 \advance \yone by 4
\xsix =#1
\ifnum#4 < 0  \ifnum#5 < #4 
\advance \xsix by -4\fi\fi
\ifnum#4 < 0  \ifnum#5 = #4 
\advance \xsix by -3\fi\fi
\ifnum#4 < 0 \xsev=#5
\multiply \xsev by -1
\ifnum\xsev < #4 
\advance \xsix by 4\fi
\ifnum\xsev = #4 
\advance \xsix by 3\fi
\fi
\ifnum#4 > 0 \ifnum#5 > #4 
\advance \xsix by -4\fi
\ifnum#5=#4 \advance \xsix by -3\fi
\fi
\ifnum#4 > 0 \xsev=#5
\multiply \xsev by -1
\ifnum\xsev > #4 
\advance \xsix by 4\fi
\ifnum\xsev = #4 
\advance \xsix by 3\fi
\fi
\put(\xsix,\yone){\makebox(0,0)[b]{$#6$}}
\advance \yone by -4
\ifnum#4 >0
\advance \yone by -\fone
\xone =#1
\advance \xone by -\eone
\put(\xone,\yone){\vector(#4,#5){#3}}
\else\advance \yone by \fone
\xone =#1
\advance \xone by \eone
\put(\xone,\yone){\vector(#4,#5){#3}}
\fi
}

\nc\Rar[4]{\slar{#1}{#2}{#3}{1}{0}{#4}
}

\def\Ra(#1,#2,#3)#4{\slar{#1}{#2}{#3}{1}{0}{#4}
}

\nc\Veq[3]{\yone=#3
\divide \yone by 2
\advance \yone by #2
\xone =#1
\advance \xone by -2
\put(\xone,\yone){\line(0,-1){#3}}
\advance \xone by 4
\put(\xone,\yone){\line(0,-1){#3}}
}

\def\Vq(#1,#2,#3){\yone=#3
\divide \yone by 2
\advance \yone by #2
\xone =#1
\advance \xone by -2
\put(\xone,\yone){\line(0,-1){#3}}
\advance \xone by 4
\put(\xone,\yone){\line(0,-1){#3}}
}

\nc\Req[3]{\yone=#3
\divide \yone by 2
\advance \yone by #1
\xone =#2
\advance \xone by -2
\put(\yone,\xone){\line(-1,0){#3}}
\advance \xone by 4
\put(\yone,\xone){\line(-1,0){#3}}
}

\def\Rq(#1,#2,#3){\yone=#3
\divide \yone by 2
\advance \yone by #1
\xone =#2
\advance \xone by -2
\put(\yone,\xone){\line(-1,0){#3}}
\advance \xone by 4
\put(\yone,\xone){\line(-1,0){#3}}
}

\nc\Vsim[3]{\fone =#3
\divide \fone by 2
\advance \fone by #2
\put(#1,\fone){\vector(0,-1){#3}}
\eone =#1
\advance \eone by -2
\put(\eone,#2){\makebox(0,0)[r]{$\scriptstyle{\mid}$}}
\advance \eone by 4
\put(\eone,#2){\makebox(0,0)[l]{$\wr$}}
}

\def\Vs(#1,#2,#3){\fone =#3
\divide \fone by 2
\advance \fone by #2
\put(#1,\fone){\vector(0,-1){#3}}
\eone =#1
\advance \eone by -2
\put(\eone,#2){\makebox(0,0)[r]{$\scriptstyle{\mid}$}}
\advance \eone by 4
\put(\eone,#2){\makebox(0,0)[l]{$\wr$}}
}

\nc\Rsim[3]{\yone=#3
\divide \yone by 2
\ytwo =#1 \advance \ytwo by -\yone
\put(\ytwo,#2){\vector(1,0){#3}}
\xone =#2
\advance \xone by 2
\put(#1,\xone){\makebox(0,0)[t]{$-$}}
\put(#1,\xone){\makebox(0,0)[b]{$\sim$}}
}

\def\Rs(#1,#2,#3){\yone=#3
\divide \yone by 2
\ytwo =#1 \advance \ytwo by -\yone
\put(\ytwo,#2){\vector(1,0){#3}}
\xone =#2
\advance \xone by 2
\put(#1,\xone){\makebox(0,0)[t]{$-$}}
\put(#1,\xone){\makebox(0,0)[b]{$\sim$}}
}

\nc\RLar[5]{\xone=#3
\divide \xone by 2
\xtwo =#1
\yone =#2
\advance \xtwo by -\xone
\advance \yone by 4
\put(\xtwo,\yone){\vector(1,0){#3}}
\advance \yone by 2
\put(#1,\yone){\makebox(0,0)[b]{$#4$}}
\advance \yone by -10
\advance \xtwo by #3
\put(\xtwo,\yone){\vector(-1,0){#3}}
\advance \yone by -2
\put(#1,\yone){\makebox(0,0)[t]{$#5$}}
}

\def\RLa(#1,#2,#3)#4#5{\xone=#3
\divide \xone by 2
\xtwo =#1
\yone =#2
\advance \xtwo by -\xone
\advance \yone by 4
\put(\xtwo,\yone){\vector(1,0){#3}}
\advance \yone by 2
\put(#1,\yone){\makebox(0,0)[b]{$#4$}}
\advance \yone by -10
\advance \xtwo by #3
\put(\xtwo,\yone){\vector(-1,0){#3}}
\advance \yone by -2
\put(#1,\yone){\makebox(0,0)[t]{$#5$}}
}

\nc\mtrx[1]{\left( \begin{array}{c}
\hbox{\raisebox{-5pt}{\usebox#1}} \end{array} \right)
}

\nc\PT[1]{\begin{array}{c}\usebox#1
\end{array}}

\nc\BP{
\begin{center}
\begin{picture}(\tinin,\tjnin)
}

\nc\EP{
\end{picture}
\end{center}
}

\def\PP(#1,#2)#3#4{
\put(#1,#2){\makebox(0,0){
$\begin{array}{ccc}
\mtrx{#3}  &  \times &  \mtrx{#4}
\end{array}$
}}}

\nc\PIP[4]{
\put(#1,#2){\makebox(0,0){
$\begin{array}{ccc}
\mtrx{#3}  &  \times &  \mtrx{#4}
\end{array}$
}}}

\def\BDI(#1,#2){
\begin{center}
\eone=#1 \fone=#2
\begin{picture}(\eone,\fone)
}

\nc\EDI{
\end{picture}
\end{center}
}

\nc\tqu{\ct^{\hbox{
\begin{picture}(6,12)
\put(-1,9){\makebox(0,0){\rm ?}}
\put(4,4){\makebox(0,0){\rm ?}}
\end{picture}}}_{*}
}

\nc\tqus[2]{\ct^{\hbox{
\begin{picture}(6,12)
\put(-1,9){\makebox(0,0){\rm ?}}
\put(4,4){\makebox(0,0){\rm ?}}
\end{picture}}}_{[{#1},{#2}]}
}

\nc\tqul{
\ct^{\hbox{\rm ?}}_{*}
}

\nc\Q{{\I\J}}

\newcount\dpone
\newcount\htone
\newcount\wdone

\nc\setstar{\setbox1=\hbox{${\bf *}$}
\dpone=\dp1
\htone=\ht1
\wdone=\wd1}

\nc\vda[3]{\uvariablearrowpos{#1}{#2}{#3}
\put(\eone ,\gone){\line(0,1){\fone}}
\put(\etwo ,\gone){\oval(4,4)[tr]}
\put(\ethr ,\gone){\oval(4,4)[tl]}
}

\nc\vla[3]{\rvariablearrowpos{#1}{#2}{#3}
\put(\gone ,\eone){\line(1,0){\fone}}
\put(\gone ,\ethr){\oval(4,4)[br]}
\put(\gone ,\etwo){\oval(4,4)[tr]}
}

\nc\vuBa[3]{
\vda{#1}{#2}{#3}
\put(\eone,\gsev){\makebox(0,0)[t]{$\bullet$}}
}

\nc\vrBa[3]{
\vla{#1}{#2}{#3}
\put(\gsev,\eone){\makebox(0,0)[r]{$\bullet$}}
}

\nc\vuaT[3]{
\vua{#1}{#2}{#3}
\put(\eone,\gfiv){\makebox(0,0){\hbox{\scriptsize \bf T}}}
}

\nc\vraT[3]{
\vrM{#1}{#2}{#3}
\put(\gfiv,\eone){\makebox(0,0){\hbox{\scriptsize \bf T}}}
}

\nc\vuaMe[3]{
\vuM{#1}{#2}{#3}
\put(\etwo ,\gsix){\oval(4,4)[br]}
\put(\ethr ,\gsix){\oval(4,4)[bl]}
}

\nc\vraMe[3]{
\vrM{#1}{#2}{#3}
\put(\gsix,\etwo){\oval(4,4)[tl]}
\put(\gsix,\ethr){\oval(4,4)[bl]}
}

\nc\vuaIe[3]{
\vuI{#1}{#2}{#3}
\put(\etwo ,\gsix){\oval(4,4)[br]}
\put(\ethr ,\gsix){\oval(4,4)[bl]}
}

\nc\vraIe[3]{
\vrI{#1}{#2}{#3}
\put(\gsix,\etwo){\oval(4,4)[tl]}
\put(\gsix,\ethr){\oval(4,4)[bl]}
}

\nc\vuat[3]{
\vua{#1}{#2}{#3}
\put(\eone,\gfiv){\makebox(0,0){$\scriptstyle \ct$}}
}

\nc\vrat[3]{
\vra{#1}{#2}{#3}
\put(\gfiv,\eone){\makebox(0,0){$\scriptstyle \ct$}}
}

\nc\vuM[3]{\setstar
\divide\htone by 2
\advance\dpone by \htone
\divide\dpone by 65536
\uvariablearrowpos{#1}{#2}{#3}
\advance\fone by -\dpone
\advance\gone by \dpone
\put(\eone ,\gone){\line(0,1){\fone}}
\put(\etwo ,\gsev){\oval(4,4)[br]}
\put(\ethr ,\gsev){\oval(4,4)[bl]}
\put(\eone ,\gone){\makebox(0,0){\box1}}
}

\nc\vrM[3]{\setstar
\divide\wdone by 2
\divide\wdone by 65536
\rvariablearrowpos{#1}{#2}{#3}
\advance\fone by -\wdone
\advance\gone by \wdone
\put(\gone,\eone){\line(1,0){\fone}}
\put(\gsev,\etwo){\oval(4,4)[tl]}
\put(\gsev,\ethr){\oval(4,4)[bl]}
\put(\gone,\eone){\makebox(0,0){\box1}}
}

\nc\vuMs[3]{
\vuM{#1}{#2}{#3}
\advance \gsix by -4
\put(\etwo,\gsix){\oval(4,4)[tr]}
\put(\ethr,\gsix){\oval(4,4)[tl]}
\advance \gsix by 2
\put(\etwo,\gsix){\line(0,1){2}}
\put(\ethr,\gsix){\line(0,1){2}}
}

\nc\vrMs[3]{
\vrM{#1}{#2}{#3}
\advance \gsix by -4
\put(\gsix,\etwo){\oval(4,4)[tr]}
\put(\gsix,\ethr){\oval(4,4)[br]}
\advance \gsix by 2
\put(\gsix,\etwo){\line(1,0){2}}
\put(\gsix,\ethr){\line(1,0){2}}
}

\nc\vuE[3]{\setstar
\divide\htone by 2
\divide\htone by 65536
\uvariablearrowpos{#1}{#2}{#3}
\advance\fone by -\htone
\advance\gsev by -\htone
\put(\eone ,\gone){\line(0,1){\fone}}
\put(\eone ,\gsev){\makebox(0,0){\box1}}
}

\nc\vrE[3]{\setstar
\divide\wdone by 2
\divide\wdone by 65536
\rvariablearrowpos{#1}{#2}{#3}
\advance\fone by -\wdone
\advance\gsev by -\wdone
\put(\gone,\eone){\line(1,0){\fone}}
\put(\gsev,\eone){\makebox(0,0){\box1}}
}

\nc\vuI[3]{\setstar
\vua{#1}{#2}{#3}
\put(\eone,\gfiv){\makebox(0,0){\box1}}
}

\nc\vui[3]{\setstar
\uvariablearrowpos{#1}{#2}{#3}
\put(\eone ,\gone){\vector(0,1){\fone}}
\put(\eone,\gfiv){\makebox(0,0){\box1}}
}

\nc\vrI[3]{\setstar
\vra{#1}{#2}{#3}
\put(\gfiv,\eone){\makebox(0,0){\box1}}
}

\nc\vri[3]{
\setstar
\rvariablearrowpos{#1}{#2}{#3}
\put(\gone ,\eone){\vector(1,0){\fone}}
\put(\gfiv,\eone){\makebox(0,0){\box1}}
}

\nc\isoar{
\put(6,5.5){\vector(2,1){14}}
\put(12,8.5){\makebox(0,0){\hbox{${\bf *}$}}}
}

\nc\barbi[2]{\hbox{\begin{picture}(21,21)
\setua25{16}{#1}
\setra52{16}{#2}
\isoar
\end{picture}}}

\nc\Callet[1]{
\if#1s{\cal S}\fi
\if#1r{\cal R}\fi
\if#1t{\cal T}\fi
\if#1S{\bf S}\fi
\if#1R{\bf R}\fi
\if#1T{\bf T}\fi
}

\nc\BX[5]{\icoo{#1}\jcoo{#2}
\put(\varone ,\vartwo){\framebox(\varthr ,\varfou){$\Callet{#3}^{\barbi{#4}{#5}}$}}
}

\nc\BXS[6]{\icoo{#1}\jcoo{#2}
\put(\varone ,\vartwo){\framebox(\varthr ,\varfou){$\Callet{#3}^{\barbi{#4}{#5}}{#6}$}}
}

\nc\HBX[7]{\icoo{#1}\jcoo{#2}
\put(\varone ,\vartwo){\framebox(\varthr ,\varfou){$\Callet{#3}^{\barbi{#4}{#5}}{#6}$}}
\hone=\varone \advance \hone by 5
\htwo=\varone \advance \htwo by \varthr
\advance \htwo by -5
\hthr =\vartwo \advance \hthr by 5
\hfou =\vartwo \advance \hfou by \varfou
\advance \hfou by -5
\hfiv =\varthr \advance \hfiv by -10
\hsix =\varfou \advance \hsix by -10
\if#7r\put(\hone ,\hthr){\vector(1,0){\hfiv}}
\else\if#7l\put(\htwo ,\hthr){\vector(-1,0){\hfiv}}
\else\if#7u\put(\hone ,\hthr){\vector(0,1){\hsix}}
\else\if#7d\put(\hone ,\hfou){\vector(0,-1){\hsix}}
\fi\fi\fi\fi
}

\nc\arubi[1]{\hbox{\begin{picture}(20,21)
\setua25{16}{#1}
\isoar
\end{picture}}}

\nc\arrbi[1]{\hbox{\begin{picture}(21,13)(3,0)
\setra52{16}{#1}
\isoar
\end{picture}}}

\def\nothing{}

\nc\LI[6]{\icoo{#1}\jcoo{#2}
\if#6r
\advance \varone by \varthr
\put(\varone,\vartwo){\line(0,1){\varfou}}
\advance\varone by -2
\def\test{#5}
\ifx\test\nothing
\put(\varone,\vartwo){\makebox(0,\varfou)[r]{${\Callet{#3}}^{\arubi{#4}}$}}
\else
\put(\varone,\vartwo){\makebox(0,\varfou)[r]{${\Callet{#3}}^{\arubi{#4}}\!\!\!{#5}$}}\fi
\else\if#6l
\put(\varone,\vartwo){\line(0,1){\varfou}}
\advance \varone by 2
\put(\varone,\vartwo){\makebox(0,\varfou)[l]{${\Callet{#3}}^{\arubi{#4}}\!\!\!{#5}$}}
\else\if#6b
\put(\varone,\vartwo){\line(1,0){\varthr}}
\advance \vartwo by 6
\put(\varone,\vartwo){\makebox(\varthr,0)[b]{${\Callet{#3}}^{\arrbi{#4}}{#5}$}}
\else\if#6t
\advance \vartwo by \varfou
\put(\varone,\vartwo){\line(1,0){\varthr}}
\advance \vartwo by -2
\put(\varone,\vartwo){\makebox(\varthr,0)[t]{${\Callet{#3}}^{\arrbi{#4}}{#5}$}}
\fi\fi\fi\fi
}

\nc\xdifr[7]{\coosq{#1}{#2}{#3}{#4}
\fone=\varone \advance \fone by \varthr
\advance \fone by 5
\eone=\varfou \divide \eone by 2
\advance \eone by \vartwo
\advance \eone by -10
\etwo=\varfou \advance \etwo by \vartwo
\advance \varfou by 20
\put(\fone,\eone){\oval(30,\varfou)[r]}
\put(\fone,\etwo){\line(-1,0){#7}}
\advance \vartwo by -5
\multiply \varthr by 2
\advance \varthr by 10
\put(\fone,\vartwo){\oval(\varthr,30)[bl]}
\put(\varone,\vartwo){\vector(0,1){#6}}
\advance \vartwo by -20
\advance \fone by 20
\put(\fone,\vartwo){\makebox(0,0){$#5$}}
}

\nc\vdifr[7]{\jcoo#5
\xone=\vartwo
\jcoo#4
\advance\xone by -\vartwo
\advance\xone by 3
\icoo#3
\yone=\varone
\advance\yone by \varthr
\icoo#2
\advance\varone by \varthr
\advance\yone by -\varone
\advance\yone by 3
\xdifr{#1}{#3}{#4}{#6}{#7}{\xone}{\yone}
}

\nc\hcurlyline[3]{
\xone=#1
\yone=#2
\xtwo=#3
\xthr=\xtwo \divide\xthr by 10
\xfou=\xone 
\advance\xfou by \xtwo
\put(\xfou,\yone){\oval(5,5)[tl]}
\put(\xone,\yone){\oval(5,5)[tr]}
\advance\xone by 5
\multiput(\xone,\yone)(10,0){\xthr}{\oval(5,5)[b]}
\advance\xthr by -1
\advance\xone by 5
\multiput(\xone,\yone)(10,0){\xthr}{\oval(5,5)[t]}
}

\nc\hcurlyxline[3]{
\hcurlyline{#1}{#2}{#3}
\divide \xtwo by 2
\advance \xtwo by \xone
\advance \xtwo by -30
\advance \yone by -5
\put(\xtwo,\yone){\line(4,1){40}}
\advance \yone by 10
\put(\xtwo,\yone){\line(4,-1){40}}
}

\nc\vcurlyline[3]{
\xone=#2
\yone=#1
\xtwo=#3
\xthr=\xtwo \divide\xthr by 10
\xfou=\xone 
\advance\xfou by \xtwo
\put(\yone,\xfou){\oval(5,5)[br]}
\put(\yone,\xone){\oval(5,5)[tr]}
\advance\xone by 5
\multiput(\yone,\xone)(0,10){\xthr}{\oval(5,5)[l]}
\advance\xthr by -1
\advance\xone by 5
\multiput(\yone,\xone)(0,10){\xthr}{\oval(5,5)[r]}
}

\nc\vcurlyxline[3]{
\vcurlyline{#1}{#2}{#3}
\divide \xtwo by 2
\advance \xtwo by \xone
\advance \xtwo by -30
\advance \yone by -5
\put(\yone,\xtwo){\line(1,4){10}}
\advance \yone by 10
\put(\yone,\xtwo){\line(-1,4){10}}
}

\nc\rcurline{
\advance\varone by \varthr
\advance\varone by -2
\vcurlyline{\varone}{\vartwo}{\varfou}
\advance\varone by -4
\put(\varone,\vartwo){\makebox(0,\varfou)[r]{$\LabeL$}}
}

\nc\lcurline{
\advance\varone by 2
\vcurlyline{\varone}{\vartwo}{\varfou}
\advance\varone by 4
\put(\varone,\vartwo){\makebox(0,\varfou)[l]{$\LabeL$}}
}

\nc\bcurline{
\advance\vartwo by 2
\hcurlyline{\varone}{\vartwo}{\varthr}
\advance\vartwo by 4
\put(\varone,\vartwo){\makebox(\varthr,0)[b]{$\LAbeL$}}
}

\nc\tcurline{
\advance\vartwo by \varfou
\advance\vartwo by -2
\hcurlyline{\varone}{\vartwo}{\varthr}
\advance\vartwo by -4
\put(\varone,\vartwo){\makebox(\varthr,0)[t]{$\LAbeL$}}
}

\nc\makeLAbeL[2]{
\def\LAbeL{
\Callet{#1}^{\arrbi{#2}}}
}

\nc\makeLabeL[2]{
\def\LabeL{
\Callet{#1}^{\arubi{#2}}}
}

\nc\crli[5]{\icoo{#1} \jcoo{#2}
\if#5r
\makeLabeL{#3}{#4}\rcurline\fi
\if#5l
\makeLabeL{#3}{#4}\lcurline\fi
\if#5t
\makeLAbeL{#3}{#4}\tcurline\fi
\if#5b
\makeLAbeL{#3}{#4}\bcurline\fi
}

\nc\rxcurline{
\advance\varone by \varthr
\advance\varone by -2
\vcurlyxline{\varone}{\vartwo}{\varfou}
}

\nc\lxcurline{
\advance\varone by 2
\vcurlyxline{\varone}{\vartwo}{\varfou}
}

\nc\bxcurline{
\advance\vartwo by 2
\hcurlyxline{\varone}{\vartwo}{\varthr}
}

\nc\txcurline{
\advance\vartwo by \varfou
\advance\vartwo by -2
\hcurlyxline{\varone}{\vartwo}{\varthr}
}

\nc\xcrli[3]{\icoo{#1} \jcoo{#2}
\if#3r\rxcurline\fi
\if#3l\lxcurline\fi
\if#3t\txcurline\fi
\if#3b\bxcurline\fi
}

\nc\Ttrui[4]{
\trui{#1}{#2}{{\bf T}^{\barbi {#3}{#4}}}
}

\nc\Ttruj[4]{
\truj{#1}{#2}{{\bf T}^{\barbi {#3}{#4}}}
}

\nc\Ttrdi[4]{
\trdi{#1}{#2}{{\bf T}^{\barbi {#3}{#4}}}
}

\nc\Ttrdj[4]{
\trdj{#1}{#2}{{\bf T}^{\barbi {#3}{#4}}}
}

\nc\vertarr[3]{
\icoo{#1} \jcoo{#2}\fone=\vartwo\advance
\fone by \varfou
\jcoo{#3}\advance \vartwo by -\fone
\varfou=\vartwo \vartwo=\fone}

\nc\horzarr[3]{\icoo{#2} \fone=\varone
\icoo{#1} \advance\varone by \varthr
\varthr=\fone \advance\varthr by -\varone \jcoo{#3}
}

\nc\fposar[5]{
\if#5t
\horzarr{#1}{#2}{#3}
\fone=\vartwo \advance \fone by \varfou
\advance \varthr by -4
\advance \varone by 2
\setra\varone\fone\varthr{#4}
\else\if#5n
\horzarr{#1}{#2}{#3}
\fone=\vartwo \advance \fone by \varfou
\advance \varthr by -4
\setra\varone\fone\varthr{#4}
\else\if#5l
\vertarr{#1}{#2}{#3}
\advance \vartwo by 2
\advance\varfou by -4
\setua\varone\vartwo\varfou{#4}
\else\if#5w
\vertarr{#1}{#2}{#3}
\advance\varfou by -4
\setua\varone\vartwo\varfou{#4}
\else\if#5b
\horzarr{#1}{#2}{#3}
\advance \varone by 2
\advance \varthr by -4
\setra\varone\vartwo\varthr{#4}
\else\if#5s
\horzarr{#1}{#2}{#3}
\advance \varone by 4
\advance \varthr by -4
\setra\varone\vartwo\varthr{#4}
\else\if#5r
\vertarr{#1}{#2}{#3}
\eone=\varone \advance \eone by \varthr
\fone=\vartwo \advance \fone by 2
\advance \varfou by -4
\setua\eone\fone\varfou{#4}
\else\if#5e
\vertarr{#1}{#2}{#3}
\eone=\varone \advance \eone by \varthr
\fone=\vartwo \advance \fone by 4
\advance \varfou by -4
\setua\eone\fone\varfou{#4}
\fi\fi\fi\fi\fi\fi\fi\fi
}

\nc\curmon[5]{\coosq{#1}{#2}{#3}{#4}
\gone=#5
\eone=\varone \advance \eone by 3 
\ethr=\eone \advance\ethr by \gone 
\advance \varthr by -\gone 
\advance \varthr by -3 
\fone=\vartwo
\advance\fone by \varfou
\put(\eone,\fone){\line(1,0){\gone}}
\advance\fone by 2 
\put(\eone,\fone){\oval(4,4)[br]}
\advance\fone by -4 
\put(\eone,\fone){\oval(4,4)[tr]}
\advance\vartwo by 14
\put(\ethr,\vartwo){\oval(4,4)[tl]}
\advance\vartwo by 4
\put(\ethr,\vartwo){\oval(4,4)[bl]}
\advance\vartwo by -2  
\advance\varfou by -16
\divide\varfou by 2
\advance\vartwo by \varfou
\multiply\varfou by 2
\put(\ethr,\vartwo){\oval(\varthr,\varfou)[r]}
}

\nc\tri[1]{\begin{triviality}
\label{#1}}
\nc\cnj[1]{\begin{conj}
\label{#1}}
\nc\prt[1]{\begin{proto}
\label{#1}}
\nc\lem[1]{\begin{lemma}
\label{#1}}
\nc\pro[1]{\begin{prop}
\label{#1}}
\nc\thm[1]{\begin{theorem}
\label{#1}}
\nc\cor[1]{\begin{corollary}
\label{#1}}
\nc\dfn[1]{\begin{defin}
\label{#1}}
\nc\sthm[1]{\begin{subth}
\label{#1}}
\nc\exm[1]{\begin{example}
\label{#1}
\begin{em}}
\nc\plm[1]{\begin{prblm}
\label{#1}
\begin{em}}
\nc\rmk[1]{\begin{remark}
\label{#1}
\begin{em}}
\nc\ntn[1]{\begin{notation}
\label{#1}
\begin{em}}
\nc\cau[1]{\begin{caution}
\label{#1}
\begin{em}}
\nc\imn[1]{\begin{importnota}
\label{#1}
\begin{em}}
\nc\cax[1]{\begin{cauex}
\label{#1}
\begin{em}}
\nc\con[1]{\begin{construction}
\label{#1}
\begin{em}}
\nc\ssthm[1]{\begin{ssubth}
\label{#1}
\begin{em}}
\nc\cnc[1]{\begin{conclusion}
\label{#1}
\begin{em}}

\nc\elem{\end{lemma}}
\nc\ecnj{\end{conj}}
\nc\eprt{\end{proto}}
\nc\ethm{\end{theorem}}
\nc\ecor{\end{corollary}}
\nc\edfn{\end{defin}}
\nc\esthm{\end{subth}}
\nc\epro{\end{prop}}
\nc\etri{\end{triviality}}
\nc\eexm{\end{em}
\end{example}}
\nc\ermk{\end{em}
\end{remark}}
\nc\eplm{\end{em}
\end{prblm}}
\nc\ecau{\end{em}
\end{caution}}
\nc\ecax{\end{em}
\end{cauex}}
\nc\eimn{\end{em}
\end{importnota}}
\nc\entn{\end{em}
\end{notation}}
\nc\econ{\end{em}
\end{construction}}
\nc\ecnc{\end{em}
\end{conclusion}}
\nc\essthm{\end{em}
\end{ssubth}}

\nc\nin{\noindent}
\nc\sea{\searrow}
\nc\swa{\swarrow}
\nc\nwa{\nwarrow}
\nc\nea{\nearrow}
\nc\da{\hbox{
\begin{picture}(0,20)
\put(0,20){\vector(0,-1){20}}
\put(-7,0){\makebox(7,20){$\theta$}}
\end{picture}}}

\def\kth{{\it K--}theory}
\def\kthi{{\it K--theory for triangulated categories I}}
\def\kthii{{\it K--theory for triangulated categories II}}
\def\kthiii{{\it K--theory for triangulated categories III}}
\nc\tst{{\it t--}structure}

\nc\pf{\par \noindent {\bf Proof.}\ \ }
\nc\car{\cal R}

\nc\la{\longrightarrow}

\begin{document}

\slugline{AJM}{2}{3}{495--589}{September}{1998}{004}
\setcounter{page}{495}


\title{K--Theory for Triangulated Categories III(a): The Theorem 
of the Heart\thanks{Received April 30, 1997; 
accepted for publication September 20, 1997.}}
\author{Amnon Neeman\thanks{Department of Mathematics,
University of Virginia, Charlottesville, Virginia 22903 USA
(an3r@virginia.edu).}}

\pagestyle{myheadings}
\thispagestyle{plain}
\markboth{AMNON NEEMAN}{K--THEORY FOR TRIANGULATED CATEGORIES III(A)}

\maketitle


\setcounter{section}{-1}

\section{Introduction}

\siv500000000
\sjv500000000

\SB\boxone{
\fbx11\ct aa}
\siv250000000
\sjv520000000

\SB\boxtwo{
\fbx21\ct aa \diff2211{}}

This is the fourth installment of a series.
The main point of the entire series
is the following: given a triangulated category $\ct$, it is
possible to attach to it a \kth\ space. Its delooping
will be denoted $\PT\boxone$. Note that, starting
with the present article, we no longer wish to consider
the construction without the differentials. In the earlier
parts of this series, we considered two simplicial
sets, namely
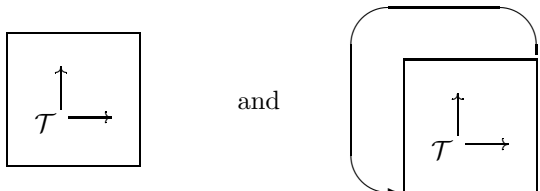
\BDI(0,70)
\UB(-70,35)\boxone  \put(0,35){\makebox(0,0){and}}
\UB(70,35)\boxtwo
\EDI
 From now on, we wish to consider only $\PT\boxtwo$;
all the simplicial sets will be the ones with
coherent differentials. We will feel free to
omit the differentials in the symbol for the 
simplicial set. The simplicial sets
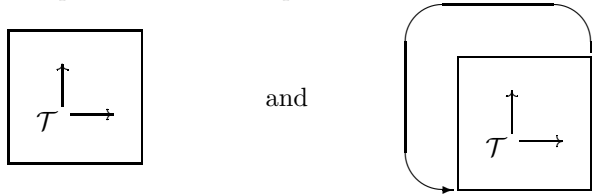
\BDI(0,70)
\UB(-80,35)\boxone  \put(0,35){\makebox(0,0){and}}
\UB(80,35)\boxtwo
\EDI
are henceforth to be viewed as identical. 

The reason for this is that, after Section~\ref{S10},
we have nothing more to say about the construction
without differentials. In Section~\ref{S10}, we
proved some significant facts about the simplicial set without
the differentials. The reader is referred to the introduction
of \kthii\ \ for more detail. Anyway, from Section~\ref{S11}
on, all our simplicial sets come with coherent
differentials.
 
The
key theorem of this series of articles is

\siv050000000
\sjv500000000
\SB\boxtwo{
\fbx21{\ca}aa 
}

\medskip

\nin 
{\bf Strong Theorem~\ref{T7.1}}. {\em Let ${\cal T}$ 
be a small triangulated category with a non--degener-ate
$t$--structure.  Let $\ca$ be the heart of
the {\it t}--structure.
With the simplicial set 
\BC
\usebox\boxone
\EC
defined appropriately, the natural map
\BDI(0,50)
\UB(-80,25)\boxtwo  \Ri(0,25,80){}
\UB(80,25)\boxone
\EDI
induces a homotopy equivalence.}

\medskip

\nin
In this theorem, $\PT\boxtwo$ is homotopy equivalent
to Quillen's $Q$--construction of the abelian category
$\ca$. The precise definition of $\PT\boxone$ is
a somewhat delicate point, discussed in some
detail in the introduction to \kthi.
Delicate points aside, a consequence of this
theorem is that, given two non--degenerate {\it t}--structures
on the same triangulated category, the two
hearts have isomorphic {\it K}--theories.
We will prove
that many abelian categories have isomorphic {\it K--}theories.

The proof of Strong Theorem~\ref{T7.1} is the bulk
of \kthiii. But in fact, we will be proving more. Let me
state for the reader another theorem, which will
follow from the same proof. We begin with definitions.

\dfn{D30.0.1} 
Let $\cal E$ be an exact category. A sequence
\[
x\la y \la z
\]
in $\cal E$ is called {\em exact at $y$} if 
\begin{description}
\item
\sthm{D30.0.1.1}
The map $x\la y$ factors as 
\[x\la y'\la y,\] 
with
$x\la y'$ an admissible epi, and $y'\la y$ an admissible
mono.
\esthm
\item
\sthm{D30.0.1.2}
The map $y\la z$ factors as 
\[y\la y''\la z,\] 
with
$y\la y''$ an admissible epi, and $y''\la z$ an admissible
mono.
\esthm
\item
\sthm{D30.0.1.3}
$y'\la y\la y''$ is an admissible short exact sequence.
\esthm
\end{description}
\edfn

Let $Gr^b({\cal E})$ be the category of bounded, 
$\Bbb Z$--graded
objects in $\cal E$. Let $\Sigma:
Gr^b({\cal E})\la 
Gr^b({\cal E})$
be the shift map. Next, we define a simplicial set.

\siv800000000
\sjv500000000
\SB\boxthr{
\fbx11{Gr^b({\cal E})}aa
}
\siv500000000
\sjv500000000
\SB\boxfou{
\fbx11{{\cal E}}aa
}
\dfn{D30.0.2}
The bisimplicial set
\BC
\usebox\boxthr
\EC
is defined as follows. A $(p,q)$--simplex is a diagram
in $Gr^b({\cal E})$
\siv{13}00000000
\sjv{13}00000000
\BP
\Sbx11aa{X_\Q}0q0p
\EP
together with a coherent differential $X_{pq}\la\Sigma X_{00}$.
The condition is that, for every $0\leq i\leq i'\leq p$,\ \ 
$0\leq j\leq j'\leq q$,
the sequence
\[
\Sigma^{-1}X_{i'j'}\la
X_{ij}\la X_{i'j}\oplus X_{ij'} \la X_{i'j'}\la
\Sigma X_{ij}
\]
gives, in each degree, an exact sequence in $\cal E$.
\edfn

\nin
The proof in this article, which will establish Theorem~\ref{T7.1},
will also prove the following fact.

\medskip

\nin 
{\bf Strong Theorem~\ref{T4.8}}. {\em Let ${\cal E}$ 
be a small exact category. The natural map
\BDI(0,50)
\UB(-80,25)\boxfou  \Ri(0,25,80){}
\UB(95,25)\boxthr
\EDI
induces a homotopy equivalence.}

\medskip

\nin
Here, $\PT\boxfou$ is homotopy equivalent
to Quillen's {\it Q}--construction on
the exact category $\cal E$.

Now for a review of the earlier parts of this series.
\kthi\ \ contains a proof of the special case of
Theorem~\ref{T4.8}, where $\cal E$ is an abelian category.
\kthii\ \ contains a proof of the special case of Theorem~\ref{T7.1},
where $\ct$ is $D^b(\ca)$, the bounded
derived category of an abelian
category $\ca$, and the \tst\ is the standard one. In a very precise
sense, the current article is better. It proves
the sharpest and most general results of the series.

This raises the question: what is the point of the
earlier articles? Let me try to answer it briefly.

First of all, the three proofs are all different. They
look at quite different chains of intermediate simplicial
sets. Let us agree that the current theory is unsatisfactory,
and that it is to be hoped that there will, some day, be a
simpler and more general treatment. Then surely different
arguments are of interest. It is unclear which
will lead to the better generalisations and
simplifications.

The second reason that the earlier articles are of interest,
is that they are simpler. The theorems they prove are
not optimal; but there is virtue in seeing first
a simple proof of a less general statement. The simpler
argument is also easier to motivate. Finally, it lends itself
more to careful study of alternatives, such as the
construction without the differentials. In
\kthi\ \ and {\it II,} we do more than just give 
the proofs of special cases of the strong theorems
stated above. We explain how and why the proofs work.

In the 
introduction to \kthi, I divided up the readers of any piece of mathematics
into three broad groups, listed in order of probable size:

\begin{description}

\item[Group 1:]  The people who want a rough idea of the contents of the 
article, and at the very most a sketch of the proofs in an easy special case. 

\item[Group 2:]  The people who want to check the result, because they might 
consider using it in their own work. 

\item[Group 3:]  The people reading the article because they might work on the 
problem themselves.

\end{description}

\nin
The first two parts of this series, \kthi{\it(A)\/} and {\it I(B)}, 
were intended for a Group 1
audience. The third, \kthii, is emphatically for the benefit of
Group 3. The present part is primarily for Group~2.

In \kthi{\it(A)},
we introduce the definitions and notation (this takes us some
88 pages). Then in \kthi{\it(B)},
we give the simplest proof of the simplest version of
our theorem. All the readers of subsequent parts are assumed to be
familiar with the notation. So you should have read at least 
\kthi{\it(A)},
if you proceed beyond this word. In fact, it is highly advisable
to have skimmed through the rest of \kthi. There is a little more notation
introduced in the last two sections, but even more relevant, there
is a relatively gentle introduction to the way the proofs work, and the type
of simplicial sets one constructs. 

The first section of this part, Section~\ref{S3.0}, is again quite soft.
There are two types of homotopy that I know, for the simplicial sets that
come up in triangulated \kth. The first type is the trivial homotopies.
These are the triangulated analogues of contractions to an initial or
a terminal object. The second type of homotopy is the non--trivial 
homotopies. And one of the key features of this theory is that there is 
really only one of the non--trivial homotopies.

In Section~\ref{S3.0}, we make this very precise, showing with explicit
examples how to reduce a typical non--trivial homotopy 
in this theory to a blueprint.
This section is really a must for anyone who reads beyond \kthi. Although
not compulsory, it is strongly recommended that the reader also look
at Section~\ref{S9}, the first section of \kthii. Although \kthii\ was
written with a Group 3 audience in mind, Section~\ref{S9} is only at
Group $2\frac12$ level. It is quite soft. It discusses, in a general way,
the type of simplicial sets and homotopies that come up in the proof, and
it also discusses why the various homotopies are well defined. 
Section~\ref{S3.0} of this part, being written for Group 2, 
focuses on the non--trivial 
homotopy. It turns out that one can give a very satisfactory treatment of
it, and explain why checking that it is well--defined can be reduced
to verifying it on a blueprint. Section~\ref{S9} is for Group 3, and 
therefore it tends to focus on potential problems. It turns out that one
of the so--called trivial homotopies is less trivial than it seems to be
at first sight. This homotopy is the truncation.

Having familiarised himself with the notation, and the type of argument 
used to show that the homotopies are well--defined, the reader will
discover that he has read almost two thirds of this \kthiii. The remaining
sections, Sections \ref{S3.2} and \ref{S3.3}, contain the proof of the
main theorem of the article. The proof we give here is very businesslike.
It demonstrates the best theorem I have about the \kth\ of 
triangulated categories,
and does so as directly as possible. It is very difficult to say much about
the proof, that would be in any way instructive. It is a sequence of maps
and homotopies, that get us where we want to be.

This completes the discussion of all the theorems in the article.
There are also three conjectures. They may be found in
Appendices~\ref{App2}, \ref{App3} and \ref{App4}. Appendix~\ref{App2}
explores the natural map from Waldhausen's \kth\ of a Waldhausen
category, to the triangulated \kth\ of the associated triangulated
category. The map need not be a homotopy equivalence. But
there is an intermediate space, with a description similar to
Waldhausen's, whose homotopy type is conjecturally the same
as triangulated \kth. In Appendix~\ref{App2}, I state the conjecture,
and show that if true, it implies that for any exact category
$\cal E$, the \kth\ of $D^b({\cal E})$ agrees with Quillen's
\kth\ of $\cal E$.

Appendix~\ref{App3} states a conjecture, generalising Quillen's
localisation theorem. The conjecture is straightforward enough
to state. In the appendix, I also explain my attempts (so far
quite unsuccessful) to generalise Quillen's proof.

Finally, Appendix~\ref{App4} gives a vaguely--stated conjecture,
generalising Quillen's devissage theorem to triangulated
categories. 

I tried to keep the conjectural appendices very short. Appendix~\ref{App2},
on the relation between Waldhausen's \kth\ and triangulated \kth,
is the longest. It is the problem I thought about the most.
At one point I even thought I had a proof. There is an error
in the manuscript I wrote at the time. But even after finding
the error, I was under the impression it was fixable. I never
checked this carefully, and in the five years that have since
passed, I have forgotten all the subtle points. There still
exists a 300 page document, {\it K--theory for triangulated
categories IV}, with an outline of what the proof might(?)
look like.

I was never under the impression, that I knew how to prove
the localisation conjecture of Appendix~\ref{App3}. All the
appendix offers is the statement of the conjecture, and
the statements that would need to be proved, to
generalise Quillen's argument from abelian to triangulated
categories. Finally, in Appendix~\ref{App4} there is
not even a clearly formulated conjecture; all there is
is an idea.

One historical note. The proofs of Theorems \ref{T4.8} and \ref{T7.1}
in the weak, special cases date to September, 1988. However,
I did not have 
the proof given here until the following spring, some six months
later.

\section{Why the Non Trivial Homotopy Is  Always Well--Defined}
\label{S3.0}

In \kthi, the reader saw a very simple proof of a weak version of our 
main theorem.
But since we introduced the notation as the proof progressed, there was no
possibility of giving a general discussion of our homotopies. It is now time 
to rectify this shortcoming.

Our homotopies fall into two groups: the trivial and the non--trivial. The
trivial homotopies are contractions to an initial or terminal object. Our 
notation for such homotopies was something like

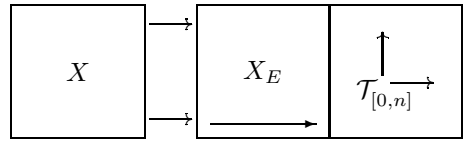
\siv525500000
\sjv500000000

\BP
\Fbx11X  \br21a \Fhbx31{X_E}r \bx41aa0n
\EP

\nin
Sometimes the matter was a little more delicate. For instance, if we consider
the simplicial set

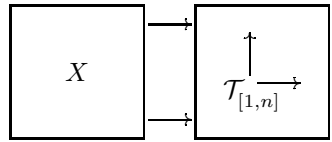
\siv525000000
\sjv500000000

\BP
\Fbx11X \br21a  \bx31aa1n
\EP

\siv500000000

\SB\boxone{
\Fbx11X
}

\SB\boxtwo{
\bx11ha0n}

\nin
where $\PT\boxone$ lies in the simplicial set $\PT\boxtwo,$ there is a
homotopy

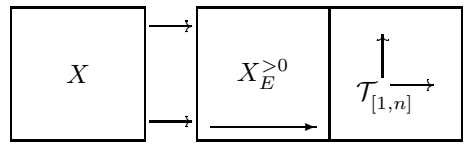
\siv525500000
\sjv500000000

\BP
\Fbx11X  \br21a \Fhbx31{X_E^{>0}}r \bx41aa1n
\EP

\nin
which is very like the contraction to the initial object, except it need not
be a contraction. This is because of the non--uniqueness of the differential 
from the truncation.
But those are subtle points for Group 3 readers. I give a very thorough 
discussion of the non--triviality of the so--called trivial homotopy, in
\kthii. A reader with an interest in this point is referred to Section~\ref{S9}.
 The present section occupies itself with the triviality of the non--trivial
homotopy. As befits the type of material we present our Group 2
audience, it is possible to give a very satisfactory general discussion,
of why the cells in the non--trivial homotopy always do as they should.

We remind the reader of our non--trivial homotopy. It is the homotopy whose 
shorthand came to something like

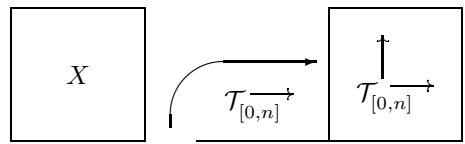
\BP
\Fbx11X   \htp31a0nr \bx41aa0n
\EP

\nin
In other words, the cells of the homotopy would be

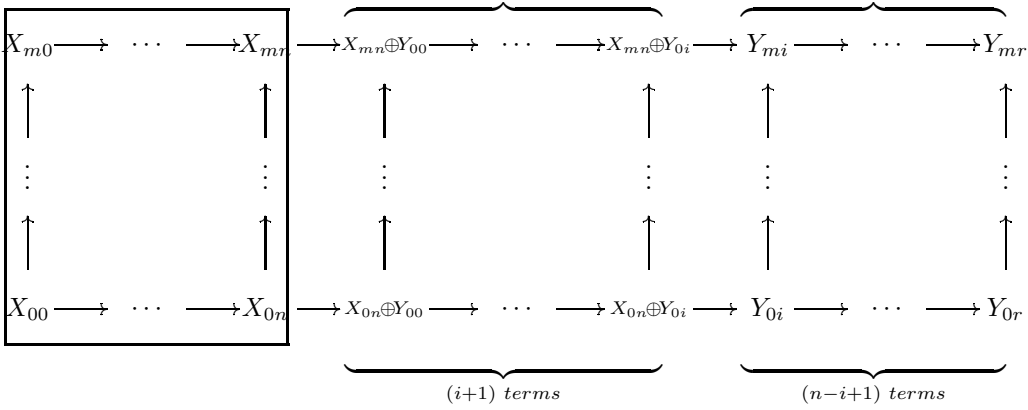
\siv{11}2{13}2{11}0000
\sjv4{13}0000000

\A10n
\A30i
\A5ir
\B20m

\SB\boxone{
\Sbxp12aa{X_\Q} \Sbrsmall22a  
\Sbxpsmall32aa{\scriptstyle X_{\I n}\!\oplus\! Y_{0\J}} 
\Sbrsmall42a \Sbxp52aa{Y_\Q}
\ovrbrc32{}0 \ovrbrc52{}0 \undrbrc32{(i+1)\,\,terms}0  
\undrbrc52{(n-i+1)\,\,terms}0
\fmbx1122
}

\BDI(0,185)
\UB(0,85)\boxone
\EDI

\siv500000000
\sjv500000000
\SB\boxone{
\Fbx11\tqu}

\SB\boxtwo{
\fbx11\ct aa}    

\nin
In practice, we almost never apply this homotopy as written above.
Usually, there are restrictions on the objects allowed, and
on the morphisms permitted among them.
Following the conventions of Section~\ref{S9}, we even have a notation which
reminds us that the objects and morphisms are restricted,
without specifying the restrictions.
The 
symbol $\PT\boxone$ stands for the
simplicial set, in which the exact 
subcategory (the subscript) is unspecified, and the horizontal and vertical 
morphisms are also left ambiguous (hence the question marks). In this section,
we propose to give general arguments for why the non--trivial homotopy above,
and its various close cousins whom we have met in \kthi, are all
well--defined. We need to establish that the purported cells
of the homotopy are genuine simplices. For this, it
is simpler to assume there are no restrictions on morphisms
and objects. In other words, it is easiest to handle the case where
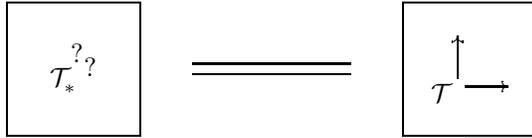
\BDI(0,50)
\UB(-75,25)\boxone \Rq(0,25,60)  \UB(75,25)\boxtwo
\EDI
Thus, we will really not be proving any theorems about $\PT\boxone$.
We will be showing that certain homotopies are well--defined in
$\PT\boxtwo$. Specifically, we will show that the 
homotopy

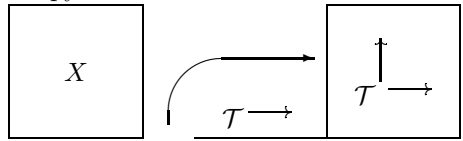
\siv525500000
\sjv500000000

\BP
\Fbx11X   \fhtpr3311 \fli31\ct a{}b \fbx41\ct aa
\EP 

\nin
together with its various analogues, 
are well defined. The particular subcategory 
of $\ct$ that we are
dealing with, and the precise nature of 
the morphisms, play absolutely no role 
in the main theorem 
of this section; they are, after all, the source of the subtleties. If we
avoid the subtleties like the plague, we can state and prove an honest theorem.
After doing that, we will return to showing how this theorem can be applied,
even in the subtle situation where the simplicial set is $\PT\boxone$, that is
where the objects and the morphisms are restricted.

\siv{19}15151500
\sjv525252500

\SB\boxeig{
                                     \trui570
                                      \bu56a
         \br25a  \fbx35\ct aa  \br45a    \fbx55\ct aa    \br65a   \fbx75\ct aa
                 \bu34a                \bu54a               
         \br23a   \fbx33\ct aa  \br43a    \fbx53\ct aa    \br63a   \fbx73\ct aa
                  \bu32a               \bu52a                   \bu72a
\trui110 \br21a   \fbx31\ct aa  \br41a    \fbx51\ct aa    \br61a   \Fbx71X
}

\thm{well defined homotopies} Consider the simplicial set

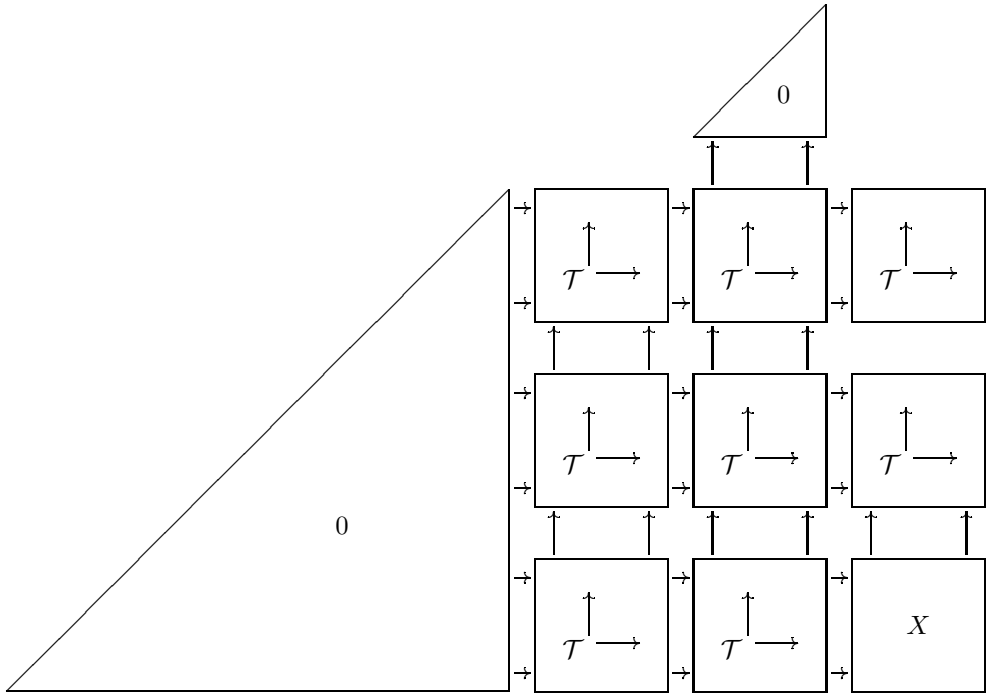
\BDI(0,260)
\UB(0,130)\boxeig
\EDI

\nin
Then on it, the homotopy whose symbol would be

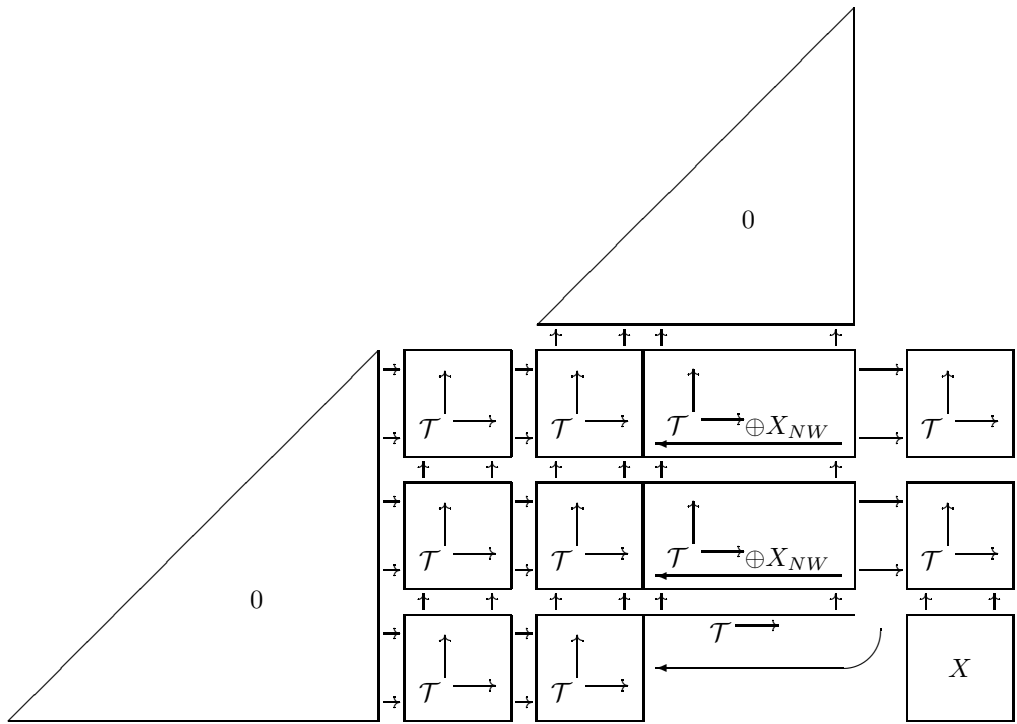
\siv{14}14148240
\sjv414141{12}00

\SB\boxten{
                                     \truj670
                                      \bu56a  \bu66a
         \br25a  \fbx35\ct aa  \br45a    \fbx55\ct aa 
\fhbx65\ct aa{ \oplus X_{NW}}l   \br75a   \fbx85\ct aa
                 \bu34a                \bu54a        \bu64a       
         \br23a   \fbx33\ct aa  \br43a    \fbx53\ct aa 
\fhbx63\ct aa{\oplus X_{NW}} l   \br73a   \fbx83\ct aa
                  \bu32a               \bu52a        \bu62a                \bu82a
\trui110 \br21a   \fbx31\ct aa  \br41a    \fbx51\ct aa  
\fhtpl6611 \fli61\ct a{} t            \Fbx81X
}

\BDI(0,270)
\UB(0,135)\boxten
\EDI

\nin
is well defined; its cells are all simplices for the simplicial set.
\ethm

{\it Proof.}
In this proof, when we say that a square is Mayer--Vietoris,
this will mean that, with an obvious choice of a
(coherent) differential,
the square ``folds'' to give a semi--triangle.
We want to prove all the squares in some large diagram
Mayer--Vietoris in this sense.
Let us begin by labeling all the regions of this diagram. 
We 
rewrite it as

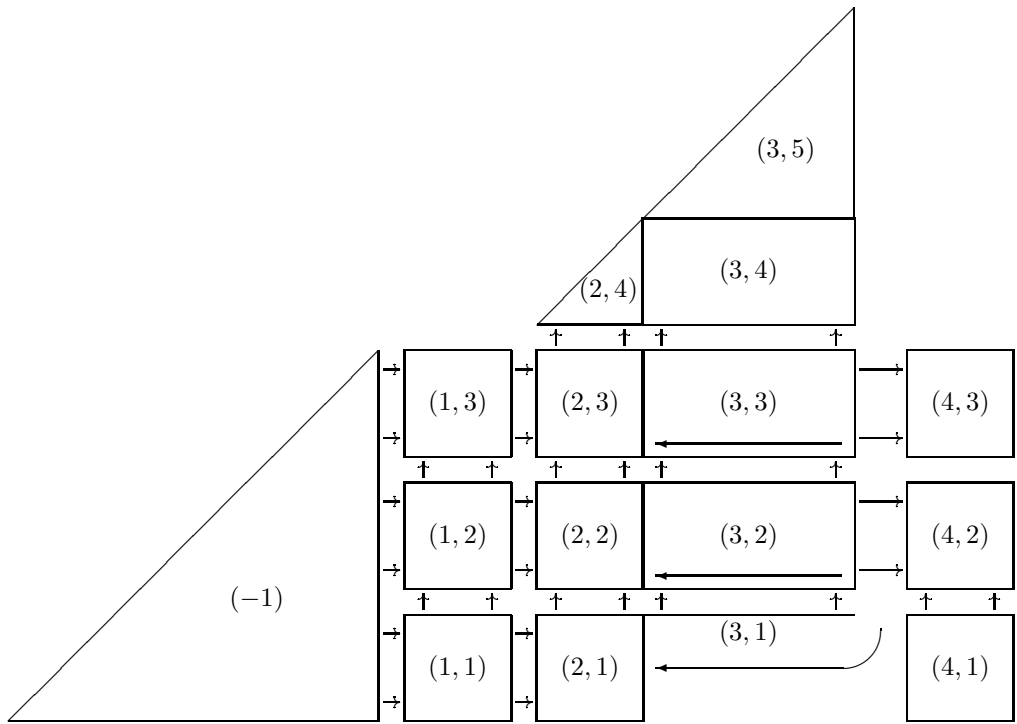
\siv{14}14148240
\sjv414141480

\SB\boxten{
                                     \truj68{(3,5)}
                                      \truj57{(2,4)} \Fbx67{(3,4)}
                                      \bu56a  \bu66a
         \br25a  \Fbx35{(1,3)} \br45a    \Fbx55{(2,3)}\Fhbx65{(3,3)}l   
\br75a   \Fbx85{(4,3)}
                 \bu34a                \bu54a        \bu64a       
         \br23a   \Fbx33{(1,2)} \br43a    \Fbx53{(2,2)}\Fhbx63{(3,2)}l  
 \br73a   \Fbx83{(4,2)}
                  \bu32a               \bu52a        \bu62a                \bu82a
\trui11{(-1)} \br21a   \Fbx31{(1,1)} \br41a    \Fbx51{(2,1)} \fhtpl6611 
\Fli61{(3,1)} t            \Fbx81{(4,1)}
}

\BDI(0,270)
\UB(0,135)\boxten
\EDI 

\nin
The labels are to allow us to easily refer to any particular region of the 
homotopy.

What we must show first, is that every square that occurs 
in this gigantic simplex 
at the very least folds to give a semitriangle.
More precisely, we must show that there is a way to
choose the differentials, so that all the squares in the diagram will 
fold to
mapping cones on triangles, or at the
very worst direct summands of mapping cones.
 Furthermore, we will do so in such a way that all the maps,
differentials included, will be given maps; they are matrices in maps defining 
the starting cell of the homotopy.

Any square in the diagram which does not meet Column 3 is a square in

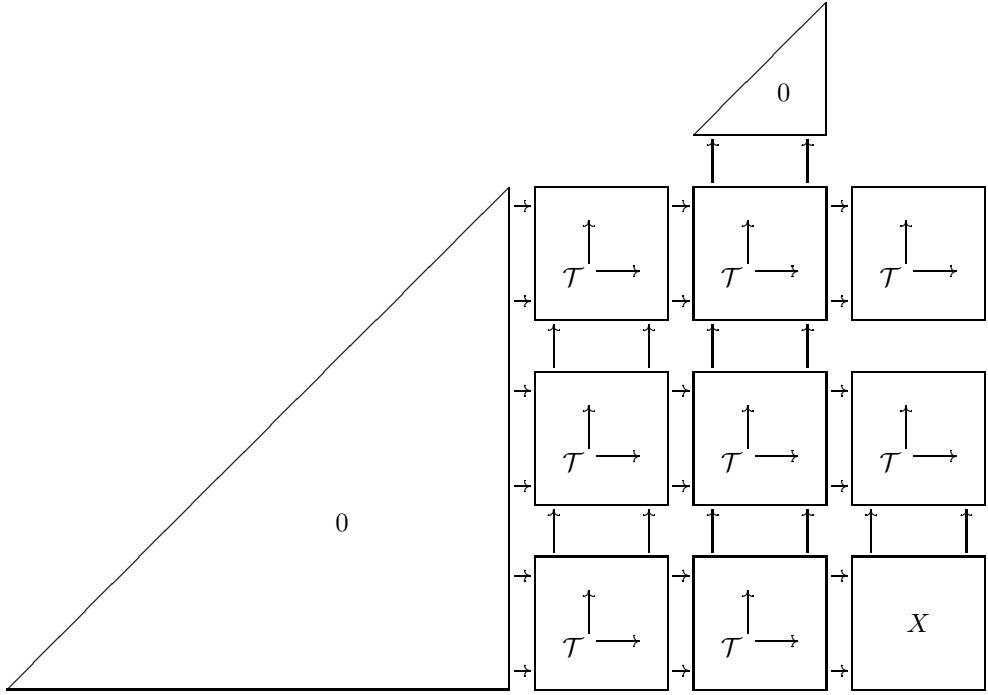
\siv{19}15151500
\sjv525252500

\SB\boxone{
                                     \trui570
                                      \bu56a
         \br25a  \fbx35\ct aa  \br45a    \fbx55\ct aa    \br65a   \fbx75\ct aa
                 \bu34a                \bu54a               
         \br23a   \fbx33\ct aa  \br43a    \fbx53\ct aa    \br63a   \fbx73\ct aa
                  \bu32a               \bu52a                   \bu72a
\trui110 \br21a   \fbx31\ct aa  \br41a    \fbx51\ct aa    \br61a   \Fbx71X
}

\BDI(0,260)
\UB(0,130)\boxone
\EDI

\nin
and therefore it is Mayer--Vietoris because by hypothesis, the homotopy started
with a simplex. Thus we need only concern ourselves with squares which meet
 Column 3.

Any square meeting Column 3 must also meet Column (--1) (i.e. the triangle on the
right), Column 1, Column 2, Column 3 or Column 4. Columns (--1),  1 and  2 behave
identically as far as the following argument goes. We will therefore treat only
the cases  of Column 1, Column 3 and Column 4.

\sthm{cols13} {\bf Case of Column 1.} Suppose we are trying 
to prove that a candidate square in  the
union of Columns 1 and 3 can be chosen canonically to be a Mayer--Vietoris 
square. Then there are  3 cases to consider. Either the square is contained
in the first row, or it does not meet the first row, or it meets the first row
as well as some other row. We discuss these cases separately.
\esthm

\ssthm{cols13rows11} Suppose the square is 
entirely contained in the first row.
Then it consists of taking a column in the (1,1) 
gridbox in the labeled diagram,
and pairing it with a column in the (3,1) gridbox. The result is a square

\siv337000000
\sjv333000000

\BP
\postr13{A_{i'j}}c \Sposar23ar \postr33{B_{tj'}\oplus X_{i'0}}c
\Sposar12au                    \Sposar32au
\postr11{A_{ij}}c \Sposar21ar \postr31{B_{tj'}\oplus X_{i0}}c
\EP

\nin
where the $A$'s are a column in the (1,1) box, the $X$'s come from the west
column of the (4,1) box and the $B_{tj'}$ is in the north face of the (3,1)
box. But

\siv333000000
\sjv333000000

\SB\boxsix{
\postr13{A_{i'j}}c \Sposar23ar \postr33{X_{i'0}}c
\Sposar12au                    \Sposar32au
\postr11{A_{ij}}c \Sposar21ar \postr31{X_{i0}}c
}

\SB\boxsev{
\postr13{0}c \Sposar23ar \postr33{B_{tj'}}c
\Sposar12au                    \Sposar32au
\postr11{0}c \Sposar21ar \postr31{B_{tj'}}c
}

\BC
\usebox\boxsix
\EC

\nin
is a square in the simplicial set

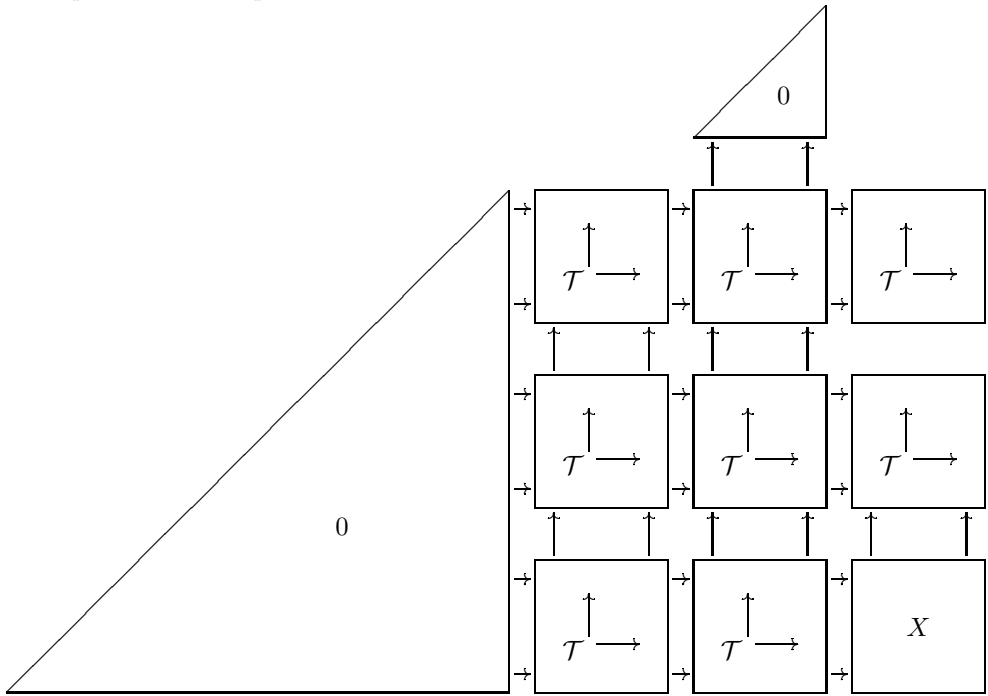
\BDI(0,260)
\UB(0,130)\boxone
\EDI

\nin
and is hence automatically Mayer--Vietoris.  The square

\siv337000000
\sjv333000000

\BP
\postr13{A_{i'j}}c \Sposar23ar \postr33{B_{tj'}\oplus X_{i'0}}c
\Sposar12au                    \Sposar32au
\postr11{A_{ij}}c \Sposar21ar \postr31{B_{tj'}\oplus X_{i0}}c
\EP

\nin
is isomorphic to the direct sum of

\barr
\PT\boxsix &\qquad{\rm and}\qquad&\PT\boxsev
\earr

\nin
and is therefore also Mayer--Vietoris.   

\essthm

\ssthm{cols13rows22} The next 
possibility is that the square does not meet Row 1
at all. In that case, it must be of the form

\siv337000000
\sjv333000000

\BP
\postr13{A_{i'j}}c \Sposar23ar \postr33{B_{i'j'}\oplus X_{NW}}c
\Sposar12au                    \Sposar32au
\postr11{A_{ij}}c \Sposar21ar \postr31{B_{ij'}\oplus X_{NW}}c
\EP

\nin
and can be expressed as a sum of

\siv333000000
\sjv333000000

\SB\boxsix{
\postr13{A_{i'j}}c \Sposar23ar \postr33{B_{i'j'}}c
\Sposar12au                    \Sposar32au
\postr11{A_{ij}}c \Sposar21ar \postr31{B_{ij'}}c
}

\SB\boxsev{
\postr13{0}c \Sposar23ar \postr33{X_{NW}}c
\Sposar12au                    \Sposar32au
\postr11{0}c \Sposar21ar \postr31{X_{NW}}c
}

\barr
\PT\boxsix &\qquad{\rm and}\qquad&\PT\boxsev
\earr

\nin
where, once again, each of the summands is trivially Mayer--Vietoris.

\essthm
\ssthm{cols13rows12}
The first subtle case, is when both the first row and some other are allowed 
to occur. In this case, our square takes the form

\siv337000000
\sjv333000000

\BP
\postr13{C_{i'j}}c \Sposar23ar \postr33{D_{i'j'}\oplus X_{t0}}c
\Sposar12au                    \Sposar32au
\postr11{A_{ij}}c \Sposar21ar \postr31{B_{tj'}\oplus X_{i0}}c
\EP

\nin
where the $A_{ij}$ comes from the (1,1) box, the $B_{tj'}$ from the north of the
(3,1) box, the $X_{i0}$ from the west face of the (4,1) box, the $X_{t0}=X_{NW}$
is the north--west corner of the (4,1) box and the $C_{i'j}$ is in the
$(1,x)$ box and $D_{i'j'}$ in the $(3,x)$ box, where $x=2$ or 3. Although
we have already discussed this in \kthi, we remind the reader of the argument
that shows why this square is Mayer--Vietoris. We have a diagram of squares

\siv333330000
\sjv333330000

\BP
\postr15{C_{i'j}}c \Sposar25ar \postr35{D_{i'j'}}c
\Sposar14au                    \Sposar34au
\postr13{A_{tj}}c \Sposar23ar \postr33{B_{tj'}}c \Sposar43ar \postr53{ X_{t0}}c
\Sposar12au                    \Sposar32au                    \Sposar52au
\postr11{A_{ij}}c \Sposar21ar \postr31{B_{ij'}}c \Sposar41ar \postr51{X_{i0}}c
\EP

\nin
This gives us three semi--triangles

\siv337333500
\sjv300000000

\SB\boxtwo{
\postr11{A_{ij}}c \Sposar21ar \postr31{B_{ij'}\oplus A_{tj}}c \Sposar41ar 
\postr51{B_{tj'}}c  \Sposar61ar 
\postr71{\Sigma A_{ij}}c
}

\SB\boxthr{
\postr11{A_{ij}}c \Sposar21ar \postr31{X_{i0}\oplus A_{tj}}c \Sposar41ar 
\postr51{X_{t0}}c  \Sposar61ar 
\postr71{\Sigma A_{ij}}c
}

\SB\boxfou{
\postr11{A_{ij}}c \Sposar21ar \postr31{B_{ij'}\oplus C_{i'j}}c \Sposar41ar 
\postr51{D_{i'j'}}c  \Sposar61ar 
\postr71{\Sigma A_{ij}}c
}

\begin{equation}
\label{tr1}
\PT\boxtwo
\end{equation}

\begin{equation}
\label{tr2}
\PT\boxthr
\end{equation}

\nin
and

\begin{equation}
\label{tr3}
\PT\boxfou
\end{equation}

\nin
But now the square

\siv337000000
\sjv333000000

\BP
\postr13{C_{i'j}}c \Sposar23ar \postr33{D_{i'j'}\oplus X_{t0}}c
\Sposar12au                    \Sposar32au
\postr11{A_{ij}}c \Sposar21ar \postr31{B_{tj'}\oplus X_{i0}}c
\EP

\nin
will be Mayer--Vietoris exactly if the following candidate triangle (\ref{tr4})
 is a semi--triangle

\siv33{9}373300
\sjv300000000

\SB\boxtwo{
\postr11{A_{ij}}c \Sposar21ar \postr31{C_{i'j}\oplus B_{tj'}\oplus X_{i0}}c 
\Sposar41ar 
\postr51{D_{i'j'}\oplus X_{t0}}c  \Sposar61ar 
\postr71{\Sigma A_{ij}}c
}

\begin{equation}
\label{tr4}
\qquad\PT\boxtwo
\end{equation}

\nin
and the real point is that (\ref{tr4}) can be obtained as a 
direct summand of the mapping 
cone on the natural map of 
semitriangles

$$
(\ref{tr1})\longrightarrow (\ref{tr2})\oplus(\ref{tr3}).
$$

\nin
What is more, the other direct summand is a contractible triangle.
This is left as an exercise to the reader, who can also find a discussion in
Remark~\ref{R5.3}.

One remark should be made now. If it so happens that the semi--triangles

\siv337333500
\sjv300000000

\SB\boxtwo{
\postr11{A_{ij}}c \Sposar21ar \postr31{B_{ij'}\oplus A_{tj}}c \Sposar41ar 
\postr51{B_{tj'}}c  \Sposar61ar 
\postr71{\Sigma A_{ij}}c
}

\barr
\PT\boxtwo
\earr

\barr
\PT\boxthr
\earr

\nin
and

\barr
\PT\boxfou
\earr

\nin
are really triangles, or more precisely triangles that
come from three short exact sequences

\siv333373333
\sjv300000000

\SB\boxtwo{
\postr11{0}c \Sposar21ar
\postr31{A_{ij}}c \Sposar41ar \postr51{B_{ij'}\oplus A_{tj}}c \Sposar61ar 
\postr71{B_{tj'}}c  \Sposar81ar 
\postr91{0}c
}

\SB\boxthr{
\postr11{0}c \Sposar21ar
\postr31{A_{ij}}c \Sposar41ar \postr51{X_{i0}\oplus A_{tj}}c \Sposar61ar 
\postr71{X_{t0}}c  \Sposar81ar 
\postr91{0}c
}

\SB\boxfou{
\postr11{0}c \Sposar21ar
\postr31{A_{ij}}c \Sposar41ar \postr51{B_{ij'}\oplus C_{i'j}}c \Sposar61ar 
\postr71{D_{i'j'}}c  \Sposar81ar 
\postr91{0}c
}

\barr
\PT\boxtwo
\earr

\barr
\PT\boxthr
\earr

\nin
and

\barr
\PT\boxfou
\earr

\nin
of objects in some heart $\cal C$ of the triangulated category $\ct$, and
where the differentials are the unique ones possible, then the map of triangles

$$
(\ref{tr1})\longrightarrow (\ref{tr2})\oplus(\ref{tr3})
$$

\nin
is really a map of short exact sequences, and it is easy to show that the 
mapping cone 
is indeed a triangle; maps of short exact sequences are good maps of triangles.
Thus the direct summand is also a triangle, and it follows that the differential
in

\siv337000000
\sjv333000000

\BP
\postr13{C_{i'j}}c \Sposar23ar \postr33{D_{i'j'}\oplus X_{t0}}c
\Sposar12au                    \Sposar32au
\postr11{A_{ij}}c \Sposar21ar \postr31{B_{tj'}\oplus X_{i0}}c
\EP

\nin
is the {\em unique} map making the above a genuine triangle.

\essthm

\sthm{cols33} {\bf Case of Column 3.}
The next case to consider, is where we have a square entirely contained in 
Column 3. Once again, this divides into cases, depending on what rows occur.
We distinguish six cases. Case~\ref{cols33rows11} is where only Row 1 occurs.
Case~\ref{cols33rows12} is where Row 1 is paired with any of Rows 2, 3 or 4. 
Case~\ref{cols33rows15} is where Row 1 and Row 5 are paired, 
Case~\ref{cols33rows22} pairs any of the Rows 2, 3 or 4,
Case~\ref{cols33rows55} assumes the square is contained in Row 5, 
while the last case,
Case~\ref{cols33rows25}, pairs one of Rows 2, 3 or 4 with Row 5. 

\esthm

\ssthm{cols33rows11} Suppose the square we want to show Mayer--Vietoris is 
embedded in the (3,1) box. Then it takes the form

\siv737000000
\sjv333000000

\BP
\postr13{B_{tj}\oplus X_{i'0}}c \Sposar23ar \postr33{B_{tj'}\oplus X_{i'0}}c
\Sposar12au                    \Sposar32au
\postr11{B_{tj}\oplus X_{i0}}c \Sposar21ar \postr31{B_{tj'}\oplus X_{i0}}c
\EP

\nin
and is therefore the sum of

\siv333000000
\sjv333000000

\SB\boxsix{
\postr13{B_{tj}}c \Sposar23ar \postr33{B_{tj'}}c
\Sposar12au                    \Sposar32au
\postr11{B_{tj}}c \Sposar21ar \postr31{B_{tj'}}c
}

\SB\boxsev{
\postr13{X_{i'0}}c \Sposar23ar \postr33{X_{i'0}}c
\Sposar12au                    \Sposar32au
\postr11{X_{i0}}c \Sposar21ar \postr31{X_{i0}}c
}

\barr
\PT\boxsix &\qquad{\rm and}\qquad&\PT\boxsev
\earr

\nin
which are both contractible triangles.

\essthm

\ssthm{cols33rows12} Suppose that we are dealing with a square, 
which is in the 
union of the (3,1) box and the $(3,x)$ box, 
where $x=2$, 3 or 4. Then the square 
is of the form

\siv737000000
\sjv333000000

\BP
\postr13{C_{i'j}\oplus X_{NW}}c \Sposar23ar \postr33{C_{i'j'}\oplus X_{NW}}c
\Sposar12au                    \Sposar32au
\postr11{B_{tj}\oplus X_{i0}}c \Sposar21ar \postr31{B_{tj'}\oplus X_{i0}}c
\EP

\nin
and is therefore the sum of

\siv333000000
\sjv333000000

\SB\boxsix{
\postr13{C_{i'j}}c \Sposar23ar \postr33{C_{i'j'}}c
\Sposar12au                    \Sposar32au
\postr11{B_{tj}}c \Sposar21ar \postr31{B_{tj'}}c
}

\SB\boxsev{
\postr13{X_{NW}}c \Sposar23ar \postr33{X_{NW}}c
\Sposar12au                    \Sposar32au
\postr11{X_{i0}}c \Sposar21ar \postr31{X_{i0}}c
}

\barr
\PT\boxsix &\qquad{\rm and}\qquad&\PT\boxsev
\earr

\nin
The first of these is part of the simplex $s_n$, while the second is
contractible.

\essthm

\ssthm{cols33rows15} If our square is contained in the union of the (3,1) and
(3,5) boxes, it is of the form

\siv737000000
\sjv333000000

\BP
\postr13{C_{i'j}}c \Sposar23ar \postr33{C_{i'j'}}c
\Sposar12au                    \Sposar32au
\postr11{B_{tj}\oplus X_{i0}}c \Sposar21ar \postr31{B_{tj'}\oplus X_{i0}}c
\EP

\nin
and is therefore the sum of

\siv333000000
\sjv333000000

\SB\boxsix{
\postr13{C_{i'j}}c \Sposar23ar \postr33{C_{i'j'}}c
\Sposar12au                    \Sposar32au
\postr11{B_{tj}}c \Sposar21ar \postr31{B_{tj'}}c
}

\SB\boxsev{
\postr13{0}c \Sposar23ar \postr33{0}c
\Sposar12au                    \Sposar32au
\postr11{X_{i0}}c \Sposar21ar \postr31{X_{i0}}c
}

\barr
\PT\boxsix &\qquad{\rm and}\qquad&\PT\boxsev
\earr

\nin
Once again, the first is part of the simplex $s_n$, while the second is
contractible.

\essthm

\ssthm{cols33rows22} Now suppose that we are dealing with a square inside the
third column, and somewhere in rows 2, 3 or 4. Then it has the form

\siv737000000
\sjv333000000

\BP
\postr13{C_{i'j}\oplus X_{NW}}c \Sposar23ar \postr33{C_{i'j'}\oplus X_{NW}}c
\Sposar12au                    \Sposar32au
\postr11{B_{ij}\oplus X_{NW}}c \Sposar21ar \postr31{B_{ij'}\oplus X_{NW}}c
\EP

\nin
and is therefore the sum of

\siv333000000
\sjv333000000

\SB\boxsix{
\postr13{C_{i'j}}c \Sposar23ar \postr33{C_{i'j'}}c
\Sposar12au                    \Sposar32au
\postr11{B_{ij}}c \Sposar21ar \postr31{B_{ij'}}c
}

\SB\boxsev{
\postr13{X_{NW}}c \Sposar23ar \postr33{X_{NW}}c
\Sposar12au                    \Sposar32au
\postr11{X_{NW}}c \Sposar21ar \postr31{X_{NW}}c
}

\barr
\PT\boxsix &\qquad{\rm and}\qquad&\PT\boxsev
\earr

\nin
The first of these is part of the simplex $s_n$, while the second is
contractible.

\essthm

\ssthm{cols33rows55} This case assumes that the entire square is in the
(3,5) box. In that case, it is just a square in the simplex $s_n$, quite
pure and unadulterated.
\essthm 

\ssthm{cols33rows25} The last case to consider is where we are dealing 
with a square inside the
third column, containing one of Rows 2, 3 or 4, and Row 5. Then it has the form

\siv737000000
\sjv333000000

\BP
\postr13{C_{i'j}}c \Sposar23ar \postr33{C_{i'j'}}c
\Sposar12au                    \Sposar32au
\postr11{B_{tj}\oplus X_{NW}}c \Sposar21ar \postr31{B_{tj'}\oplus X_{NW}}c
\EP

\nin
and is therefore the sum of

\siv333000000
\sjv333000000

\SB\boxsix{
\postr13{C_{i'j}}c \Sposar23ar \postr33{C_{i'j'}}c
\Sposar12au                    \Sposar32au
\postr11{B_{tj}}c \Sposar21ar \postr31{B_{tj'}}c
}

\SB\boxsev{
\postr13{0}c \Sposar23ar \postr33{0}c
\Sposar12au                    \Sposar32au
\postr11{X_{NW}}c \Sposar21ar \postr31{X_{NW}}c
}

\barr
\PT\boxsix &\qquad{\rm and}\qquad&\PT\boxsev
\earr

\nin
and, yet again, the first of these is part of the simplex 
$s_n$, while the second is
contractible.

\essthm

\sthm{cols34} {\bf Case of Column 4.} The last series of cases to consider 
is where the fourth column also occurs. In other words, we are concerned 
with the possibility that our square is in the union of the third and fourth
columns. Again, depending on which rows occur, we deal with cases. 
Case~\ref{cols34rows11} is where the square is entirely in the first row. 
Case~\ref{cols34rows12} is where Row 1 and Row 2 occur.  The final case,
Case~\ref{cols34rows22},  is where the square is either entirely in the 
second row or entirely in the third row (these cases turn out to be identical).

\esthm

\ssthm{cols34rows11} Suppose the square we want to show Mayer--Vietoris is 
embedded in the union of the
(3,1) and (4,1) boxes. Then it takes the form

\siv733000000
\sjv333000000

\BP
\postr13{B_{tj}\oplus X_{i'0}}c \Sposar23ar \postr33{X_{i'j'}}c
\Sposar12au                    \Sposar32au
\postr11{B_{tj}\oplus X_{i0}}c \Sposar21ar \postr31{X_{ij'}}c
\EP

\nin
and is therefore the sum of

\siv333000000
\sjv333000000

\SB\boxsix{
\postr13{X_{i'0}}c \Sposar23ar \postr33{X_{i'j'}}c
\Sposar12au                    \Sposar32au
\postr11{X_{i0}}c \Sposar21ar \postr31{X_{ij'}}c
}

\SB\boxsev{
\postr13{B_{tj}}c \Sposar23ar \postr33{0}c
\Sposar12au                    \Sposar32au
\postr11{B_{tj}}c \Sposar21ar \postr31{0}c
}

\barr
\PT\boxsix &\qquad{\rm and}\qquad&\PT\boxsev
\earr

\nin
where the first square is part of the simplex $s_n$, while the second is 
contractible.

\essthm

\ssthm{cols34rows12} Suppose the square we want to show Mayer--Vietoris  
involves Row 1 and Row 2. Then it takes the form

\siv733000000
\sjv333000000

\BP
\postr13{C_{i'j}\oplus X_{t0}}c \Sposar23ar \postr33{D_{i'j'}}c
\Sposar12au                    \Sposar32au
\postr11{A_{tj}\oplus X_{i0}}c \Sposar21ar \postr31{X_{ij'}}c
\EP

\nin
and  this is the second delicate square. The point is that we have a diagram
of squares

\siv333330000
\sjv333330000

\BP
\postr15{C_{i'j}}c \Sposar25ar \postr35{D_{i'0}}c  \Sposar45ar \postr55{ D_{i'j'}}c
\Sposar14au                    \Sposar34au                      \Sposar54au
\postr13{A_{tj}}c \Sposar23ar \postr33{X_{t0}}c \Sposar43ar \postr53{ X_{tj'}}c
                               \Sposar32au                    \Sposar52au
                            \postr31{X_{i0}}c \Sposar41ar \postr51{X_{ij'}}c
\EP

\nin
and the dual of the argument in Case~\ref{cols13rows12}
applies. Precisely, we have three semi--triangles

\siv337333500
\sjv300000000

\SB\boxtwo{
\postr11{A_{tj}}c \Sposar21ar \postr31{X_{tj'}\oplus C_{i'j}}c \Sposar41ar 
\postr51{D_{i'j'}}c  \Sposar61ar 
\postr71{\Sigma A_{tj}}c
}

\SB\boxthr{
\postr11{X_{i0}}c \Sposar21ar \postr31{X_{ij'}\oplus D_{i'0}}c \Sposar41ar 
\postr51{D_{i'j'}}c  \Sposar61ar 
\postr71{\Sigma X_{i0}}c
}

\SB\boxfou{
\postr11{X_{t0}}c \Sposar21ar \postr31{X_{tj'}\oplus D_{i'0}}c \Sposar41ar 
\postr51{D_{i'j'}}c  \Sposar61ar 
\postr71{\Sigma X_{t0}}c
}

\begin{equation}
\label{tr1'}
\PT\boxtwo
\end{equation}

\begin{equation}
\label{tr2'}
\PT\boxthr
\end{equation}

\nin
and

\begin{equation}
\label{tr3'}
\PT\boxfou
\end{equation}

\nin
and, exactly dually to Case~\ref{cols13rows12}, to prove that 

\siv733000000
\sjv333000000

\BP
\postr13{C_{i'j}\oplus X_{t0}}c \Sposar23ar \postr33{D_{i'j'}}c
\Sposar12au                    \Sposar32au
\postr11{A_{tj}\oplus X_{i0}}c \Sposar21ar \postr31{X_{ij'}}c
\EP

\nin
is Mayer--Vietoris, it suffices to establish that there is a semi--triangle

\siv53{9}333{7}00
\sjv300000000

\SB\boxtwo{
\postr11{A_{tj}\oplus X_{i0}}c \Sposar21ar 
\postr31{X_{ij'}\oplus C_{i'j}\oplus X_{t0}}c 
\Sposar41ar 
\postr51{D_{i'j'}}c  \Sposar61ar 
\postr71{\Sigma \{A_{tj}\oplus X_{i0}\}}c
}

\begin{equation}
\label{tr4'}
\qquad\PT\boxtwo
\end{equation}

\nin
But (\ref{tr4'}) can be obtained as a direct summand of the mapping 
cone on the natural map of 
semi--triangles

$$
(\ref{tr1'})\oplus(\ref{tr2'})\longrightarrow (\ref{tr3'}).
$$

\nin
What is more, the other direct summand is a contractible triangle.
This is also left as an exercise to the reader, who can once again
find a discussion in
Remark~\ref{R5.3}.

As in the proof of Case \ref{cols13rows12}, the reader will see that if we
start out with a simplex, where all the objects are in the heart of some
\tst, and all the  squares fold to give genuine triangles, then the 
square produced by the homotopy also folds to give a genuine triangle.
In particular, the differential is unique.
\essthm

\ssthm{cols34rows22} If Row 1 does not occur, our square must take the form

\siv733000000
\sjv333000000

\BP
\postr13{C_{i'j}\oplus X_{NW}}c \Sposar23ar \postr33{D_{i'j'}}c
\Sposar12au                    \Sposar32au
\postr11{C_{ij}\oplus X_{NW}}c \Sposar21ar \postr31{D_{ij'}}c
\EP

\nin
and is therefore the sum of

\siv333000000
\sjv333000000

\SB\boxsix{
\postr13{C_{i'j}}c \Sposar23ar \postr33{D_{i'j'}}c
\Sposar12au                    \Sposar32au
\postr11{C_{ij}}c \Sposar21ar \postr31{D_{ij'}}c
}

\SB\boxsev{
\postr13{X_{NW}}c \Sposar23ar \postr33{0}c
\Sposar12au                    \Sposar32au
\postr11{X_{NW}}c \Sposar21ar \postr31{0}c
}

\barr
\PT\boxsix &\qquad{\rm and}\qquad&\PT\boxsev
\earr

\nin
For the last time, the first of these is part of the simplex 
$s_n$, while the second is
contractible.

\essthm

This completes the discussion of the cases. We now know that there
is a way to choose the maps and differentials, so that at least every
square folds to give a semitriangle, and the maps that arise are some 
universally given matrices in the maps defining the original simplex. 
They are quite explicitly computable, from the various mapping cones we used.

Two things are not immediately clear. The first is why the maps and differentials,
that we have just shown how to choose, must be coherent. The second problem is 
how to show, that if we start with a somewhat restricted diagram, for 
instance a diagram having a lifting to model 
categories, the result of the homotopy
is another diagram with a lift.

\nin
{\bf Point 1: The differentials are coherent}

It is, of course, possible to compute this by brute force. But here is a cleaner
argument. Suppose that we start with a cell $s_n$ in

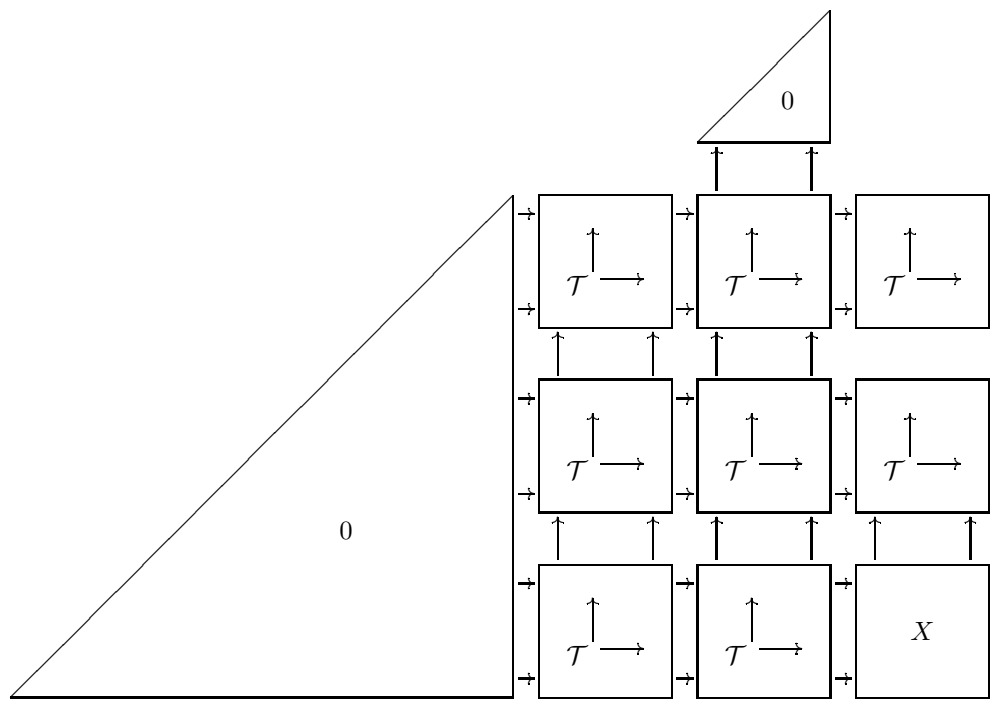
\BDI(0,260)
\UB(0,130)\boxeig
\EDI

\nin
where all the objects are in fact in some heart $\cal C$ of the triagulated 
category $\ct$, for some \tst. That is, our simplex 
happens to lie in the smaller simplicial set

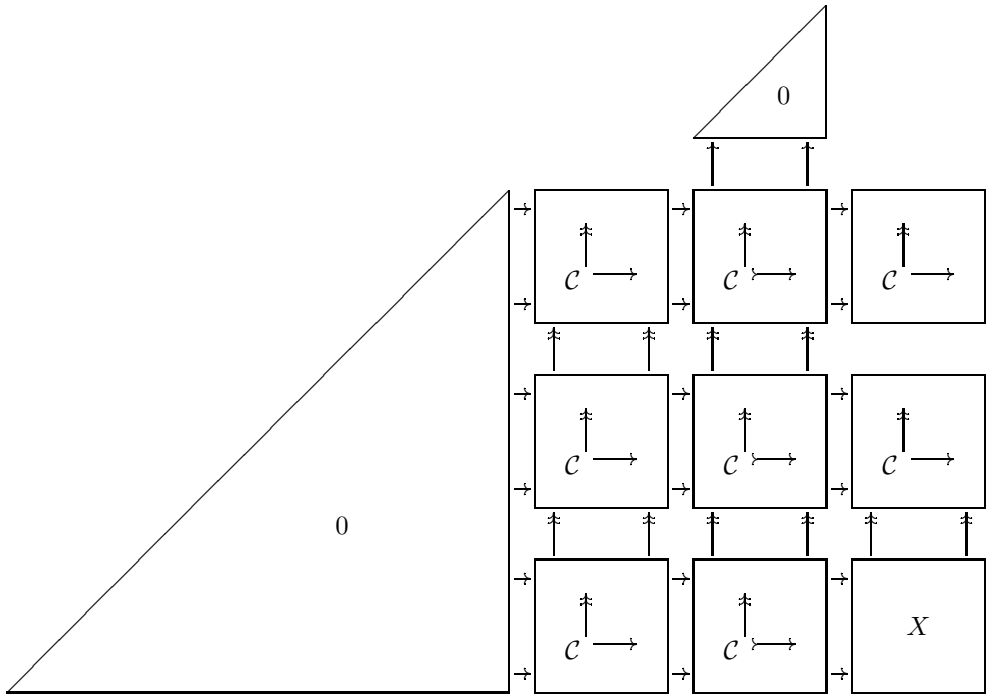
\siv{19}15151500
\sjv525252500

\SB\boxten{
                                     \trui570
                                      \bu56a
         \br25a  \fbx35{{\cal C}} ea  \br45a    \fbx55{{\cal C}} em    \br65a   \fbx75{{\cal C}} ea
                 \bu34e                \bu54e               
         \br23a   \fbx33{{\cal C}} ea  \br43a    \fbx53{{\cal C}} em    \br63a   \fbx73{{\cal C}} ea
                  \bu32e               \bu52e                   \bu72e
\trui110 \br21a   \fbx31{{\cal C}} ea  \br41a    \fbx51{{\cal C}} em    \br61a   \Fbx71X
}

\BDI(0,260)
\UB(0,130)\boxten
\EDI

\nin
Then we have just proved that the morphisms, in the diagram which should
give a cell of the homotopy

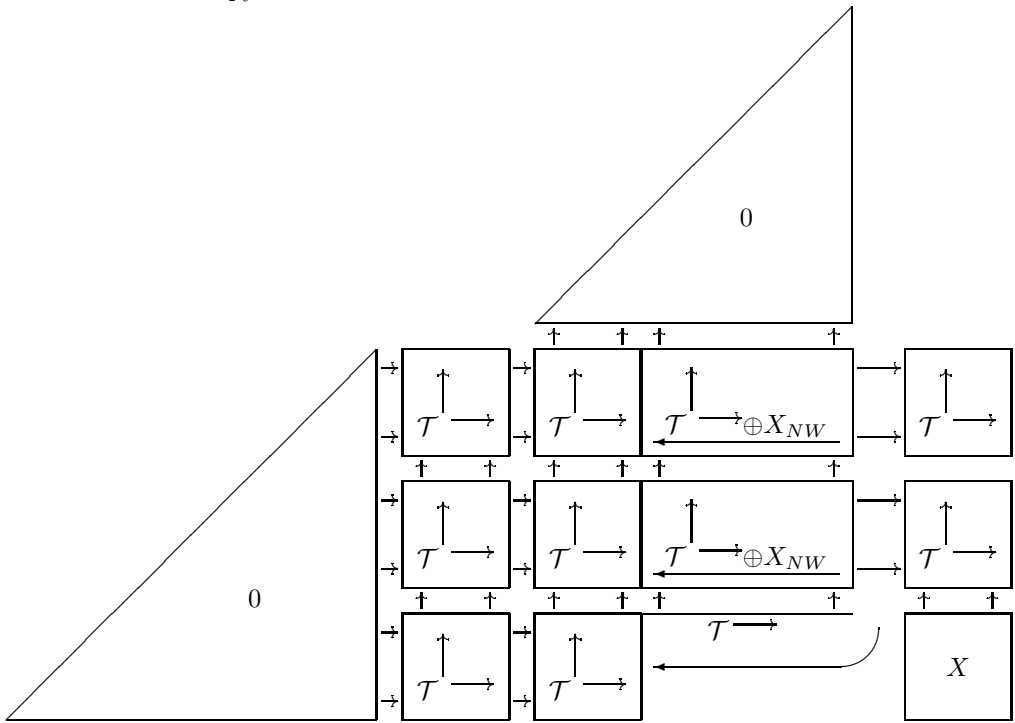
\siv{14}14148240
\sjv414141{12}00

\SB\boxten{
                                     \truj670
                                      \bu56a  \bu66a
         \br25a  \fbx35\ct aa  \br45a    \fbx55\ct aa \fhbx65\ct aa{ \oplus X_{NW}}l   \br75a   \fbx85\ct aa
                 \bu34a                \bu54a        \bu64a       
         \br23a   \fbx33\ct aa  \br43a    \fbx53\ct aa \fhbx63\ct aa{\oplus X_{NW}} l   \br73a   \fbx83\ct aa
                  \bu32a               \bu52a        \bu62a                \bu82a
\trui110 \br21a   \fbx31\ct aa  \br41a    \fbx51\ct aa  \fhtpl6611 \fli61\ct a{} t            \Fbx81X
} 
\BDI(0,270)
\UB(0,135)\boxten
\EDI
\nin
can be defined, as some universally given matrices in the structure maps of
the starting simplex, so that all squares fold to give {\em triangles}. Note
that in the two places where we took mapping cones on triangles in the proof,
we observed that in the special case above, they are mapping cones on good maps.

But of course the differentials are unique, and coherent, in this special case.
This is because a short exact sequence in a heart of a \tst\ 
corresponds to
a unique triangle; the differential is unique. To say that the differentials
are coherent, is to say that the restriction of the differential on a large
square to a smaller square is equal to the differential on the smaller square.
This gives some identity in the universally defined matrices of given maps
we have just constructed.

The identity must hold for all simplices all of whose objects are in a heart.
This means that the difference between the restriction of a large differential and
a small differential is a universally defined matrix of given maps, which 
vanishes whenever the starting simplex involves only objects of a heart. But
in the starting simplex, there was only one map connecting any two objects.
The matrix of the difference has components which can only be integer
multiples of the given map. If the integer is non--zero, it is easy to
construct a simplex of objects in a heart, on which it will not vanish. Hence
the difference between a small differential and the restriction of a large
differential must be the zero map; the differentials are coherent.

\nin
{\bf Point 2: There is a Waldhausen lift}

If we presume that we started with a simplex that had a 
lifting to a model category,
we wish to show that the universal diagram we have just constructed, in
order to define our homotopy, also has
such a lifting, in fact to the same model category.

The point is simple enough. A lifting is nothing other than some simplex
in

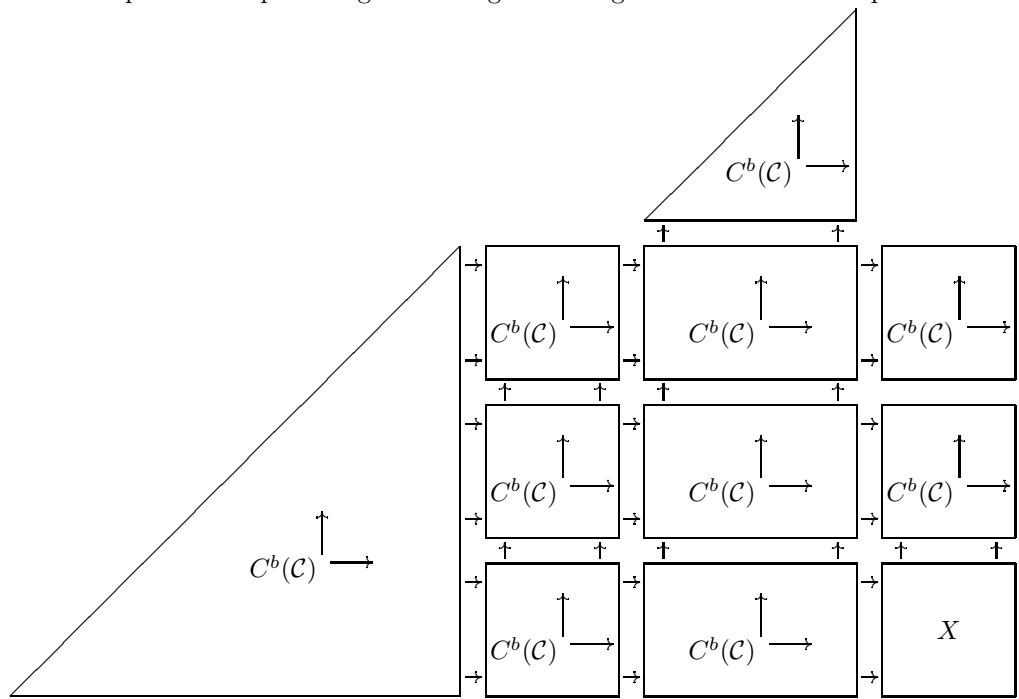
\siv{17}15181500
\sjv515151800

\SB\boxten{
                                     \trui57{C^b({\cal C})^{\barb aa}}
                                      \bu56a
         \br25a  \fbx35{C^b({\cal C})} aa  \br45a    \fbx55{C^b({\cal C})} aa 
   \br65a   \fbx75{C^b({\cal C})} aa
                 \bu34a                \bu54a               
         \br23a   \fbx33{C^b({\cal C})} aa  \br43a    \fbx53{C^b({\cal C})} aa
    \br63a   \fbx73{C^b({\cal C})} aa
                  \bu32a               \bu52a                   \bu72a
\trui11{C^b({\cal C})^{\barb aa}} \br21a   \fbx31{C^b({\cal C})} aa  \br41a  
  \fbx51{C^b({\cal C})} aa    \br61a   \Fbx71X
}
\BDI(0,260)
\UB(0,130)\boxten
\EDI

\nin
for the abelian category $C^b({\cal C}), $ which can always be 
viewed as the heart of its own derived category. 
We have shown that, with our universal
matrices between direct sums of objects, the homotopy 
will take a simplex in the heart to
a simplex in the heart, with even the differentials being right. Hence 
the existence of a lifting to the same model category follows.

This completes the proof of 
Theorem~\ref{well defined homotopies}.\hfill{$\Box$}

It seems appropriate to give a very simple illustration of the power of 
Theorem~\ref{well defined homotopies}. Let us now work carefully through the 
proof of Theorem~\ref{T5.1}. We remind the reader that in our notation, 
Theorem~\ref{T5.1} stands for Theorem 5.1 in \kthi.
For the reader's benefit, we recall the statement of the theorem.

\siv525000000
\sjv500000000

\SB\boxone{
\trui110  \br21a    \fbx31{\cs}ea
}

\siv500000000

\SB\boxfou{
\fbx11{\cs}aa
}

\nin {\bf Theoem~\ref{T5.1}.}  {\em The natural map

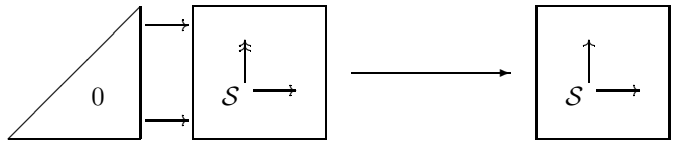
\BC
\begin{picture}(250,50)
\UBP{60}{25}\boxone    \Rar{160}{25}{60}{}       \UBP{225}{25}\boxfou
\end{picture}
\EC

\nin 
is a homotopy equivalence.}

{\it Proof of Theorem \ref{T5.1}}.  We considered the trisimplicial set and 
two projections

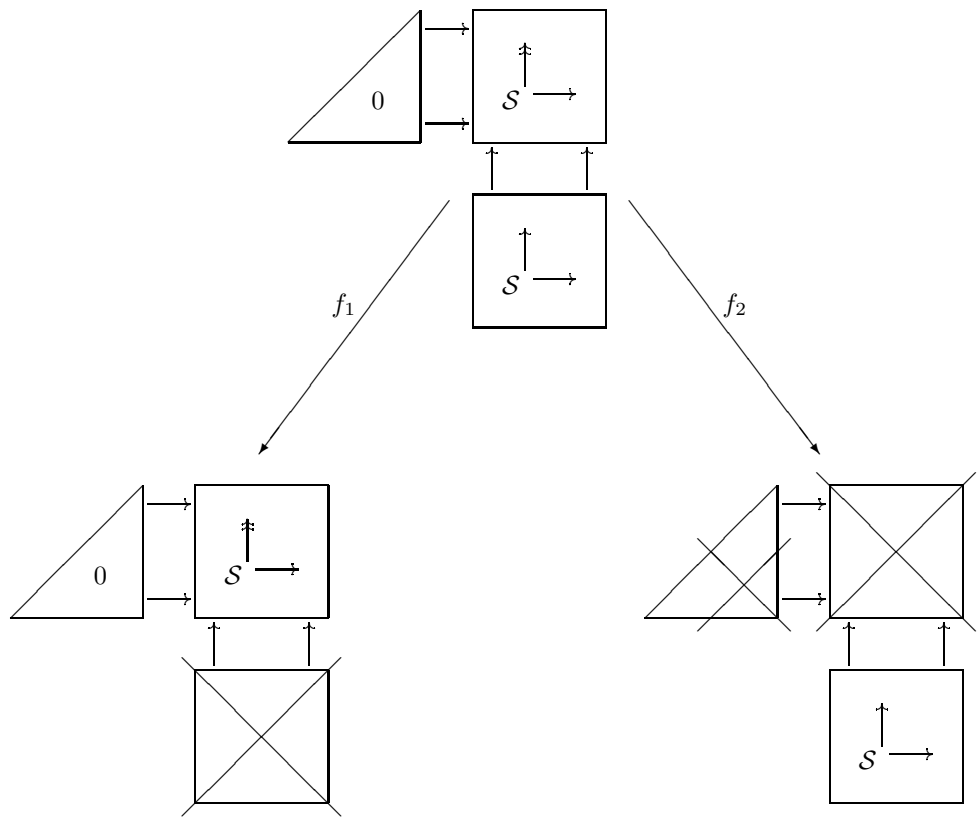
\siv525000000
\sjv525000000

\SB\boxone{
\trui130    \br23a       \fbx33{\cs}ea
                         \bu32a
                         \fbx31{\cs}aa
}

\SB\boxthr{
\xtrui13    \br23a       \xbx33
                         \bu32a
                         \fbx31{\cs}aa
}

\SB\boxtwo{
\trui130    \br23a       \fbx33{\cs}ea
                         \bu32a
                         \xbx31
}

\siv500000000

\SB\boxfou{
 \fbx13{\cs}aa
                         \bu12a
                         \fbx11{\cs}aa
}

\BC
\begin{picture}(300,300)
                             \UBP{130}{240}\boxone
       \slar{95}{180}{72}{-3}{-4}{f_1}           \slar{235}{180}{72}{3}{-4}{f_2}
\UBP{25}{60}\boxtwo                                         \UBP{265}{60}\boxthr 
\end{picture}
\EC

\siv525000000
\sjv525000000

\SB\boxfou{
\trui13Z     \br23a      \Fbx33Y
                         \bu32a
                         \fbx31{\cs}aa
}

\SB\boxfiv{
\trui130    \br23a       \fbx33{\cs}ea
                         \bu32a
                         \Fbx31X
}

\sjv552500000

\SB\boxsix{
\trui14Z     \br24a      \Fbx34Y
                         \bu33a
                         \Fhbx32{Y_S}d
                         \fbx31{\cs}aa
}

\nin
The Segal fiber of the map $f_1$ is the simplicial set

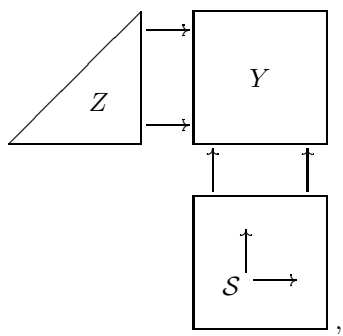
\BC
\usebox\boxfou\ ,
\EC

\nin
which is contracted by the homotopy

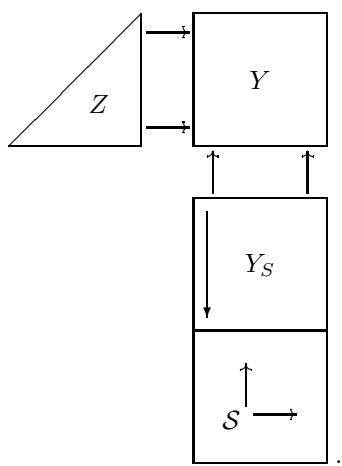
\BC
\usebox\boxsix\ .
\EC

\nin
The subtlety of the proof comes from trying to contract the Segal fiber of 
$f_2$, that is the simplicial set

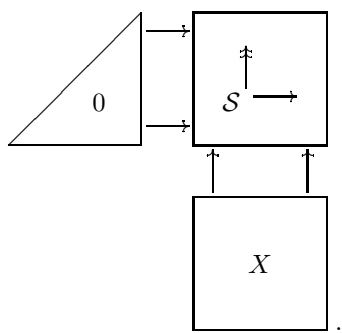
\BC
\usebox\boxfiv\ .
\EC

\nin
Recall the simplicial set

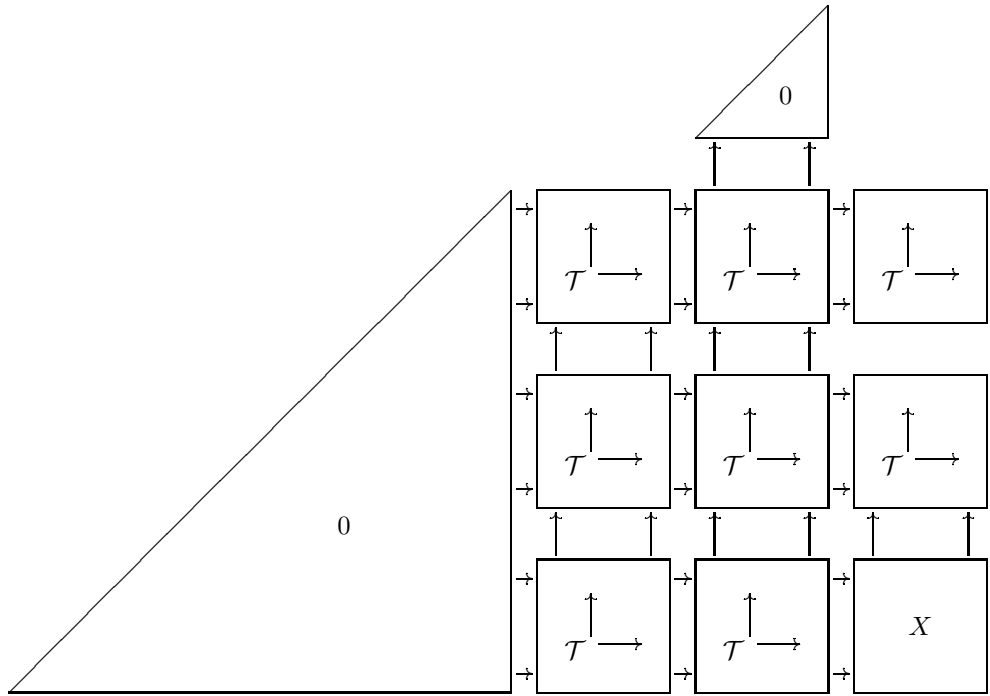
\siv{19}15151500
\sjv525252500

\SB\boxeig{
                                     \trui570
                                      \bu56a
         \br25a  \fbx35\ct aa  \br45a    \fbx55\ct aa    \br65a   \fbx75\ct aa
                 \bu34a                \bu54a               
         \br23a   \fbx33\ct aa  \br43a    \fbx53\ct aa    \br63a   \fbx73\ct aa
                  \bu32a               \bu52a                   \bu72a
\trui110 \br21a   \fbx31\ct aa  \br41a    \fbx51\ct aa    \br61a   \Fbx71X
}

\BC
\usebox\boxeig
\EC

\nin
 From now on we will refer to the above as the ``blueprint set".
By Theorem~\ref{well defined homotopies}, this blueprint set admits a homotopy

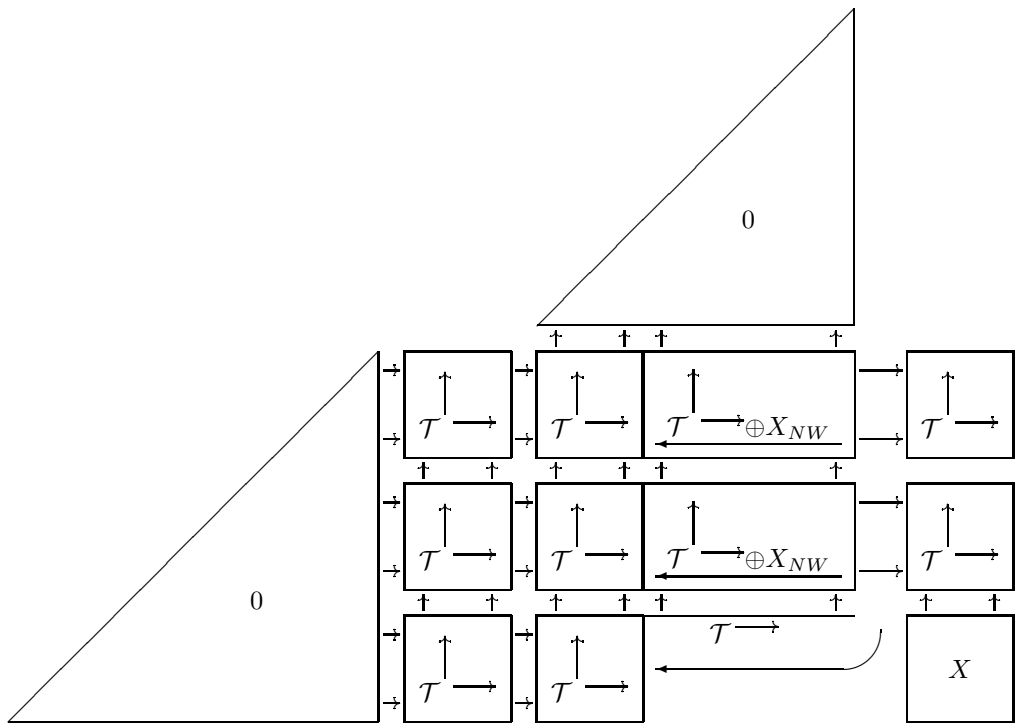
\siv{14}14148240
\sjv414141{12}00

\SB\boxten{
                                     \truj670
                                      \bu56a  \bu66a
         \br25a  \fbx35\ct aa  \br45a    \fbx55\ct aa \fhbx65\ct aa{ \oplus X_{NW}}l   \br75a   \fbx85\ct aa
                 \bu34a                \bu54a        \bu64a       
         \br23a   \fbx33\ct aa  \br43a    \fbx53\ct aa \fhbx63\ct aa{\oplus X_{NW}} l   \br73a   \fbx83\ct aa
                  \bu32a               \bu52a        \bu62a                \bu82a
\trui110 \br21a   \fbx31\ct aa  \br41a    \fbx51\ct aa  \fhtpl6611 \fli61\ct a{} t            \Fbx81X
}
\BDI(0,270)
\UB(0,135)\boxten
\EDI

\nin
which we will henceforth call the ``blueprint homotopy".
 
Consider now the smaller simplicial set

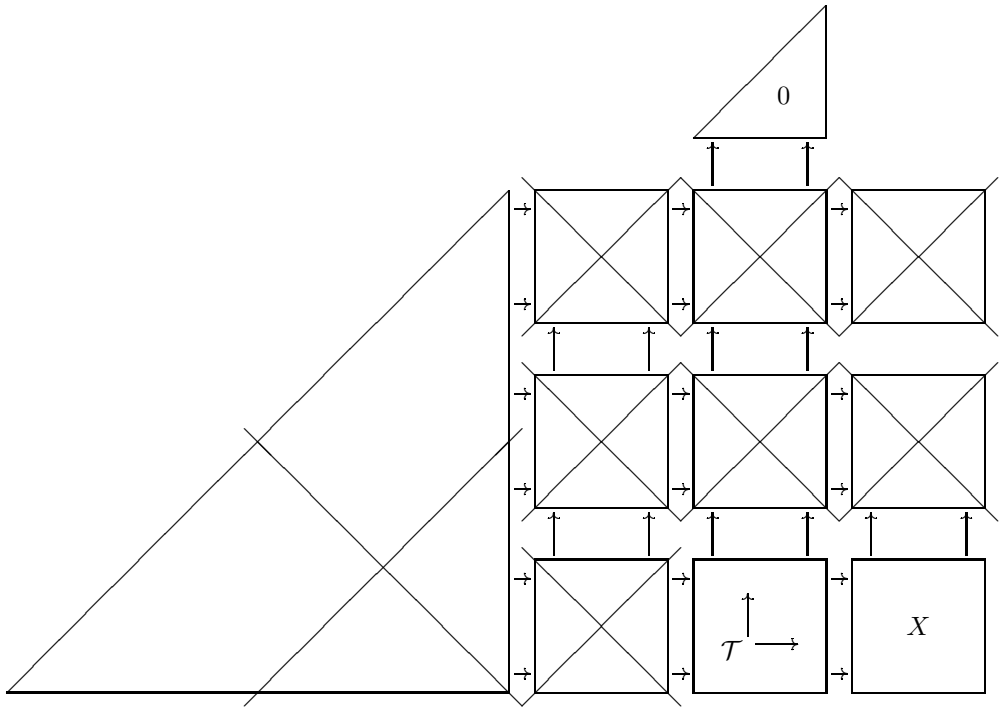
\siv{19}15151500
\sjv525252500

\BP
                                     \trui570
                                      \bu56a
         \br25a  \xbx35      \br45a    \xbx55          \br65a   \xbx75
                 \bu34a                \bu54a               
         \br23a   \xbx33     \br43a    \xbx53          \br63a   \xbx73
                  \bu32a               \bu52a                   \bu72a
\xtrui11 \br21a   \xbx31     \br41a    \fbx51\ct aa    \br61a   \Fbx71X
\EP

\nin
Because this simplicial set is obtained from the blueprint set by deletion, the
blueprint homotopy must be defined on it. This follows from the proof of 
Theorem~\ref{well defined homotopies}. The point is that the proof consisted of 
checking that certain squares are Mayer--Vietoris. In a diagram where part of 
the data is deleted, there are fewer squares to check Mayer--Vietoris, and it
follows we must have checked them all, back when we proved 
Theorem~\ref{well defined homotopies}.

Thus the homotopy, which we write more succinctly as

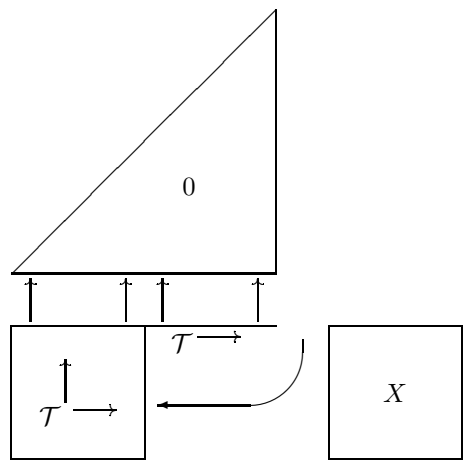
\siv552500000
\sjv52{10}000000

\BP
       \truj230
\bu12a  \bu22a
\fbx11\ct aa  \fhtpl2211  \fli21\ct a{}t  \Fbx41X
\EP

\nin
must be well--defined. Dually, we obtain a homotopy

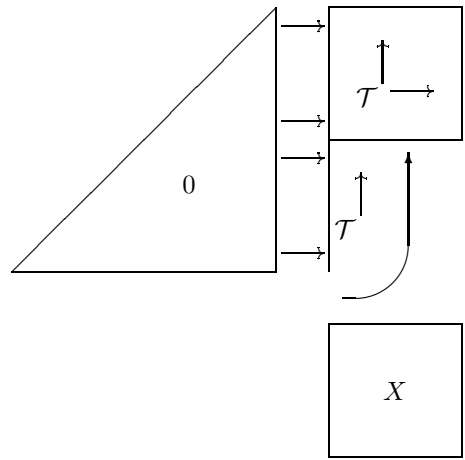
\siv{10}25000000
\sjv525500000

\BP
            \br24a  \fbx34\ct aa
\trui130    \br23a  \fhtpu3333 \fli33\ct a{}l
                    \Fbx31X
\EP

\nin
Now, if $\cs$ is an exact subcategory of the triangulated category $\ct$, we
have just shown that the homotopy

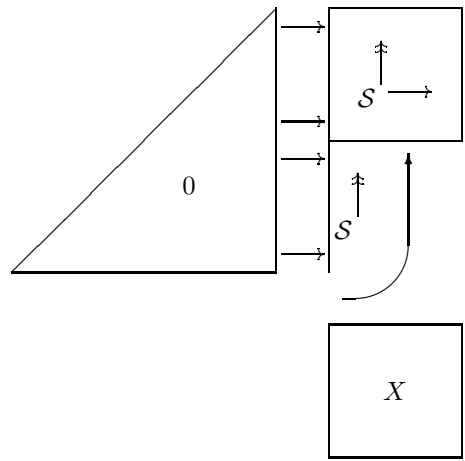
\BP
            \br24a  \fbx34\cs ea
\trui130    \br23a  \fhtpu3333 \fli33\cs e{}l
                    \Fbx31X
\EP

\nin
comes very close to being well--defined. The image of any homotopy cell in

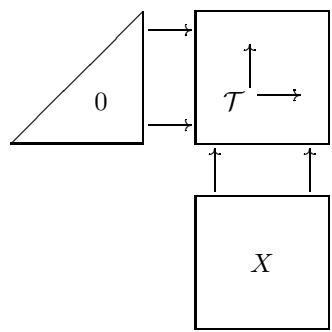
\siv525000000
\sjv525000000

\SB\boxone{
\trui130  \br23a  \fbx33\ct aa
                   \bu32a
                  \Fbx31X
}

\SB\boxtwo{
\trui130  \br23a  \fbx33\cs ea
                   \bu32a
                  \Fbx31X
}

\BC
\usebox\boxone
\EC

\nin
is unmistakeably a simplex. Checking that the homtopy is well defined amounts
to verifying only that the homotopy cells never leave

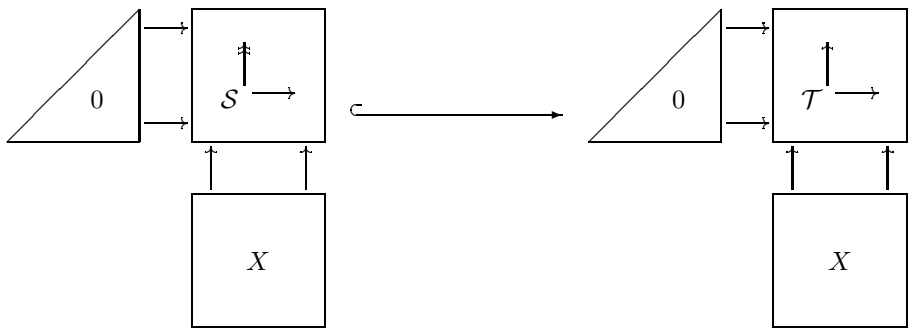
\BDI(220,120)
\UB(0,60)\boxtwo  \Ri(110,80,80){}   \UB(220,60)\boxone
\EDI

\nin
This reduces to showing that the objects stay in $\cs$, and the vertical 
morphisms are epi when they should be. 
The point of the exercise is that neither
reader nor writer need worry, whether the squares are Mayer--Vietoris, or the 
differentials coherent.  

This proves that the identity is homotopic to a map denoted
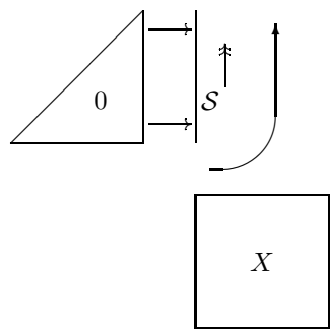
\siv525000000
\sjv525000000
\SB\boxone{
\trui130  \br23a  \fhtpu3333 \fli33\cs e{}l
                  \Fbx31X
}
\BC
\usebox\boxone
\EC
We refer the reader to the proof of Theorem~\ref{T5.1} in
\kthi, for
the argument showing that the above is homotopic to the null 
map. See also the argument starting on page~\pageref{page where}
of the present article,
under the title
``Continuation of the Proof of Theorem~\ref{safe kernels}''.
\hfill{$\Box$}

\siv{19}15151500
\sjv525252500

\SB\boxeig{
                                     \trui570
                                      \bu56a
         \br25a  \fbx35\ct aa  \br45a    \fbx55\ct aa    \br65a   \fbx75\ct aa
                 \bu34a                \bu54a               
         \br23a   \fbx33\ct aa  \br43a    \fbx53\ct aa    \br63a   \fbx73\ct aa
                  \bu32a               \bu52a                   \bu72a
\trui110 \br21a   \fbx31\ct aa  \br41a    \fbx51\ct aa    \br61a   \Fbx71X
}

\siv400000000
\sjv100000000

\SB\boxone{
\Sposar11hr}

Now that we have given a very simple--minded application of 
Theorem~\ref{well defined homotopies}, it  seems to be time to give a 
sophisticated one, which we will actually need in this article. Recall that an
arrow of type $\PT\boxone$ stands for a morphism in $\ct_{[0,n]}$ inducing an
isomorphism on $H^0.$

\siv525000000
\sjv500000000

\SB\boxone{
\trui110  \br21a  \bx31ha0n
}

\SB\boxtwo{
\xtrui11  \br21a  \bx31ha0n
}

\thm{safe kernels} Let $\ct$ be a triangulated category with a 
{\it t--}structure. Then the natural projection

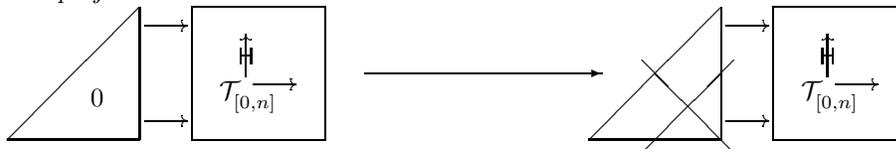
\BDI(220,50)
\UB(0,25)\boxone  \Ra(120,25,90){}   \UB(220,25)\boxtwo
\EDI

\nin
induces a homotopy equivalence.
\ethm

\pf We will first give the argument formally. This means, we will not
worry about
the subtleties, which are caused by having 
to consider the various modifications
introduced in Remark~\ref{R3.4}, to the basic simplicial set of 
Construction~\ref{C3.3}. A Group 3 reader should be able to fill in the
details himself. For Group 2 readers, we will discuss in 
Remark~\ref{group2 trunc} why the homotopies in the argument are 
all well defined, in the particular 
modification of Construction~\ref{C3.3} which
I advised the reader to adopt. In other words, Remark~\ref{group2 trunc} will
show that all the simplices that arise decompose as sums of simplices, 
each with a
lifting to a model category.

 We consider the more complicated diagram of a simplicial set and two
projections:

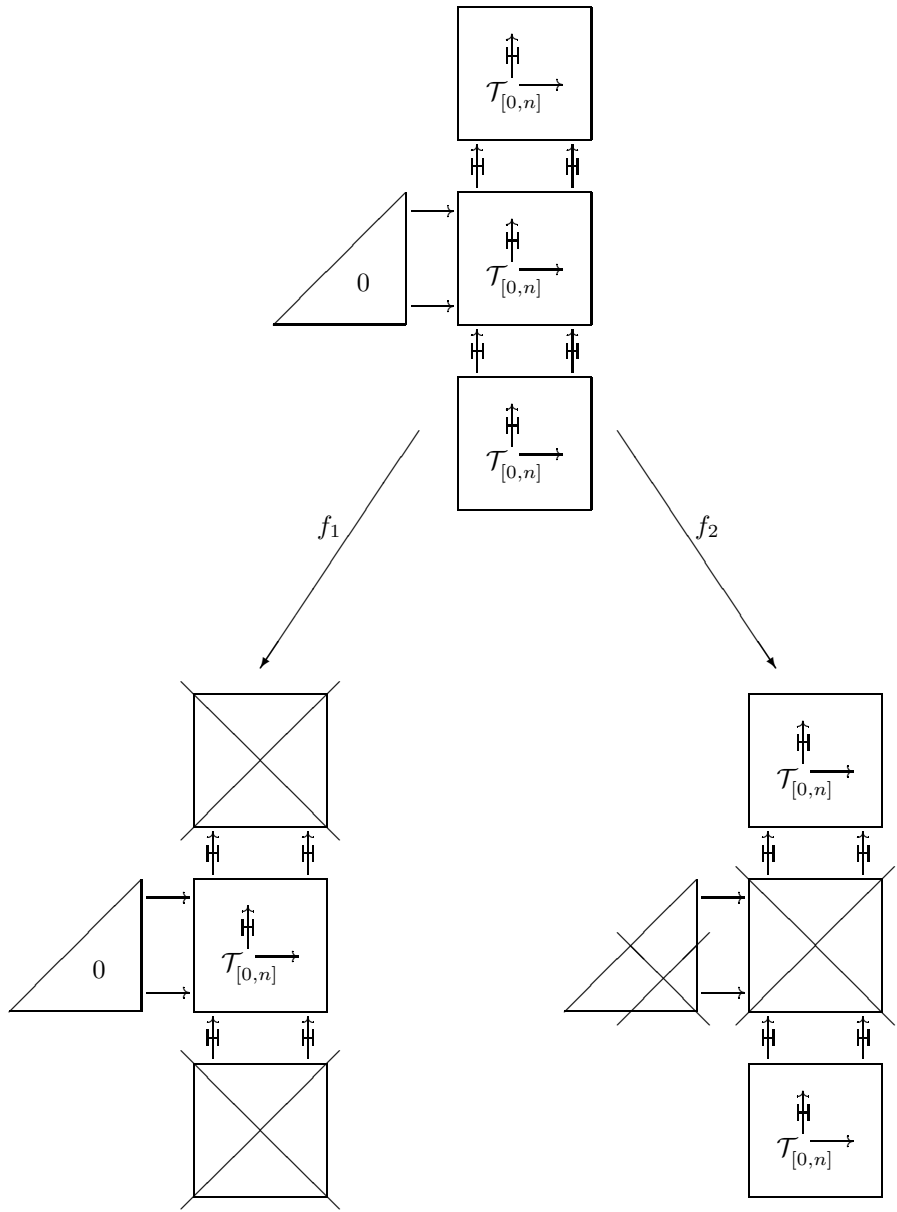
\siv525000000
\sjv525250000

\SB\boxtwo{
                     \bx35ha0n
                     \bu34h
\trui130  \br23a     \bx33ha0n
                     \bu32h
                     \bx31ha0n
}

\SB\boxthr{
                     \xbx35
                     \bu34h
\trui130  \br23a     \bx33ha0n
                     \bu32h
                     \xbx31
}

\SB\boxfou{
                     \bx35ha0n
                     \bu34h
\xtrui13  \br23a     \xbx33
                     \bu32h
                     \bx31ha0n
}

\SB\boxfiv{
                     \bx35ha0n
                     \bu34h
\trui13A  \br23a     \Fbx33B
                     \bu32h
                     \bx31ha0n
}

\SB\boxsix{
                     \Fbx35B
                     \bu34h
\trui130  \br23a     \bx33ha0n
                     \bu32h
                     \Fbx31A
}

\BDI(0,450)
                 \UB(-10,355)\boxtwo
       \sla(-45,245,60,-2,-3){f_1}  \sla(90,245,60,2,-3){f_2}
\UB(-110,95)\boxthr                         \UB(100,95)\boxfou
\EDI

\nin
and our theorem will follow, once we establish that $f_1$ and $f_2$ are 
homotopy equivalences. For $f_1$ this is very easy; the Segal fiber is the
simplicial set

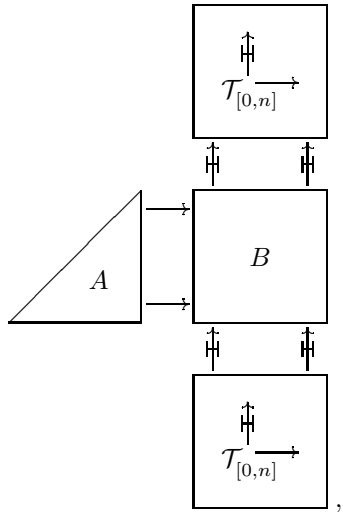
\BC
\usebox\boxfiv\ ,
\EC

\nin
which is contracted by the homotopy

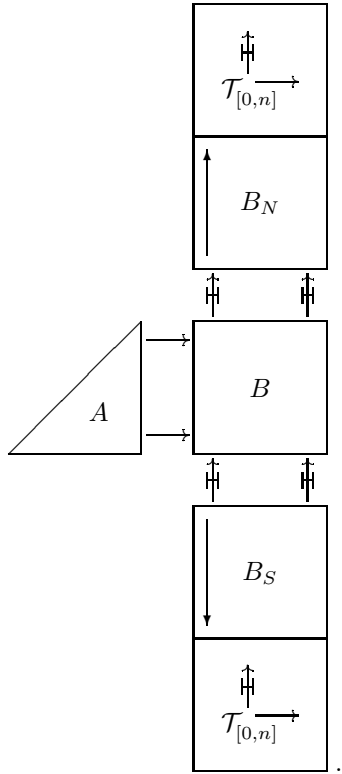
\siv525000000
\sjv552525500

\SB\boxfiv{
                     \bx37ha0n
                     \Fhbx36{B_N}u
                     \bu35h
\trui14A  \br24a     \Fbx34B
                     \bu33h
                     \Fhbx32{B_S}d
                     \bx31ha0n
}

\BC
\usebox\boxfiv\ .
\EC

\nin
The map $f_2$ is slightly subtler to prove a homotopy equivalence. The proof
amounts to using the homotopy

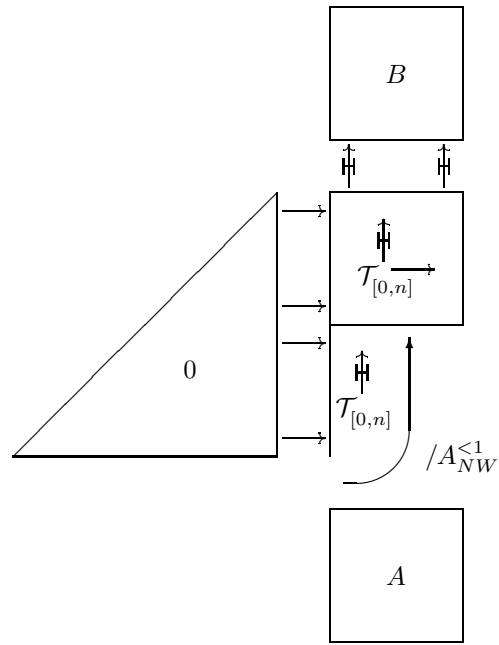
\sjv525525000
\siv{10}25000000

\SB\boxone{
                     \Fbx36B
                     \bu35h
\trui130  \br24a     \bx34ha0n
          \br23a     \htp33h0nu   \postr33{/A_{NW}^{<1}}e
                     \Fbx31A
}

\BC
\usebox\boxone
\EC

\nin
to factor the identity on the Segal fiber

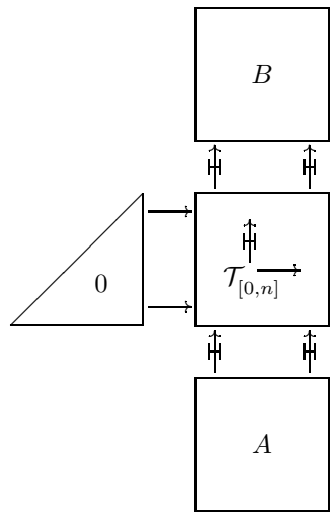
\BC
\usebox\boxsix
\EC

\nin
through the contractible simplicial set

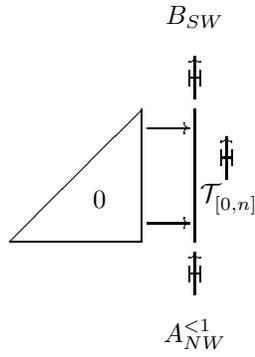
\siv522000000
\sjv325230000

\BP
                     \postr35{B_{SW}}l
                     \posar24he
\trui130   \br23a    \li33h0n{}l
                     \posar32hw
                     \postr31{A_{NW}^{<1}}l
\EP

\nin
The contractibility of this last simplicial set is by the 
contraction to the initial object.

But in this section we are in the business of worrying why the non--trivial 
homotopies are well defined. In particular, we should concern ourselves with
the homotopy

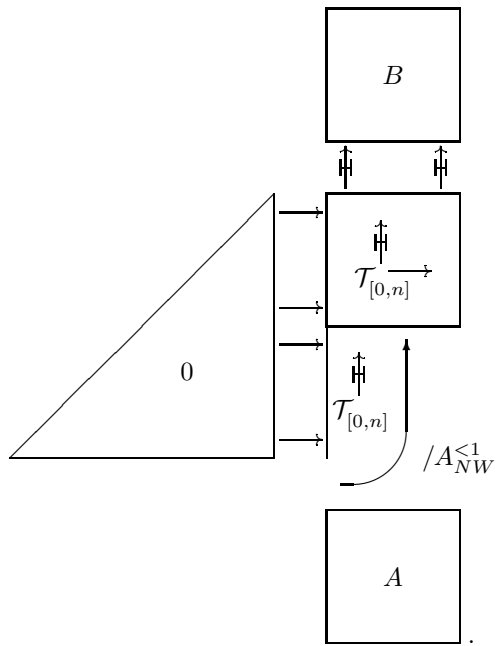
\BC
\usebox\boxone\ .
\EC

\nin
To this end, we recall our simplicial set

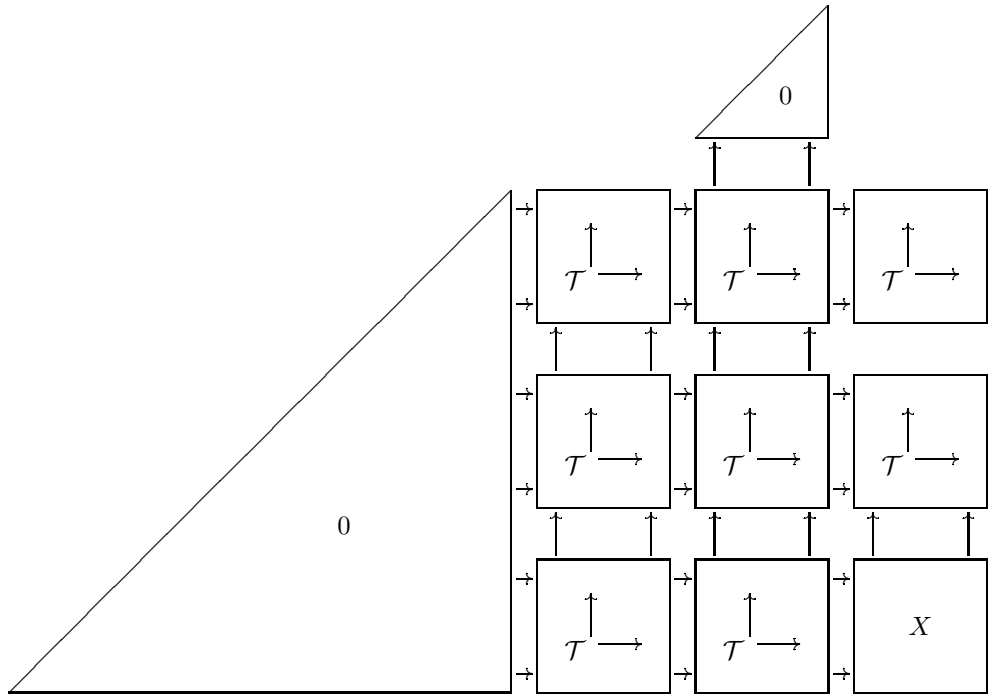
\BC
\usebox\boxeig
\EC

\nin
By deleting some of the data, we obtain the smaller simplicial set

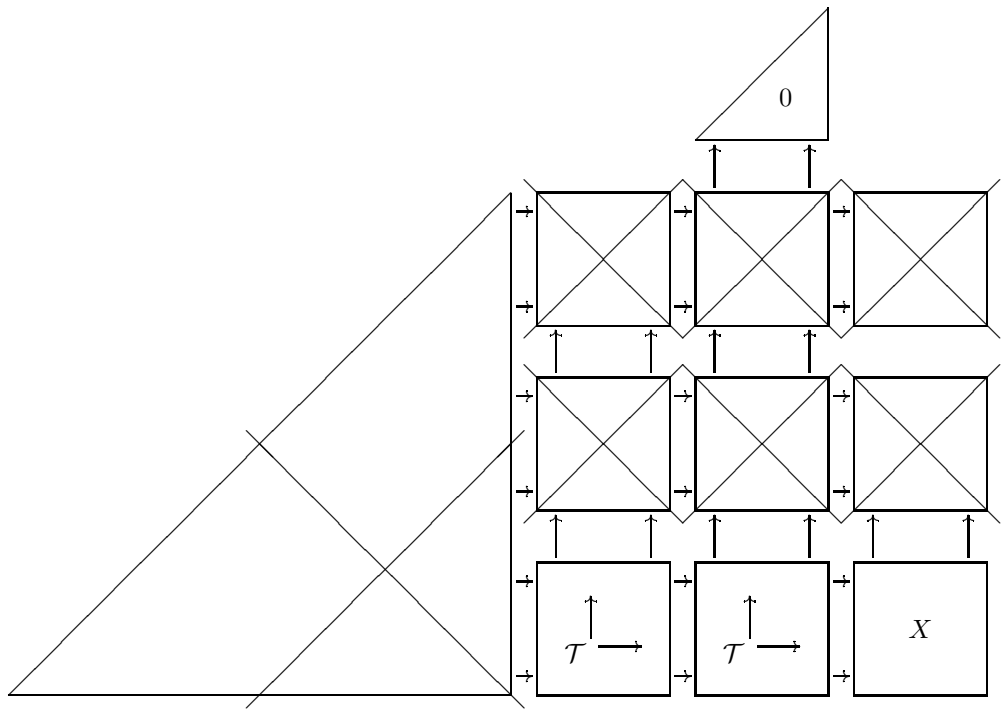
\siv{19}15151500
\sjv525252500

\BP                                     \trui570
                                      \bu56a
         \br25a  \xbx35  \br45a    \xbx55    \br65a   \xbx75
                 \bu34a                \bu54a               
         \br23a   \xbx33  \br43a    \xbx53    \br63a   \xbx73
                  \bu32a               \bu52a                   \bu72a
\xtrui11 \br21a   \fbx31\ct aa  \br41a    \fbx51\ct aa    \br61a   \Fbx71X
\EP

\nin
Theorem~\ref{well defined homotopies} guarantees that there is a well--defined
 homotopy

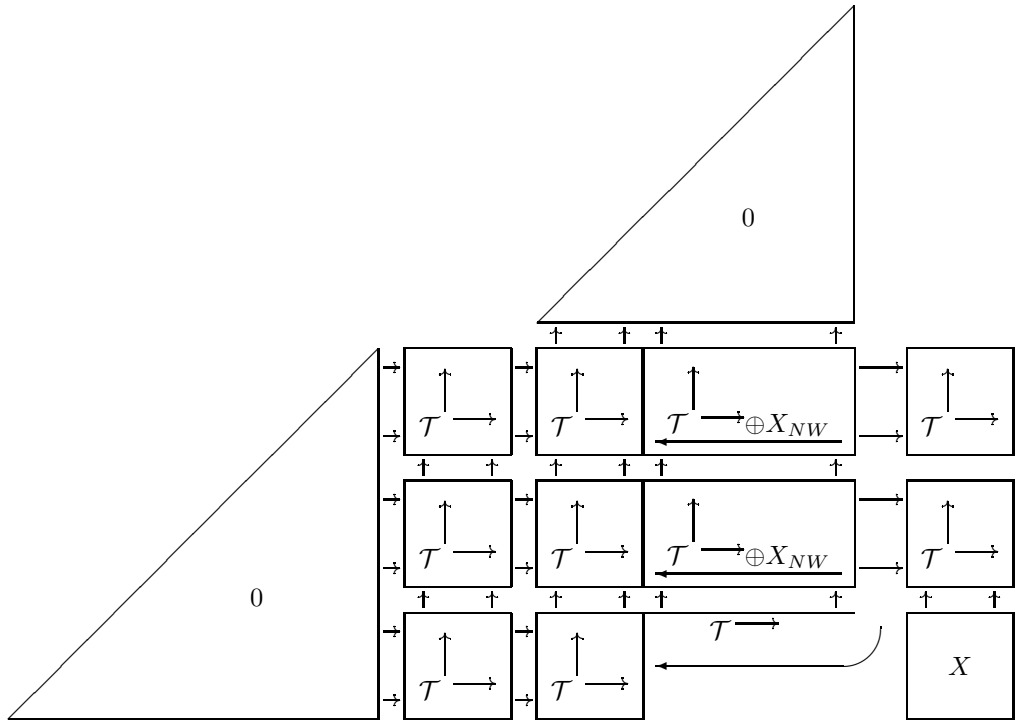
\siv{14}14148240
\sjv414141{12}00

\SB\boxten{
                                     \truj670
                                      \bu56a  \bu66a
         \br25a  \fbx35\ct aa  \br45a    \fbx55\ct aa \fhbx65\ct aa{ \oplus X_{NW}}l   \br75a   \fbx85\ct aa
                 \bu34a                \bu54a        \bu64a       
         \br23a   \fbx33\ct aa  \br43a    \fbx53\ct aa \fhbx63\ct aa{\oplus X_{NW}} l   \br73a   \fbx83\ct aa
                  \bu32a               \bu52a        \bu62a                \bu82a
\trui110 \br21a   \fbx31\ct aa  \br41a    \fbx51\ct aa  \fhtpl6611 \fli61\ct a{} t            \Fbx81X
}
\BDI(0,270)
\UB(0,135)\boxten
\EDI

\nin
and from the proof, it follows the homotopy remains well--defined when we delete
some of the structure. In particular, on our smaller simplicial set we
obtain a homotopy which we will denote

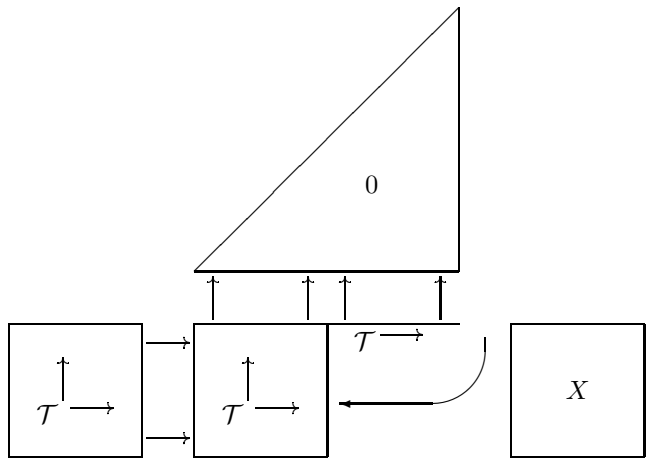
\siv525525000
\sjv52{10}000000

\BP
                   \truj430
                   \bu32a               \bu42a
\fbx11\ct aa \br21a\fbx31\ct aa \fhtpl4411 \fli41\ct a{} t            \Fbx61X
\EP

\nin
This homotopy is dual to

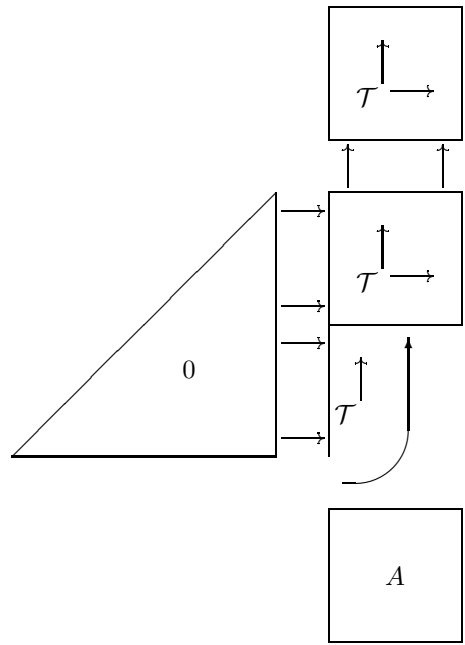
\siv{10}25000000
\sjv525525000

\BP
                     \fbx36\ct aa
                     \bu35a
\trui130  \br24a     \fbx34\ct aa
          \br23a     \fhtpu3333   \fli33\ct a{}l   
                     \Fbx31A
\EP

\nin
In the homotopy, the contents of the top box are fixed. So we may also view it
as a homotopy on the simplicial set, where the top is constrained to be fixed; 
that is, the homotopy whose symbol would be

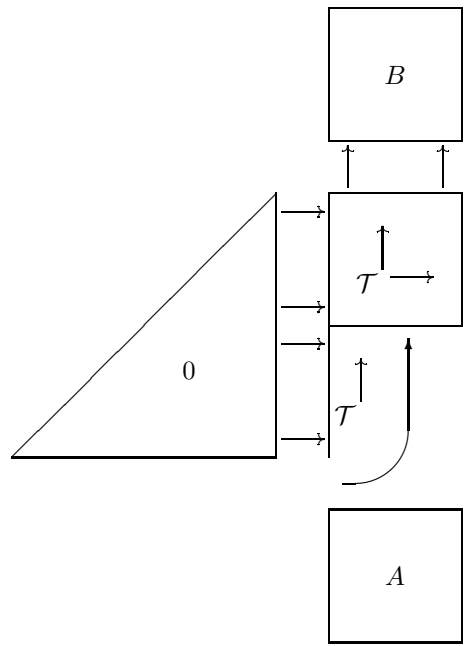
\BP
                     \Fbx36B
                     \bu35a
\trui130  \br24a     \fbx34\ct aa
          \br23a     \fhtpu3333   \fli33\ct a{}l   
                     \Fbx31A

\EP

So far this argument has been completely painless, appealing directly to the 
proof of Theorem~\ref{well defined homotopies}. To obtain the homotopy

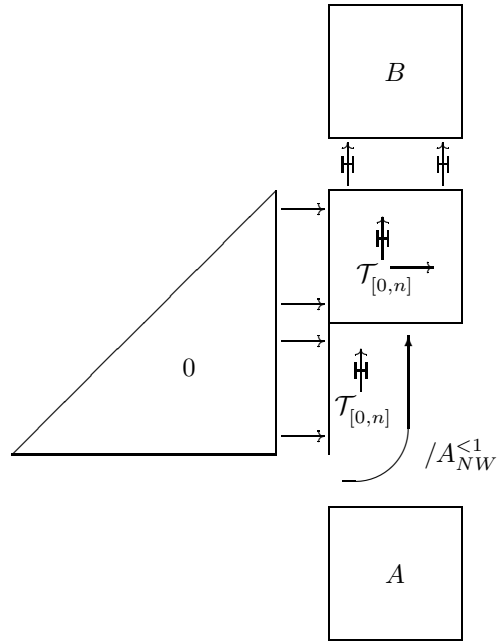
\BP
                     \Fbx36B
                     \bu35h
\trui130  \br24a     \bx34ha0n
          \br23a     \htp33h0nu   \postr33{/A_{NW}^{<1}}e
                     \Fbx31A
\EP

\nin
we have to worry a little about the \tst\ truncation.  Suppose we start with a 
simplex $s_n$ in the simplicial set

\BC
\usebox\boxsix\ .
\EC

It is given as a diagram

\siv{13}2{13}000000
\sjv{13}2{13}2{13}0000

\A10n

\A30q
\B10p
\B30n

\B50r

\SB\boxone{
                           \Sbxp35ha{B_\Q}
                           \Sbu34h
\Strup13aa{R_\Q} \Sbr23a   \Sbxp33ha{S_\Q}
                           \Sbu32h
                           \Sbxp31ha{A_\Q}
\fmbx3355 \fmbx3311
}

\barr
s_n=\mtrx\boxone
\earr

 \nin
Then the homotopy

\siv{10}25000000
\sjv525525000

\BP
                     \Fbx36B
                     \bu35a
\trui130  \br24a     \fbx34\ct aa
          \br23a     \fhtpu3333   \fli33\ct a{}l   
                     \Fbx31A

\EP

\nin
takes $s_n$ to $n+1$ distinct $(n+1)$--simplices, the $i^{th}$ of which is
 a diagram

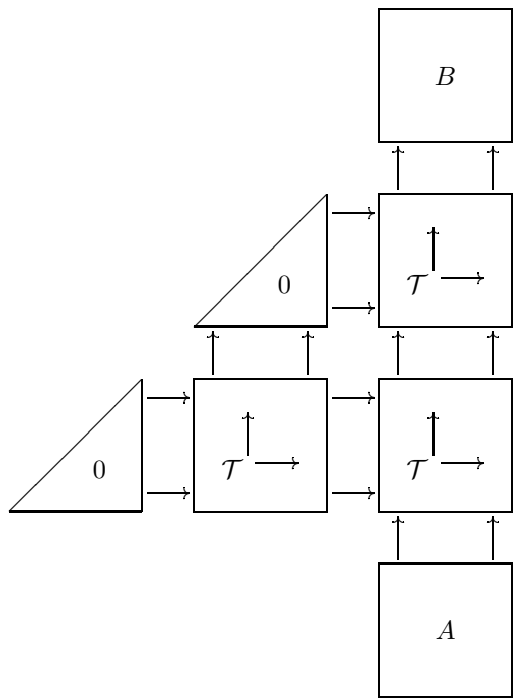
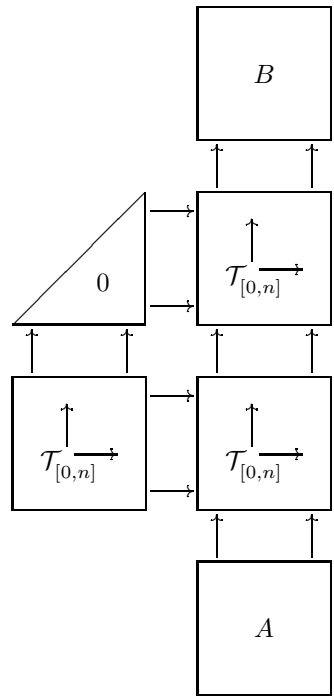
\small

\siv{10}1{15}1{14}0000
\sjv{13}2{13}2{13}2{13}00
\A10i
\A3in
\A50q
\B10p
\B30i
\B5in
\B70r

\SB\boxten{
                                      \Sbxp57ha{B_\Q}
                                      \Sbu56h
                    \Strup35aa{R_\Q}\Sbr45a \Sbxp55ha{S_\Q}
                    \Sbu34a                   \Sbu54a
\Strup13aa{R_\Q}\Sbrsmall23a
\Sbxpsmall33aa{\scriptstyle R_\Q\!\oplus\! A_{NW}}\Sbrsmall43a 
              \Sbxpsmall53ha{\scriptstyle S_{\I 0}\!\oplus\! A_{p\J}}
                                           \Sbu52a
                                         \Sbxp51ha{A_\Q}
\fmbx5511 \fmbx5577
}
\BDI(0,580)
\UB(,290)\boxten
\EDI

\normalsize

\nin
and Theorem~\ref{well defined homotopies} assures us that this is a simplex in
the simplicial set

\siv525250000
\sjv525252500

\BP
                     \Fbx57B
                     \bu56a
\trui350  \br45a     \fbx55\ct aa
\bu34a               \bu54a
\trui130\br23a\fbx33\ct aa \br43a  \fbx53\ct aa
                     \bu52a   
                     \Fbx51A

\EP

\nin
By deletion,

\siv00{19}2{17}0000
\sjv{13}2{13}2{13}2{13}00
\A10i
\A3in
\A50q
\B10p
\B30i
\B5in
\B70r

\SB\boxten{
                                      \Sbxp57ha{B_\Q}
                                      \Sbu56h
                    \Strup35aa{R_\Q}\Sbr45a \Sbxp55ha{S_\Q}
                    \Sbu34a                   \Sbu54a
                   \Sbxp33aa{R_\Q\oplus A_{NW}}\Sbr43a 
                                         \Sbxp53ha{S_{\I 0}\oplus A_{p\J}}
                                           \Sbu52a
                                         \Sbxp51ha{A_\Q}
\fmbx5511 \fmbx5577
}

\BDI(0,580)
\UB(,290)\boxten
\EDI

\nin
can certainly be viewed as a simplex in the simplicial set

\siv525000000
\sjv525252500

\BP
                     \Fbx37B
                     \bu36a
\trui150  \br25a     \bx35aa0n
\bu14a               \bu34a
\bx13aa0n \br23a  \bx33aa0n
                     \bu32a   
                     \Fbx31A

\EP

\nin
and the relevant point is that, for the special simplices above,
 there is a mono of $A_{NW}^{<1}$ into
the second row of this simplicial set, whose composite to the third row is 
zero. Precisely, the injection  $A_{NW}^{<1}\hookrightarrow A_{NW}$ gives
a compatible family of morphisms $A_{NW}^{<1}\hookrightarrow R_{ij}\oplus A_{NW}$,
whose composites with the projection $R_{ij}\oplus A_{NW}\longrightarrow R_{ij}$
is zero. This gives a map of  $A_{NW}^{<1}$ into the part of the
diagram below surrounded by a dashbox
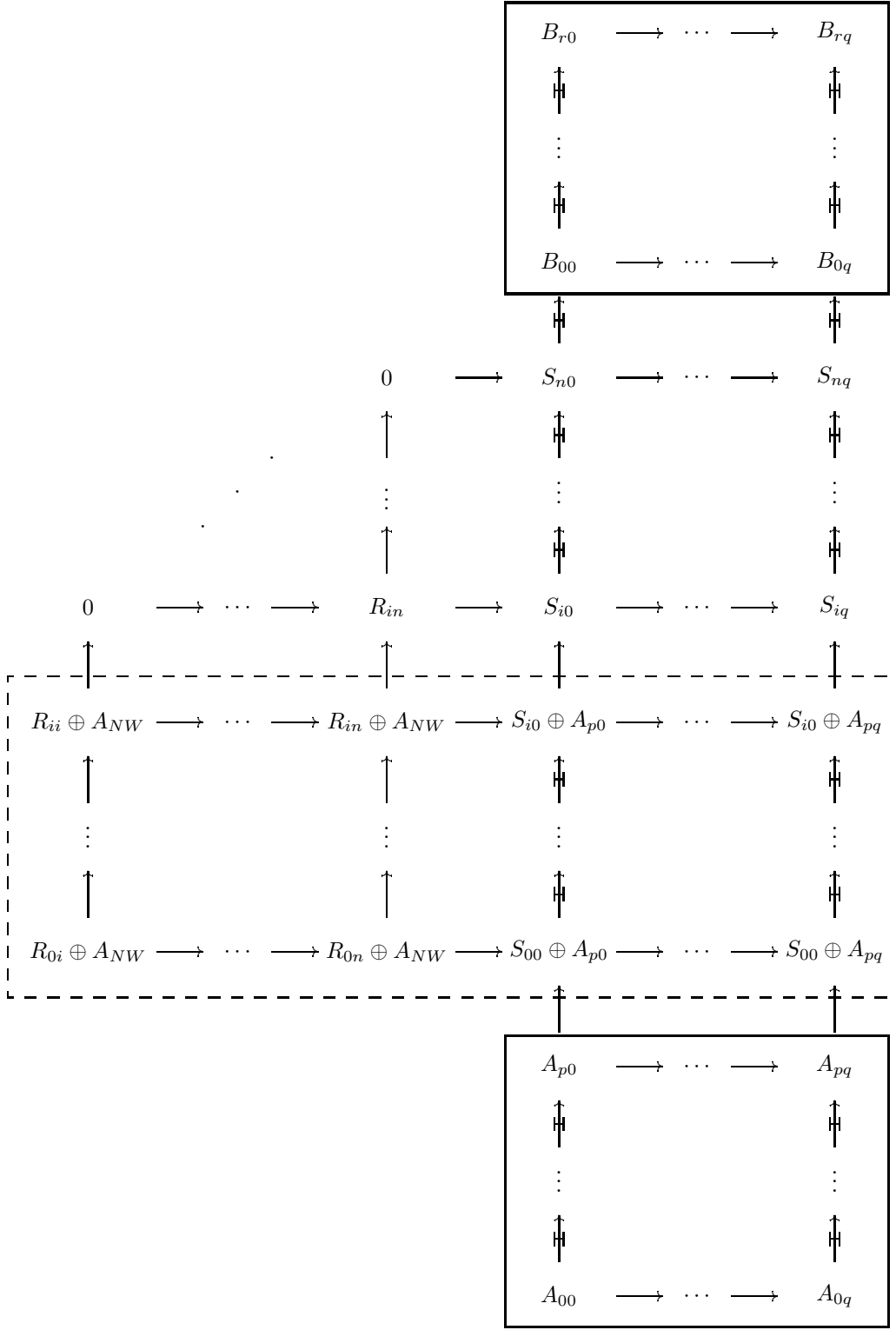
\eject
\vsize=10in
\siv00{19}2{17}0000
\sjv{13}2{13}2{13}2{13}00
\A10i
\A3in
\A50q
\B10p
\B30i
\B5in
\B70r

\SB\boxten{
                                      \Sbxp57ha{B_\Q}
                                      \Sbu56h
                    \Strup35aa{R_\Q}\Sbr45a \Sbxp55ha{S_\Q}
                    \Sbu34a                   \Sbu54a
                   \Sbxp33aa{R_\Q\oplus A_{NW}}\Sbr43a 
                                         \Sbxp53ha{S_{\I 0}\oplus A_{p\J}}
                                           \Sbu52a
                                         \Sbxp51ha{A_\Q}
\fmbx5511 \fmbx5577 \dshbx3533
}
\BDI(0,580)
\UB(,290)\boxten
\EDI

\nin
and what is especially good about this special case, is that the maps from
$A_{NW}^{<1}$ into the region within the dashbox are all monos; that is, all
the objects in the diagram are in $\ct^{\ge 0}$, and the map from the object
$A_{NW}^{<1}$ into any of the objects inside the dashbow is mono in $\ct^{\ge 0}$.
Furthermore, by the time we compose the map from $A_{NW}^{<1}$ with any 
map leaving the dashbox,
the composite vanishes.

The unsubtle truncation homotopy allows us then to define a homotopy whose
cells are
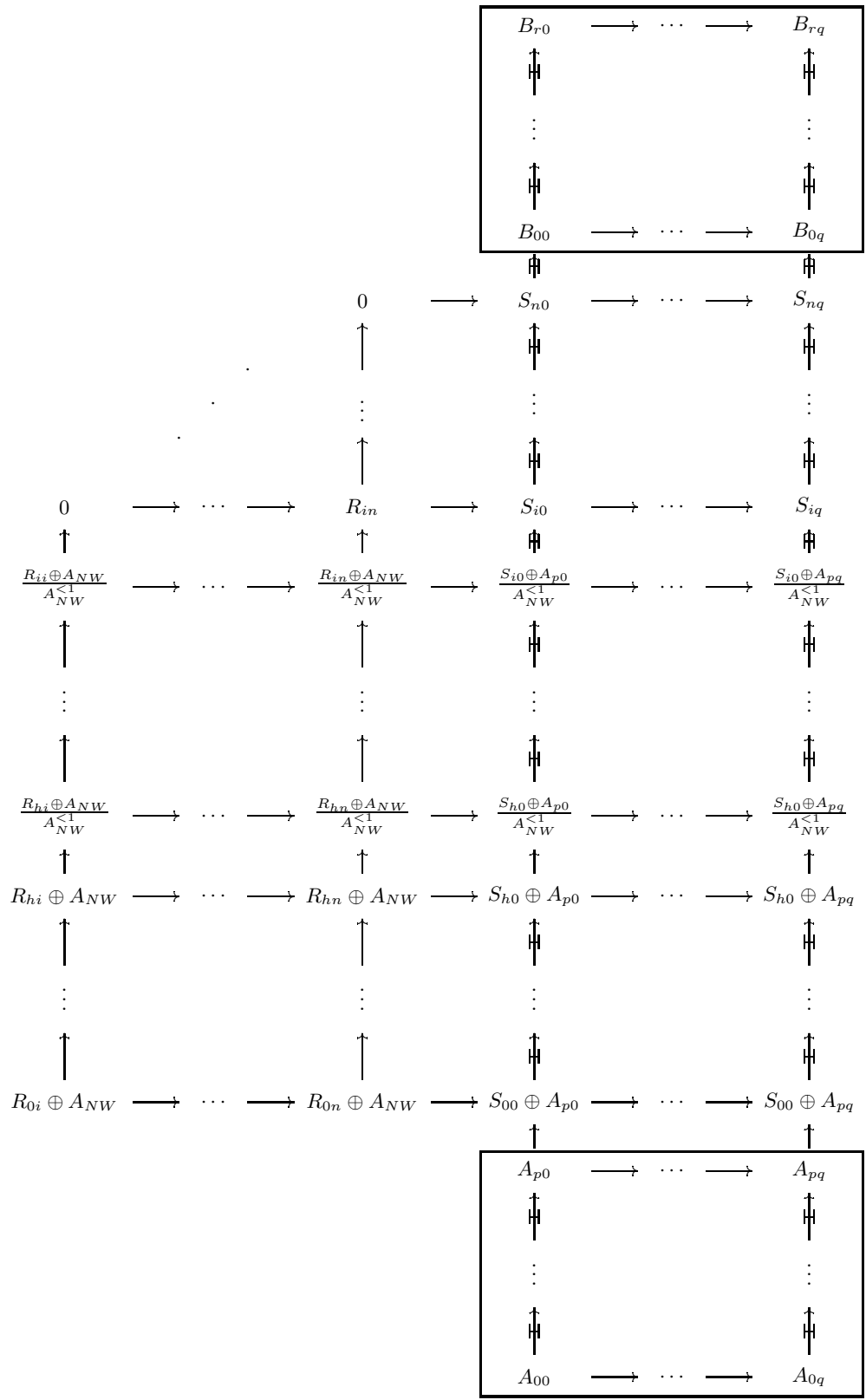
\eject
\vsize=8.4truein
\small

\siv00{19}2{17}0000
\sjv{11}1{11}1{13}1{11}1{11}
\A10i
\A3in
\A50q
\B10p
\B30h
\B5hi
\B7in
\B90r

\SB\boxten{
                                      \Sbxp59ha{B_\Q}
                                      \Sbu58h
                    \Strup37aa{R_\Q}\Sbr47a \Sbxp57ha{S_\Q}
                    \Sbu36a                   \Sbu56h
                    \Sbxp35aa{{R_\Q\oplus A_{NW}}\over{A_{NW}^{<1}}}  \Sbr45a
                           \Sbxp55ha{{S_{\I 0}\oplus A_{p\J}}\over{A_{NW}^{<1}}}
                    \Sbu34a                   \Sbu54a
                   \Sbxp33aa{R_\Q\oplus A_{NW}}\Sbr43a 
                              \Sbxp53ha{S_{\I 0}\oplus A_{p\J}}
                                           \Sbu52a
                                         \Sbxp51ha{A_\Q}
\fmbx5511 \fmbx5599 
}
\BDI(0,610)
\UB(,305)\boxten
\EDI

\normalsize

\nin
We remind the reader that, in our shorthand, this homotopy would be denoted
something like

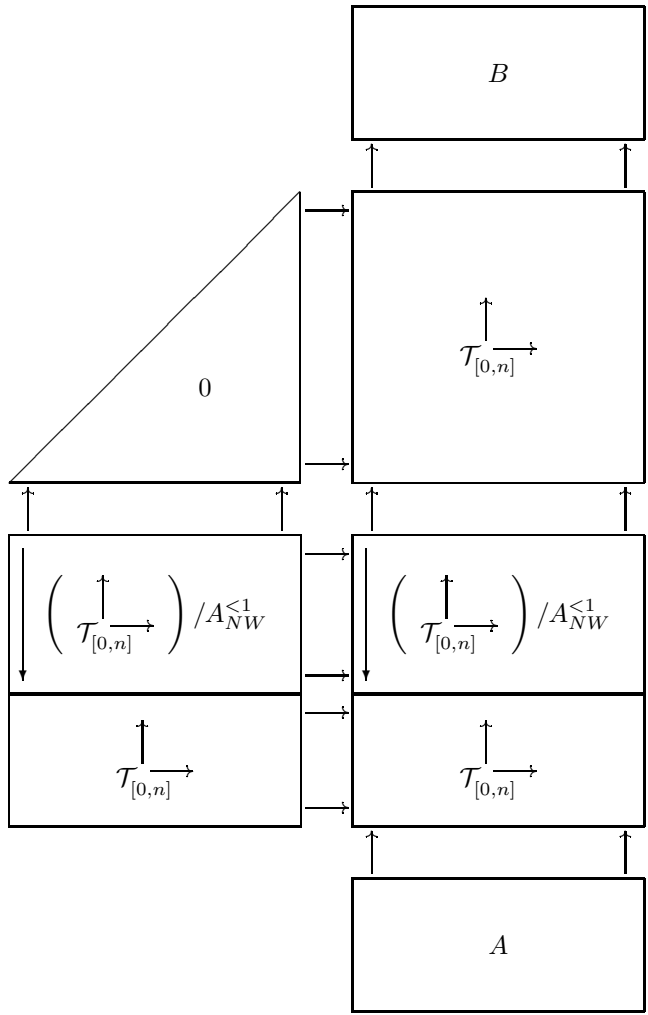
\siv{11}2{11}000000
\sjv52562{11}250

\BP
                     \Fbx38B
                     \bu37a
\trui160  \br26a     \bx36aa0n
\bu15a               \bu35a
\Fhbx14{\left(\BA \ct^{\barb{a}a}_{[0,n]}\EA\right)/A_{NW}^{<1}}d\br24a
        \Fhbx34{\left(\BA \ct^{\barb{a}a}_{[0,n]}\EA\right)/A_{NW}^{<1}}d
\bx13aa0n \br23a  \bx33aa0n
                     \bu32a   
                     \Fbx31A

\EP

\nin   
The last cell of this homotopy has, for one of its faces,   
the following diagram

\siv00{19}2{17}0000
\sjv{13}2{13}2{13}2{13}00
\A10i
\A3in
\A50q
\B10p
\B30i
\B5in
\B70r

\SB\boxten{
                                      \Sbxp57ha{B_\Q}
                                      \Sbu56h
                    \Strup35aa{R_\Q}\Sbr45a \Sbxp55ha{S_\Q}
                    \Sbu34a                   \Sbu54h
                   \Sbxp33aa{{R_\Q\oplus A_{NW}}\over{A_{NW}^{<1}}}
\Sbr43a 
                       \Sbxp53ha{{S_{\I 0}\oplus A_{p\J}}\over{A_{NW}^{<1}}}
                                           \Sbu52h
                                         \Sbxp51ha{A_\Q}
\fmbx5511 \fmbx5577 
}
\BDI(0,580)
\UB(0,290)\boxten
\EDI

\nin
which, lo and behold, is almost exactly what we need. It is practically exactly a 
cell in the homotopy whose shorthand is 
 
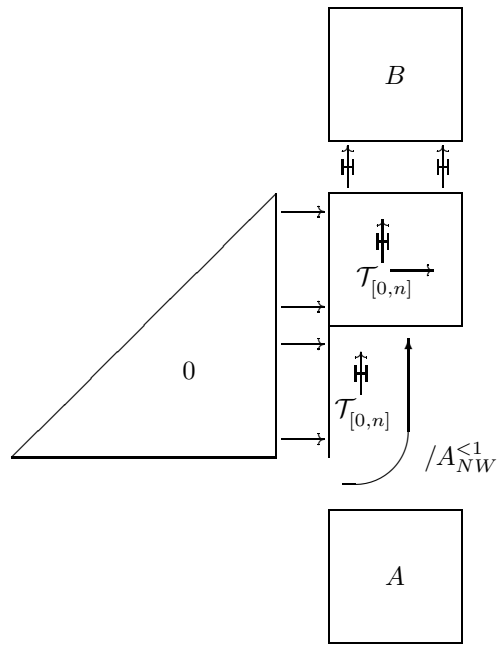
\sjv525525000
\siv{10}25000000

\SB\boxone{
                     \Fbx36B
                     \bu35h
\trui130  \br24a     \bx34ha0n
          \br23a     \htp33h0nu   \postr33{/A_{NW}^{<1}}e
                     \Fbx31A
}

\BC
\usebox\boxone
\EC

\nin
Precisely, in the diagram

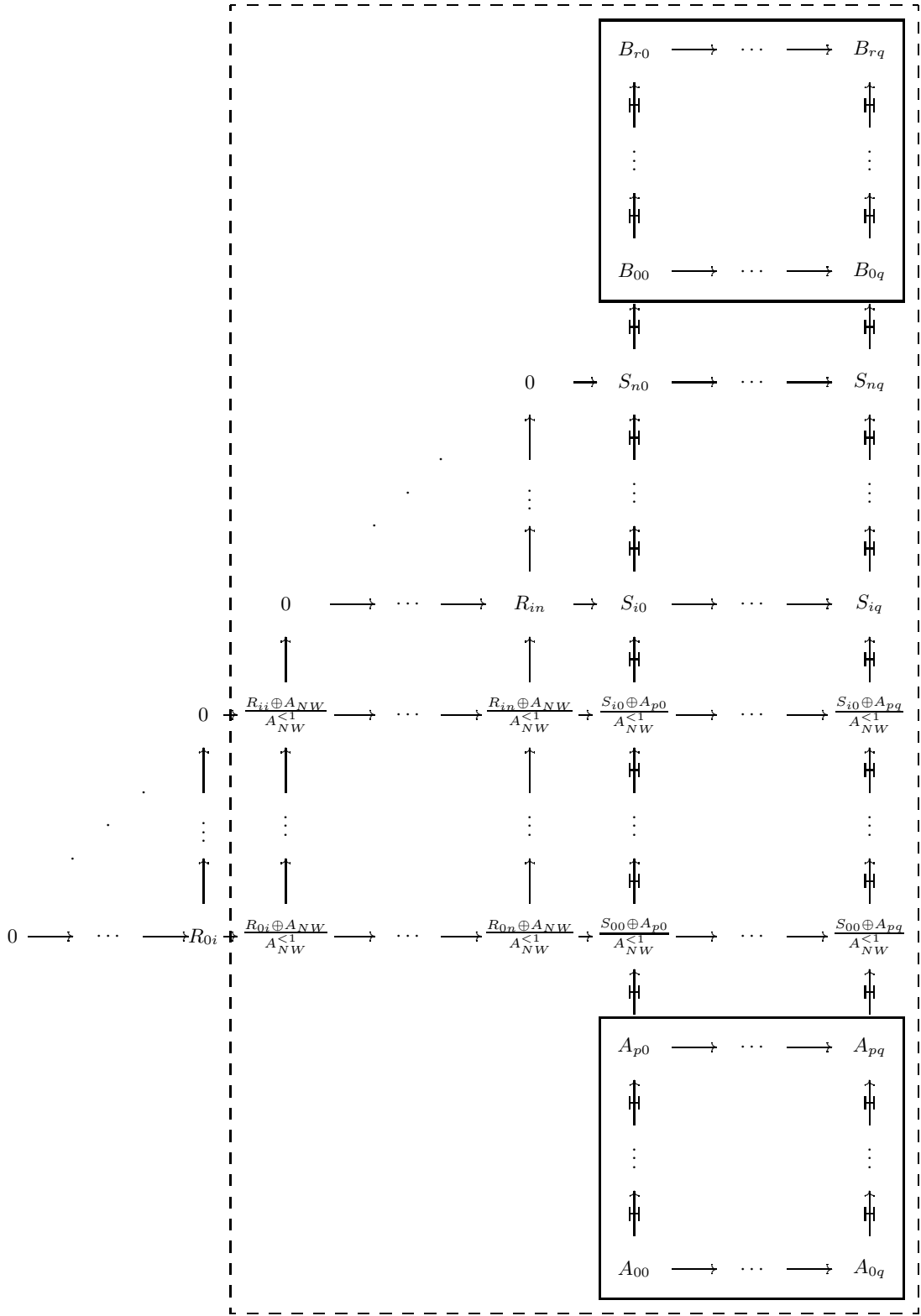
\small
\siv{10}1{15}1{14}0000
\sjv{13}2{13}2{13}2{13}00
\A10i
\A3in
\A50q
\B10p
\B30i
\B5in
\B70r

\SB\boxten{
                                      \Sbxp57ha{B_\Q}
                                      \Sbu56h
                    \Strup35aa{R_\Q}\Sbr45a \Sbxp55ha{S_\Q}
                    \Sbu34a                   \Sbu54h
\Strup13aa{R_\Q}\Sbrsmall23a   
\Sbxpsmall33aa{\scriptstyle {R_\Q\oplus A_{NW}}\over{A_{NW}^{<1}}}
\Sbrsmall43a 
\Sbxpsmall53ha{\scriptstyle {S_{\I 0}\oplus A_{p\J}}\over{A_{NW}^{<1}}}
                                           \Sbu52h
                                         \Sbxp51ha{A_\Q}
\fmbx5511 \fmbx5577  \dshbx3517
}
\BDI(0,580)
\UB(0,290)\boxten
\EDI

\normalsize

\nin
which is a typical cell in the homotopy

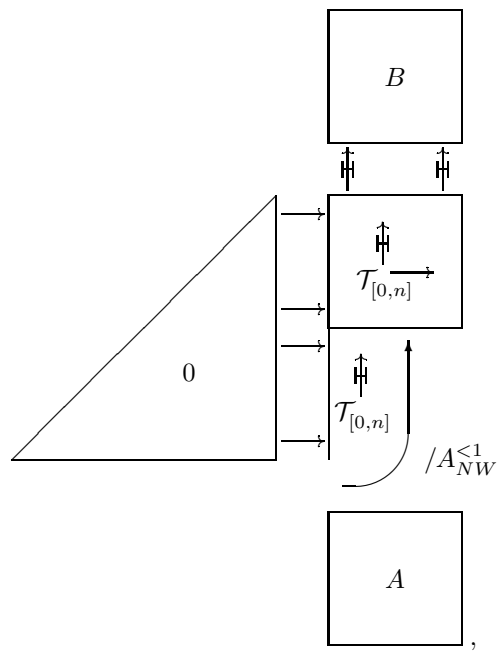
\BC
\usebox\boxone\ ,
\EC

\nin
we have shown that everything inside the dashbox is OK. All the squares are
naturally Mayer--Vietoris. But any square in the diagram is contained in
the union of the dashbox above and the dashbox below.

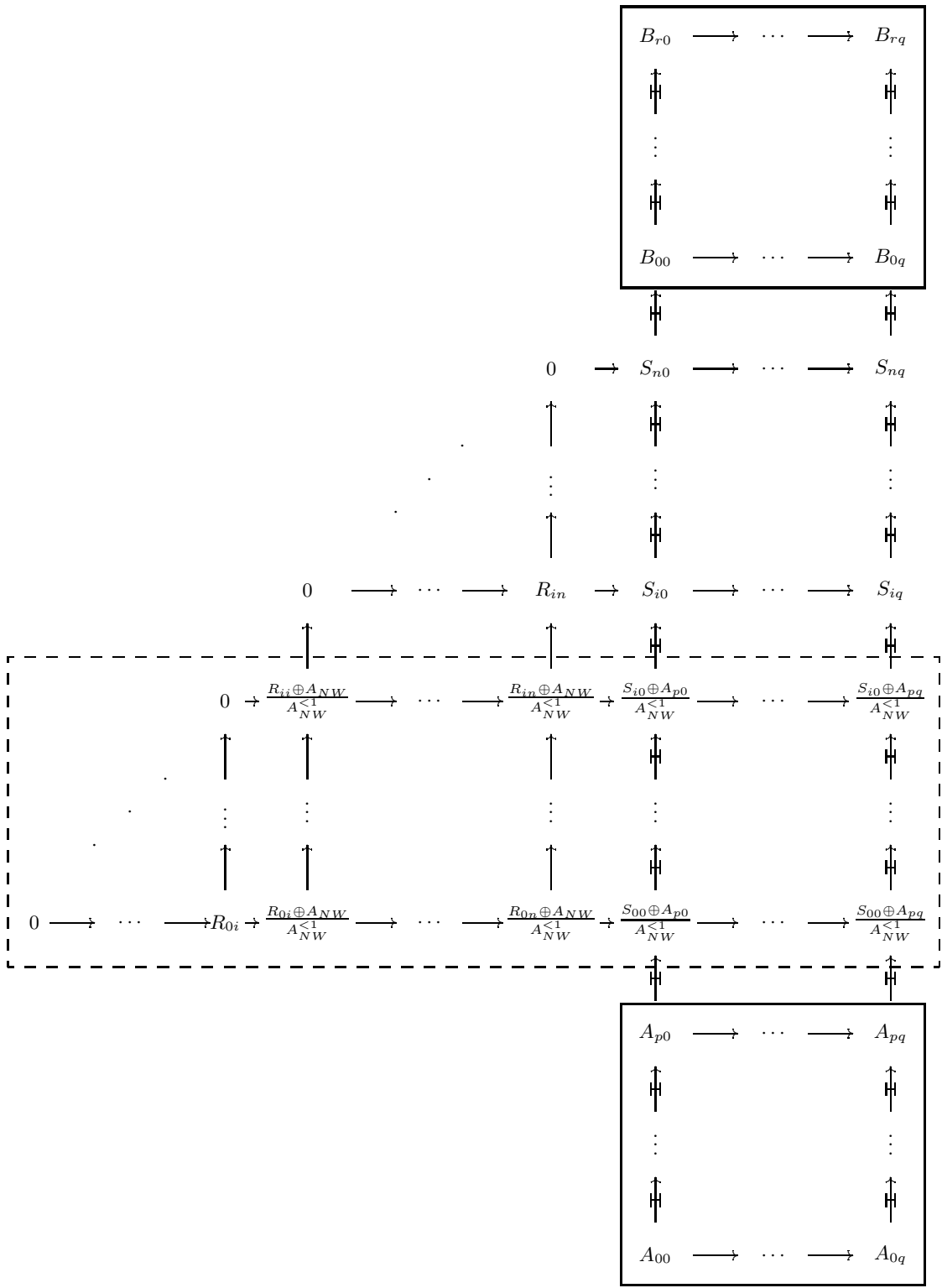
\small
\siv{10}1{15}1{14}0000
\sjv{13}2{13}2{13}2{13}00
\A10i
\A3in
\A50q
\B10p
\B30i
\B5in
\B70r

\SB\boxten{
                                      \Sbxp57ha{B_\Q}
                                      \Sbu56h
                    \Strup35aa{R_\Q}\Sbr45a \Sbxp55ha{S_\Q}
                    \Sbu34a                   \Sbu54h
\Strup13aa{R_\Q}\Sbrsmall23a   
\Sbxpsmall33aa{{R_\Q\oplus A_{NW}}\over{A_{NW}^{<1}}}
\Sbrsmall43a 
\Sbxpsmall53ha{{S_{\I 0}\oplus A_{p\J}}\over{A_{NW}^{<1}}}
                                           \Sbu52h
                                         \Sbxp51ha{A_\Q}
\fmbx5511 \fmbx5577  \dshbx1533
}
\BDI(0,580)
\UB(0,290)\boxten
\EDI

\normalsize

\nin
It suffices therefore to show that

\small
\siv{10}1{15}1{14}0000
\sjv00{13}000000
\A10i
\A3in
\A50q
\B10p
\B30i
\B5in
\B70r

\SB\boxten{
\Strup13aa{R_\Q}\Sbrsmall23a   
\Sbxpsmall33aa{{R_\Q\oplus A_{NW}}\over{A_{NW}^{<1}}}
\Sbrsmall43a 
                 \Sbxpsmall53ha{{S_{\I 0}\oplus A_{p\J}}\over{A_{NW}^{<1}}}
}
\BDI(0,130)
\UB(0,65)\boxten
\EDI

\normalsize

\nin
is a simplex (i.e. all squares are naturally Mayer--Vietoris). But

\small

\siv{10}1{15}1{14}0000
\sjv00{13}000000
\A10i
\A3in
\A50q
\B10p
\B30i
\B5in
\B70r

\SB\boxten{
\Strup13aa{R_\Q}\Sbrsmall23a   
\Sbxpsmall33aa{\scriptstyle R_\Q\oplus A_{NW}} \Sbrsmall43a 
                    \Sbxpsmall53ha{\scriptstyle S_{\I 0}\oplus A_{p\J}}
}
\BDI(0,130)
\UB(0,65)\boxten
\EDI

\normalsize

\nin
certainly is, being part of the well--defined homotopy cell

\small

\siv{10}1{15}1{14}0000
\sjv{13}2{13}2{13}2{13}00
\A10i
\A3in
\A50q
\B10p
\B30i
\B5in
\B70r

\SB\boxten{
                                      \Sbxp57ha{B_\Q}
                                      \Sbu56h
                    \Strup35aa{R_\Q}\Sbr45a \Sbxp55ha{S_\Q}
                    \Sbu34a                   \Sbu54a
\Strup13aa{R_\Q}\Sbrsmall23a
\Sbxpsmall33aa{\scriptstyle R_\Q\oplus A_{NW}}\Sbrsmall43a 
                  \Sbxpsmall53ha{\scriptstyle S_{\I 0}\oplus A_{p\J}}
                                           \Sbu52a
                                         \Sbxp51ha{A_\Q}
\fmbx5511 \fmbx5577 \dshbx1533
}
\BDI(0,580)
\UB(0,290)\boxten
\EDI

\normalsize 

\nin
and once again the cell

\small

\siv{10}1{15}1{14}0000
\sjv00{13}000000
\A10i
\A3in
\A50q
\B10p
\B30i
\B5in
\B70r

\SB\boxten{
\Strup13aa{R_\Q}\Sbr23a   \Sbxp33aa{{R_\Q\oplus A_{NW}}\over{A_{NW}^{<1}}}
\Sbr43a 
                                         \Sbxp53ha{{S_{\I 0}\oplus A_{p\J}}\over{A_{NW}^{<1}}}
}
\BDI(0,130)
\UB(0,65)\boxten
\EDI

\normalsize

\nin
is obtained from the above by a trivial truncation homotopy. It is the homotopy 
which, in our shorthand, has the symbol

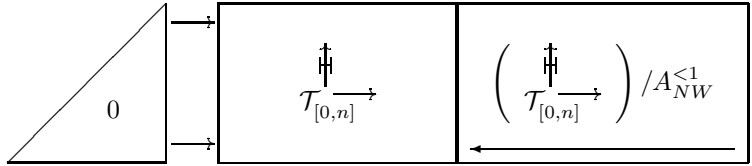
\siv629{11}00000
\sjv600000000

\BP
\trui110 \br21a \bx31ha0n \Fhbx41{\left(\BA 
\ct^{\barb{h}a}_{[0,n]}\EA\right)/A_{NW}^{<1}}l
\EP

\rmk{group2 trunc} Now we have to discuss why the above homotopy
takes cells with
a lifting to some model category 
$C^b(\cal Q)$, to bigger cells which also lift.
The problem is related to the fact that it is not entirely clear, 
why \tst\ truncations should
take a diagram with a lifting to another diagram with a lifting. 
In other words,
because the \tst\ is not assumed to be in any way related to the model category
$C^b({\cal Q})$ to which the simplex lifts, 
why should one expect the lifting to
be preserved by truncation?

\siv500000000
\sjv500000000

\SB\boxele{
\fbx11\ct aa}

\siv700000000

\SB\boxtwe{
\fbx11{C^b({\cal Q})}aa}

Consider now a simplex in $\PT\boxele$ with a lifting to $\PT\boxtwe.$
That is, we have a diagram of bicartesian squares in $C^b({\cal Q}):$

\siv{13}00000000
\sjv{13}00000000

\A10n
\B10m

\BP
\Sbxp11aa{B_\Q}
\EP

\nin
Of course, this diagram is entirely determined by pullback from the diagram

\BP
\uSbxp11aa{B_\Q}
\EP

\nin
What we want to do in this Remark, is reinterpret this obvious fact in 
terms of sheaves on some space.

Consider the partially ordered set $\cal B$, given below:

\siv{12}00000000
\sjv{11}00000000

\BP
\uSbxp11aa{b_\Q}
\EP

\nin
That is, the elements of the set $\cal B$ are 

\barr
b_{ij} &\quad\hbox{ where }\quad & \left\{\begin{array}{lcr}
i=m & \hbox{ and } &0\leq j\leq n\\
    & \hbox{\bf or } &           \\
0\leq i\leq m  & \hbox{ and } & j=n\end{array}\right.
\earr

\nin
The partial ordering is given by setting $b_{ij}<b_{i'j'}$ if there is an
arrow in the diagram from  $b_{ij}$ to $b_{i'j'}.$

Now we make $\cal B$ into a topological space. The open subsets of $\cal B$
are the subsets closed under majorisation. Precisely, $U\subset{\cal B}$ is
open if, whenever $x\in U$ and $x<y,$ then $y\in U.$

What is a sheaf over $\cal B$ with values in the abelian category $\cal Q?$
It is nothing more nor less than a diagram

\siv{13}00000000
\sjv{13}00000000

\A10n
\B10m

\BP
\Sbxp11aa{B_\Q}
\EP

\nin
where each $B_{ij}\in{\cal Q},$ and every square is a pullback square. After all,
the open sets in $\cal B$ are all of the form ${\cal B}_{ij},$
for some pair of integers $0\leq i\leq m,$
and $0\leq j\leq n.$ ${\cal B}_{ij}$ just consists of the elements

\siv{12}00000000
\sjv{11}00000000

\A1jn
\B1im

\BP
\uSbxp11aa{b_\Q}
\EP

\nin
Clearly, ${\cal B}_{ij}={\cal B}_{mj}\cup{\cal B}_{in},$ and therefore to give
an arbitrary sheaf $\cal S$ on $\cal B$, it suffices to give its values on the
open sets ${\cal B}_{mj}$  for $0\leq j\leq n,$ and ${\cal B}_{in}$ for
$0\leq i\leq m.$ Let $B_{ij}$ stand for $\Gamma({\cal B}_{ij},{\cal S}).$
The diagram

\siv{13}00000000
\sjv{13}00000000

\A10n
\B10m

\BP
\uSbxp11aa{B_\Q}
\EP

\nin
determines the values of the sheaf $\cal S$ on the open sets ${\cal B}_{mj}$
and ${\cal B}_{in},$ which form a basis for the topology. The value on 
${\cal B}_{ij}$ is determined by the pullback diagram

\siv343000000
\sjv343000000

\BP
\postr13{B_{mj}}c \Sposar23ar  \postr33{B_{mn}}c
\Sposar12au                    \Sposar32au
\postr11{B_{ij}}c \Sposar21ar  \postr31{B_{in}}c 
\EP

\nin
and hence a sheaf on this topological space really is nothing other than
a diagram

\siv{13}00000000
\sjv{13}00000000

\A10n
\B10m

\BP
\Sbxp11aa{B_\Q}
\EP

\nin
where every square is a pullback square.

The topological space $\cal B$ and all of its open subsets are finite sets.
In particular, they are paracompact,
 and the \u{C}ech cohomology agrees with
the derived functor cohomology. Given a sheaf $\cal S$ on ${\cal B},$ that
is a diagram

\BP
\Sbxp11aa{B_\Q}
\EP

\nin
one can reasonably ask to compute the cohomology of this sheaf. Given an
open set ${\cal B}_{ij}\subset{\cal B},$ what are the groups 
$H^k({\cal B}_{ij},{\cal S})$?

It is immediate that every cover of ${\cal B}_{ij}$ contains the cover
${\cal B}_{ij}={\cal B}_{mj}\cup{\cal B}_{in}$ as a refinement. The \u{C}ech
cohomology, which agrees with the derived functor cohomology, can therefore
be computed on this cover. We deduce that the cohomology 
$H^k({\cal B}_{ij},{\cal S})$ is just the cohomology of the complex

$$
B_{in}\oplus B_{mj}\longrightarrow B_{mn}.
$$

\nin
Therefore, to say that the diagram

\BP
\Sbxp11aa{B_\Q}
\EP

\nin
consists of {\em bicartesian} squares, is nothing more nor less than to require
that the sheaf $\cal S$ have vanishing higher cohomology on every open subset
${\cal B}_{ij}\subset{\cal B}.$ Let us call a sheaf $\cal S$ {\em acyclic} if its
cohomology vanishes on every ${\cal B}_{ij}\subset{\cal B}.$

Let $D^b({\cal B})$ be the bounded derived category of the category of 
all sheaves on $\cal B$ with values in $\cal Q.$ The objects of $D^b({\cal B})$
are complexes of sheaves $\cal S$ on $\cal B.$ Every such complex is
quasi--isomorphic to a complex of acyclic sheaves. Take for instance a flabby
resolution. To begin with, we know we can always replace a complex by a (possibly 
infinite) flabby resolution. But in fact, the resolution 
by acyclics may be taken finite,
because every one of the finitely many open sets in $\cal B$ has finite
homological dimension.

Thus $D^b({\cal B})$ is the derived category of a category we denote $\PT\boxtwe$,
whose objects are diagrams

\BP
\Sbxp11aa{B_\Q}
\EP

\nin
of bicartesian squares in $C^b({\cal Q}).$ Now the real point of the Remark is
that $\cal B$ is a stratified space, with every point being a stratum. Thus,
given an arbitrary perversity function on the strata, there corresponds a
\tst\ on $D^b({\cal B}).$ In particular, there corresponds a \tst\
truncation. And the point is that the truncated simplices used
in the proof of Theorem~\ref{safe kernels} are all obtainable from
simplices in $\PT\boxtwe$ by a \tst\ truncation for some perverse \tst\ on the
category  $D^b({\cal B}).$

\begin{description}

\item[Note.] It is slightly criminal to denote our category by $\PT\boxtwe,$
because we tend to think of $\PT\boxtwe$ as a simplicial set, containing diagrams

\BP
\Sbxp11aa{B_\Q}
\EP

where $m$ and $n$ are allowed to vary. Our category
has objects in which $m$ and $n$ are fixed. It would perhaps be 
better notation to denote the category
${\mtrx\boxtwe}_{mn}.$ However, the author happens to find this notation so 
\ae sthetically revolting, that he preferred the slightly misleeding notation.

\end{description}

Now I have told the reader in general terms what we will do. It is time to turn
very specific. We need to show how to construct the cell

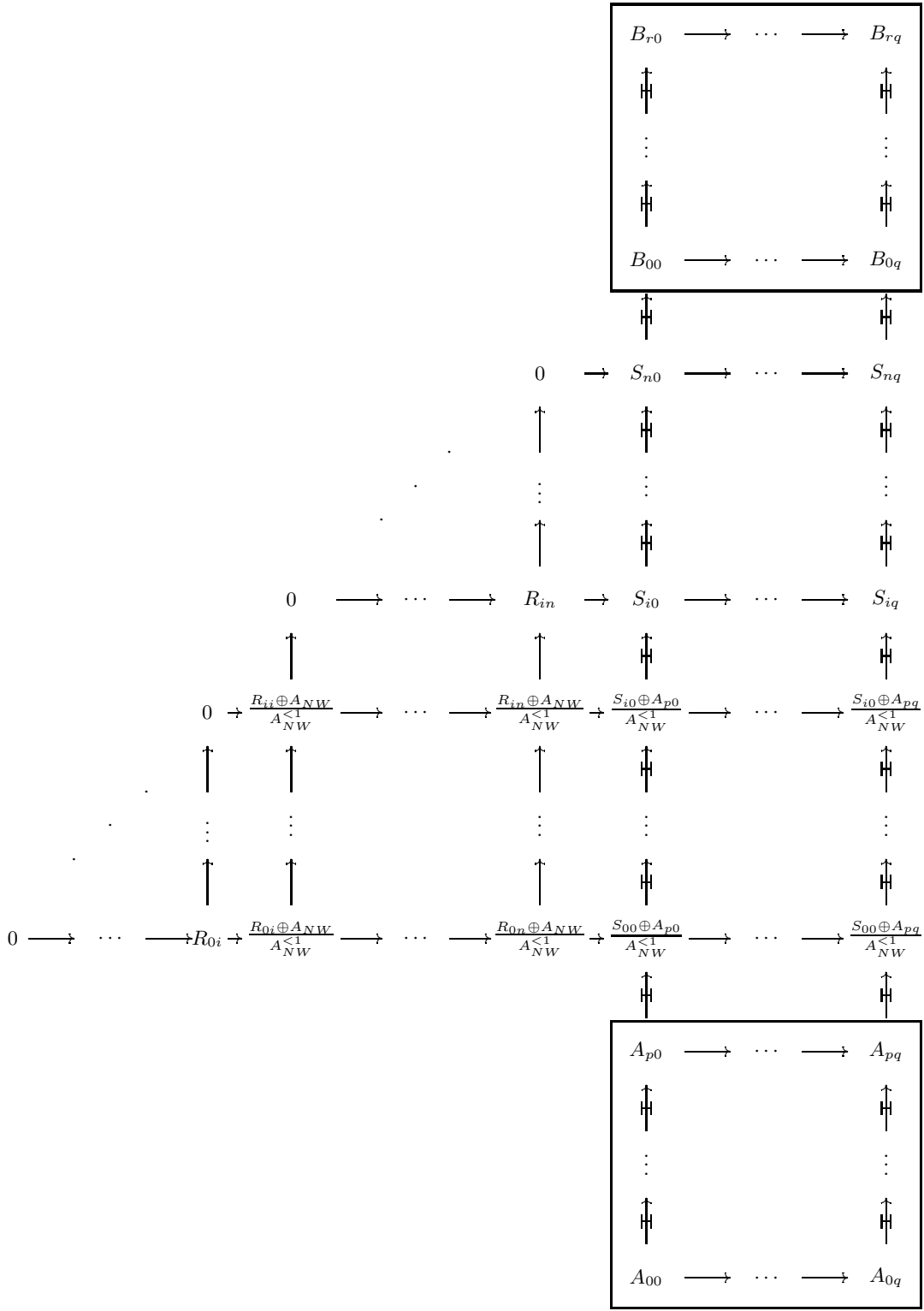
\small

\siv{10}1{15}1{14}0000
\sjv{13}2{13}2{13}2{13}00
\A10i
\A3in
\A50q
\B10p
\B30i
\B5in
\B70r

\SB\boxten{
                                      \Sbxp57ha{B_\Q}
                                      \Sbu56h
                    \Strup35aa{R_\Q}\Sbr45a \Sbxp55ha{S_\Q}
                    \Sbu34a                   \Sbu54h
\Strup13aa{R_\Q}\Sbrsmall23a   
\Sbxpsmall33aa{{R_\Q\oplus A_{NW}}\over{A_{NW}^{<1}}}
\Sbrsmall43a 
                       \Sbxpsmall53ha{{S_{\I 0}\oplus A_{p\J}}\over{A_{NW}^{<1}}}
                                           \Sbu52h
                                         \Sbxp51ha{A_\Q}
\fmbx5511 \fmbx5577  
}
\BDI(0,580)
\UB(0,290)\boxten
\EDI

\normalsize

\nin
out of the cell

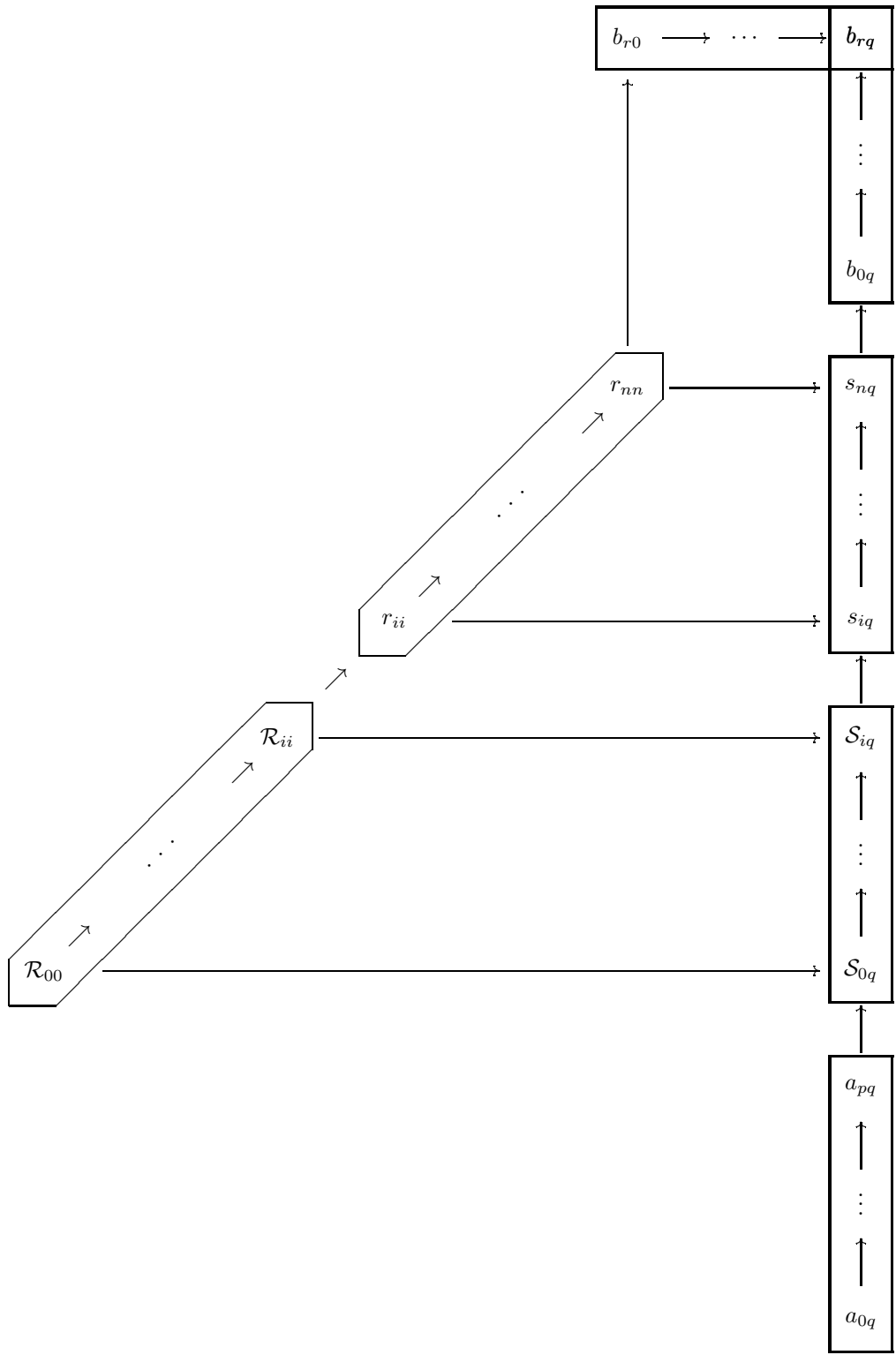
\small

\siv{10}1{15}1{14}0000
\sjv{13}2{13}2{13}2{13}00
\A10i
\A3in
\A50q
\B10p
\B30i
\B5in
\B70r

\SB\boxten{
                                      \Sbxp57ha{B_\Q}
                                      \Sbu56h
                    \Strup35aa{R_\Q}\Sbr45a \Sbxp55ha{S_\Q}
                    \Sbu34a                   \Sbu54a
\Strup13aa{R_\Q}\Sbrsmall23a
\Sbxpsmall33aa{\scriptstyle R_\Q\oplus A_{NW}}\Sbrsmall43a 
                      \Sbxpsmall53ha{\scriptstyle S_{\I 0}\oplus A_{p\J}}
                                           \Sbu52a
                                         \Sbxp51ha{A_\Q}
\fmbx5511 \fmbx5577
}
\BDI(0,580)
\UB(0,290)\boxten
\EDI

\normalsize

\nin
using nothing other than perverse \tst\ truncations on suitably chosen 
finite topological spaces. To begin with, observe that to give a lifting
of

\small

\siv{10}1{15}1{14}0000
\sjv{13}2{13}2{13}2{13}00
\A10i
\A3in
\A50q
\B10p
\B30i
\B5in
\B70r

\SB\boxten{
                                      \Sbxp57ha{B_\Q}
                                      \Sbu56h
                    \Strup35aa{R_\Q}\Sbr45a \Sbxp55ha{S_\Q}
                    \Sbu34a                   \Sbu54a
\Strup13aa{R_\Q}\Sbrsmall23a
\Sbxpsmall33aa{\scriptstyle R_\Q\oplus A_{NW}}\Sbrsmall43a 
                  \Sbxpsmall53ha{\scriptstyle S_{\I 0}\oplus A_{p\J}}
                                           \Sbu52a
                                         \Sbxp51ha{A_\Q}
\fmbx5511 \fmbx5577
}
\BDI(0,580)
\UB(0,290)\boxten
\EDI

\normalsize

\nin
\eject
\vsize=10in
is the same as to give a diagram $C^b({\cal Q})$

\small

\siv{10}1{15}1{14}0000
\sjv{13}2{13}2{13}2{13}00
\A10i
\A3in
\A50q
\B10p
\B30i
\B5in
\B70r

\SB\boxten{
                                      \Sbxp57ha{{\tilde B}_\Q}
                                      \Sbu56h
                    \fStrup35aa{{\tilde R}_\Q}\Sbr45a \Sbxp55ha{{\tilde S}_\Q}
                    \Sbu34a                   \Sbu54a
\fStrup13aa{{\tilde R}_\Q}\Sbrsmall23a
\Sbxpsmall33aa{\scriptstyle{\tilde R}_\Q\oplus {\tilde A}_{NW}}\Sbrsmall43a 
             \Sbxpsmall53ha{\scriptstyle{\tilde S}_{\I 0}\oplus\tilde A_{p\J}}
                                           \Sbu52a
                                         \Sbxp51ha{{\tilde A}_\Q}
\fmbx5511 \fmbx5577
}
\BDI(0,580)
\UB(0,290)\boxten
\EDI

\normalsize

\nin
where the ${\tilde R}_{ii}$ happen to be 
isomorphic to zero in $D^b({\cal Q}).$ To be
more precise, they come with an isomorphism to zero, and all other objects in the
lifted diagram come with isomorphisms to the objects they lift, and the
isomorphisms commute with all the maps that are part of the definition of a 
simplex. 

This diagram can be viewed as a sheaf on the topological space ${\cal X}$,
which we will 
schematically denote
/eject
/vsize=8.4truein
\siv{13}2{10}0{10}3000
\sjv{13}2{13}2{13}2{13}00
\A10i
\A3in
\A50q
\B10p
\B30i
\B5in
\B70r
\SB\boxtwe{
\Sulip67a{b_{\I q}} \fmbx6677
\Sposar66au
\Sulip65a{s_{\I q}} \fmbx6655
\Sposar64au
\Sulip63a{{\cal S}_{\I q}} \fmbx6633
\Sposar62au
\Sulip61a{a_{\I q}} \fmbx6611
}
\siv{13}2{10}0{13}0000
\sjv{13}2{13}2{13}2{10}30
\SB\boxele{
\Srlip58a{b_{r\J}}  \fmbx5588
}
\siv3732373{7}3
\sjv{15}3732373{15}
\SB\boxthi{
\postr12{{\cal R}_{00}}c 
\postr34{{\cal R}_{ii}}c
\postr56{{r}_{ii}}c
\postr78{{r}_{nn}}c  
\postr45\nearrow c
\postr12\nearrow n
\postr34\nearrow s
\postr56\nearrow n
\postr78\nearrow s
}
\siv6162616{7}3
\sjv{15}6162616{15}
\SB\boxten{
\postr23\cdot c
\postr23\cdot n
\postr23\cdot s
\postr67\cdot c
\postr67\cdot n
\postr67\cdot s
}
\SB\boxnin{
\put(0,150){\line(1,0){20}}
\put(0,150){\line(0,1){20}}
\put(130,280){\line(-1,0){20}}
\put(130,280){\line(0,-1){20}}
\put(20,150){\line(1,1){110}}
\put(0,170){\line(1,1){110}}
\put(150,300){\line(1,0){20}}
\put(150,300){\line(0,1){20}}
\put(280,430){\line(-1,0){20}}
\put(280,430){\line(0,-1){20}}
\put(170,300){\line(1,1){110}}
\put(150,320){\line(1,1){110}}
\UB(190,290)\boxele
\UB(190,290)\boxtwe
\UB(190,290)\boxthi
\UB(190,290)\boxten
\setua{265}{433}{114}a
\setra{283}{415}{64}a
\setra{133}{265}{214}a
\setra{40}{165}{307}a
\setra{190}{315}{157}a
}
\BDI(0,580)
\UB(0,290)\boxnin
\EDI

\nin
The way this diagram should be read is the following. The partially ordered set
has 
\eject
\vsize=10in
\noindent
points $b_{rj}$ where $0\leq j\leq q,$ $b_{jq}$
with $0\leq j\leq r$, $r_{jj}$
with $i\leq j\leq r$,
${\cal R}_{jj}$
with $0\leq j\leq i$,
$s_{jq}$
with $i\leq j\leq r$,
${\cal S}_{jq}$
with $0\leq j\leq i$ and
$a_{jq}$
with $0\leq \gamma\leq p.$ The labels of the points are meant to be suggestive;
each point corresponds to the object having the same name but with Roman capitals,
in the diagram of the simplex. The only exception is ${\cal S}_{jq}$,
which corresponds to $S_{j0}\oplus A_{pq}$.

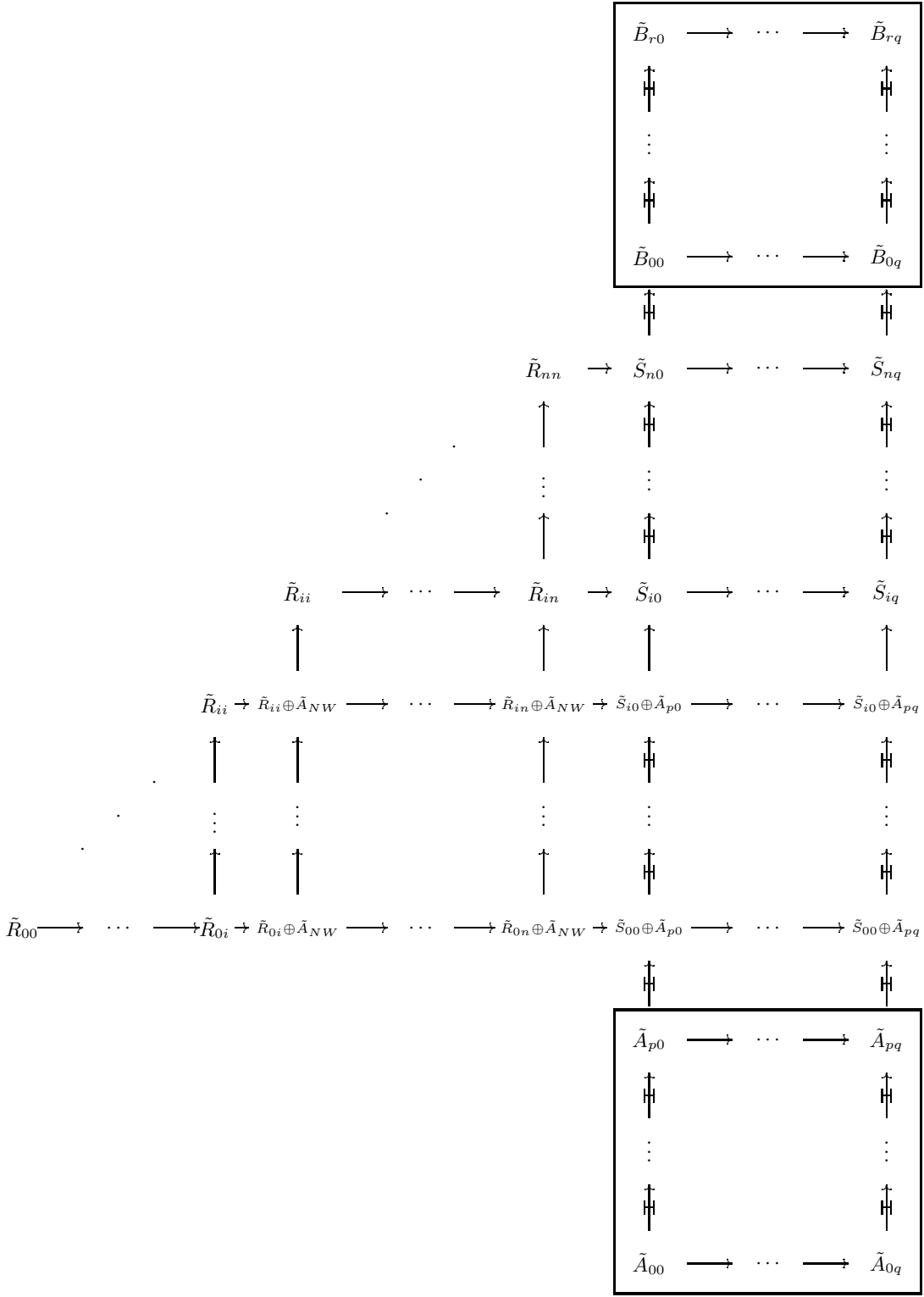
\small

\siv{10}1{15}1{14}0000
\sjv{13}2{13}2{13}2{13}00
\A10i
\A3in
\A50q
\B10p
\B30i
\B5in
\B70r

\SB\boxten{
                                      \Sbxp57ha{{\tilde B}_\Q}
                                      \Sbu56h
                    \fStrup35aa{{\tilde R}_\Q}\Sbr45a \Sbxp55ha{{\tilde S}_\Q}
                    \Sbu34a                   \Sbu54a
\fStrup13aa{{\tilde R}_\Q}\Sbrsmall23a
\Sbxpsmall33aa{\scriptstyle {\tilde R}_\Q\oplus {\tilde A}_{NW}}\Sbrsmall43a 
           \Sbxpsmall53ha{\scriptstyle {\tilde S}_{\I 0}\oplus \tilde A_{p\J}}
                                           \Sbu52h
                                         \Sbxp51ha{{\tilde A}_\Q}
\fmbx5511 \fmbx5577
}
\BDI(0,580)
\UB(0,290)\boxten
\EDI

\normalsize

\nin
The partial ordering is that $b$'s are as in the example discussed above,
 the $r$'s are smaller than the $b$'s,
the ${\cal R}$'s are smaller than the $r$'s, the $s$'s are smaller than the
$b_{iq}$'s, the $\cal S$'s are smaller than  the $s$'s, and the $a$'s
are smaller than the $\cal S$'s. Furthermore, for every $j$,
$r_{jj}<s_{jq}$ and ${\cal R}_{jj}<
 \cs_{jq}.$ Perhaps a simpler way to state this, is to say 
that the points of the partially ordered set are all the objects of the diagram
on the right and top fringe, together with the 
diagonal entries of the triangles
on the 
\eject
\vsize=19truein
\noindent
left (these are declared to be honorary members of the top fringe).
A point in this partially ordered set is less 
than another, if there is an arrow
joining the lesser to the greater.

Once again, one can turn this partially ordered set into a topological space, 
by declaring sets closed under majorisation 
to be open. I leave it to the reader
to check that an acyclic sheaf on this space is nothing more nor less than a
diagram

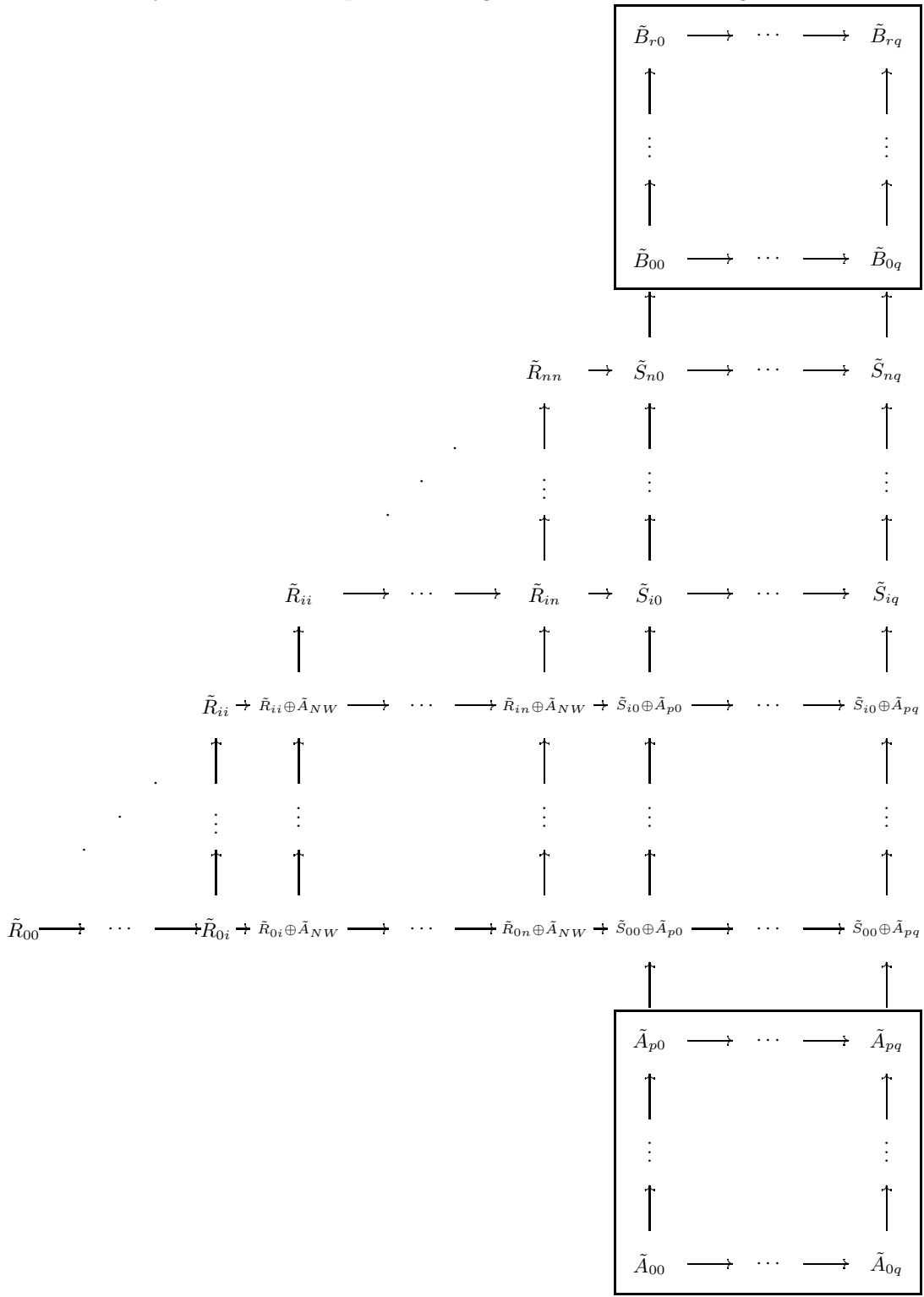
\small

\siv{10}1{15}1{14}0000
\sjv{13}2{13}2{13}2{13}00
\A10i
\A3in
\A50q
\B10p
\B30i
\B5in
\B70r

\SB\boxten{
                                      \Sbxp57aa{{\tilde B}_\Q}
                                      \Sbu56a
                    \fStrup35aa{{\tilde R}_\Q}\Sbr45a \Sbxp55aa{{\tilde S}_\Q}
                    \Sbu34a                   \Sbu54a
\fStrup13aa{{\tilde R}_\Q}\Sbrsmall23a
\Sbxpsmall33aa{\scriptstyle {\tilde R}_\Q\oplus {\tilde A}_{NW}}\Sbrsmall43a 
             \Sbxpsmall53aa{\scriptstyle {\tilde S}_{\I 0}\oplus\tilde A_{p\J}}
                                           \Sbu52a
                                         \Sbxp51aa{{\tilde A}_\Q}
\fmbx5511 \fmbx5577
}
\BDI(0,580)
\UB(0,290)\boxten
\EDI

\normalsize

\nin
of bicartesian squares of objects, in the abelian category 
$\cal Q.$ The derived
category 
\eject
\vsize=8.4truein
\noindent
$D^b({\cal X})$ has for a model a category 
which we will denote

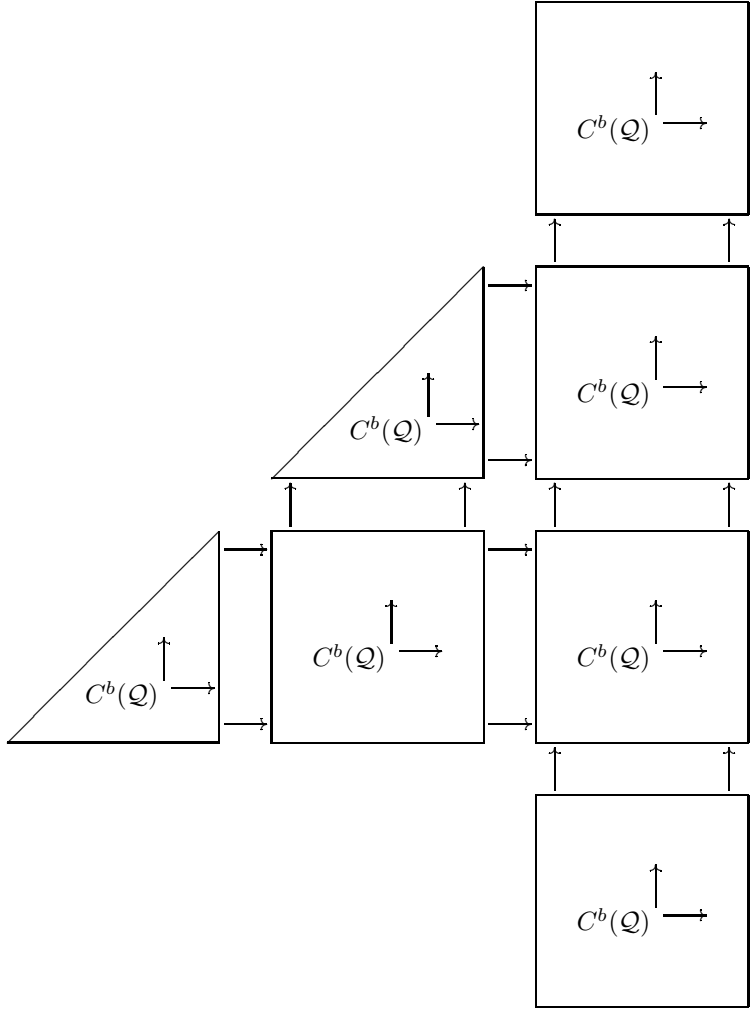
\siv828280000
\sjv828282800

\BP
                                    \fbx57{C^b({\cal Q})}aa
                                    \bu56a
             \trui35{C^b({\cal Q})^{\barb aa}} \br45a  
                               \fbx55{C^b({\cal Q})}aa
              \bu34a                                    \bu54a
\trui13{C^b({\cal Q})^{\barb aa}} \br23a  \fbx33{C^b({\cal Q})}aa
\br43a  \fbx53{C^b({\cal Q})}aa
  \bu52a
\fbx51{C^b({\cal Q})}aa
\EP

\nin
consisting of diagrams of bicartesian squares in 
$C^b({\cal Q})$, as above.

To define the \tst\ we want, we need first a 
slightly larger topological space.
We enlarge the partially ordered set $\cal X$ by 
one point, obtaining a space
we will denote ${\cal Y}.$ Let us call the extra point, 
i.e. the unique
point in ${\cal Y}-{\cal X},$ by the name $p.$ We declare $p$ to be
less than
\[
\cs_{iq}<s_{iq}<\cdots<s_{nq}<b_{0q}<\cdots<b_{rq}
\] 
and no other points. The only points less than
$p$ are the $a$'s and the $\cal R$'s.

An acyclic sheaf on $\cal Y$ with values in 
$\cal Q$ is a three dimensional
diagram of bicartesian squares. To represent it on a planar piece of 
paper, I will 
write two planar sections of it. It consists of a diagram

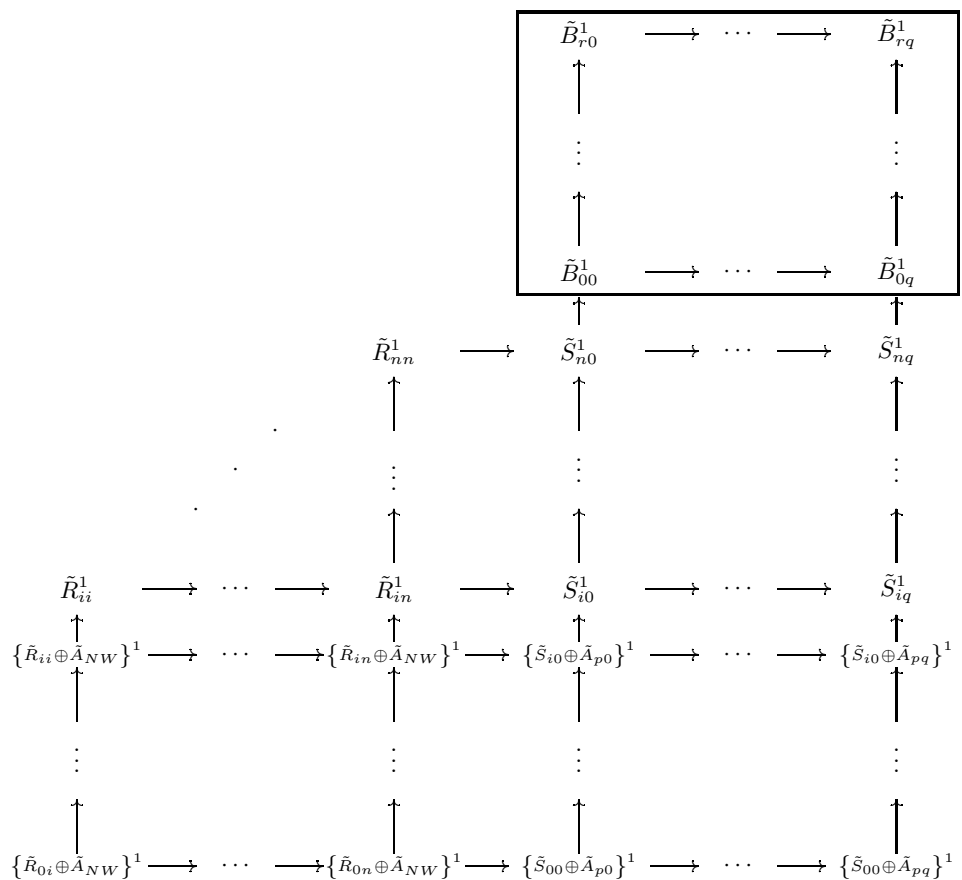
\small

\siv00{17}2{17}0000
\sjv00{9}1{11}1{11}00
\A10i
\A3in
\A50q
\B10p
\B30i
\B5in
\B70r

\SB\boxten{
                                      \Sbxp57aa{{\tilde B}^1_\Q}
                                      \Sbu56a
                    \fStrup35aa{{\tilde R}^1_\Q}\Sbr45a 
\Sbxp55aa{{\tilde S}^1_\Q}
                    \Sbu34a                   \Sbu54a
\Sbxpsmall33aa{\scriptstyle\left\{{\tilde R}_\Q\oplus 
{\tilde A}_{NW}\right\}^1}
\Sbrsmall43a 
     \Sbxpsmall53aa{\scriptstyle\left\{{\tilde S}_{\I 0}
\oplus\tilde A_{p\J}\right\}^1}
                                          \fmbx5577
}
\BDI(0,330)
\UB(0,165)\boxten
\EDI

\normalsize
\medskip

\nin
together with another diagram

\medskip

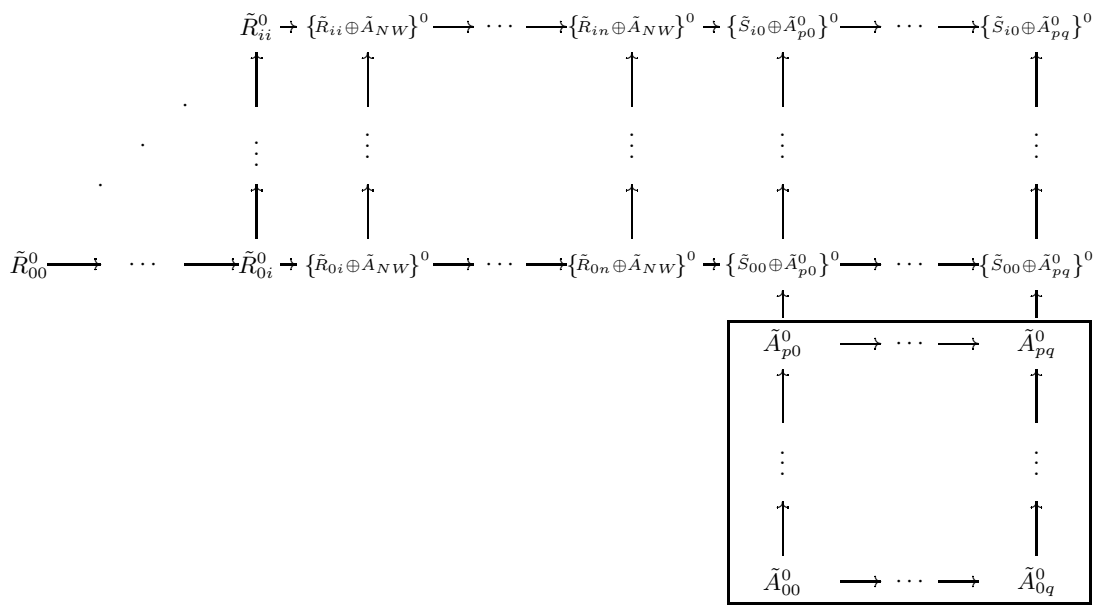
\small

\siv{10}1{15}1{14}0000
\sjv{11}1{11}000000
\A10i
\A3in
\A50q
\B10p
\B30i
\B5in
\B70r

\SB\boxele{
\fStrup13aa{{\tilde R}^0_\Q}\Sbrsmall23a
\tinySbxp33aa{\scriptstyle\left\{\!{\tilde R}_\Q\oplus {\tilde
A}_{NW}\!\right\}^0}
\Sbrsmall43a 
 \tinySbxp53aa{\scriptstyle\left\{{\tilde S}_{\I 0}
\oplus\tilde A^0_{p\J}\right\}^0}
                                           \tinySbu52a
                               \tinySbxp51aa{{\tilde A}^0_\Q}
\fmbx5511 
}
\BDI(0,230)
\UB(0,115)\boxele
\EDI

\normalsize

\nin
together with maps from any object with a 
superscript 0 to the object of the same 
\eject
\vsize=10in
\noindent
label with a superscript 1, so that every square, 
even in the third direction,
is bicartesian. Concretely, the open sets 
containing $p\in{\cal Y}-{\cal X}$
correspond to the objects with superscript 0, 
while the open sets not containing
$p$ have superscript 1.

The perversity we want to consider on $\cal Y$ is 
actually relatively simple.
Let $Z$ be the closure of $p$ in $\cal Y,$ and let 
$U={\cal Y}-Z.$ On $Z$
we take the trivial \tst, where $D^b(Z)^{\ge 0}=D^b(Z).$ 
On $U$ we take the standard
\tst, where $D^b(U)^{\ge 0}$ is the category of complexes 
of sheaves whose homology is supported in positive degrees. 
These {\it t}--structures glue to give
a \tst\ on  all of $D^b({\cal Y}).$ The actual simplex we 
started with,
in the proof of Theorem~\ref{safe kernels}, had a lifting
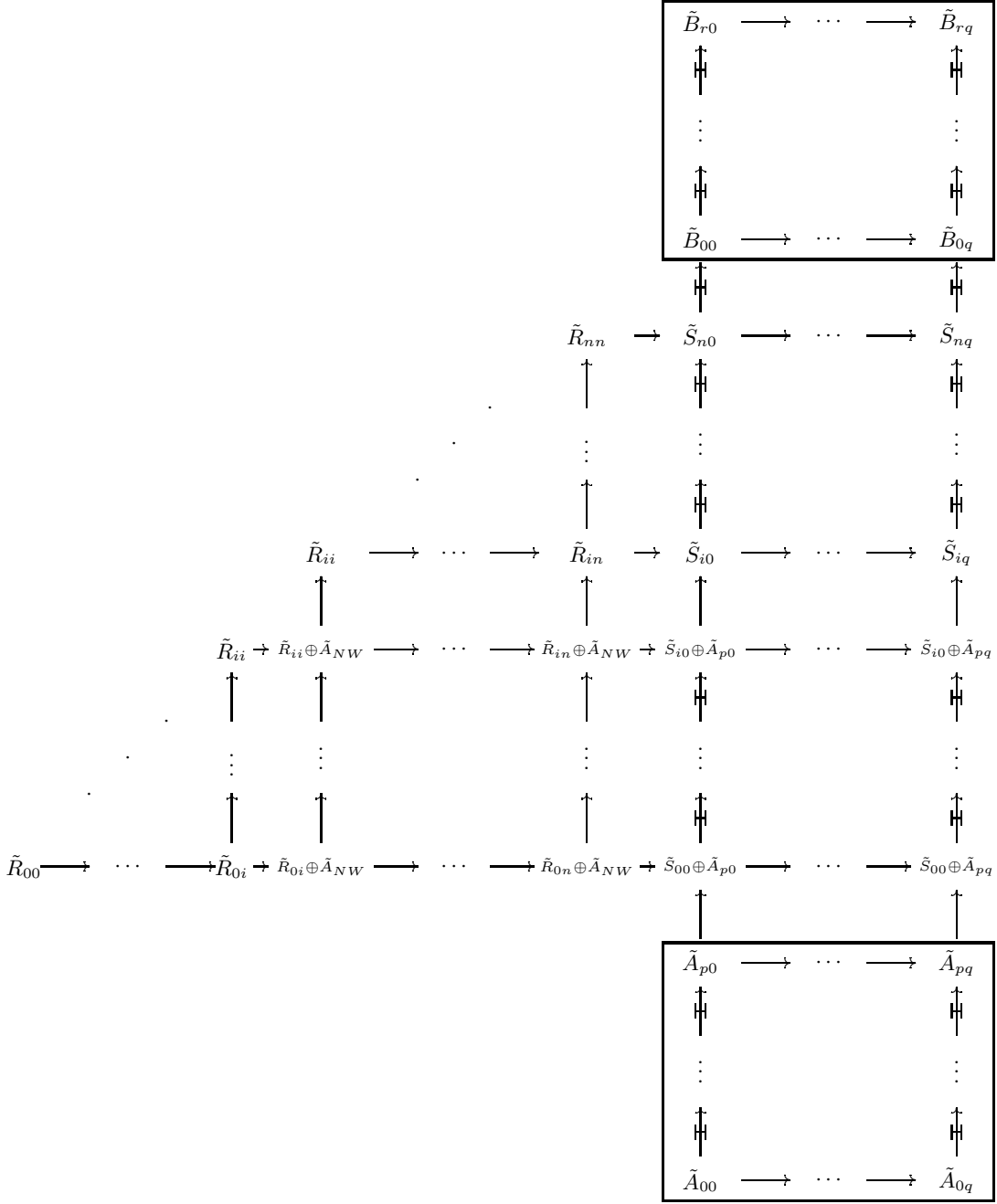
\small

\siv{10}1{15}1{14}0000
\sjv{11}2{11}2{11}2{11}00
\A10i
\A3in
\A50q
\B10p
\B30i
\B5in
\B70r

\SB\boxten{
                                      \Sbxp57ha{{\tilde B}_\Q}
                                      \Sbu56h
                    \fStrup35aa{{\tilde R}_\Q}\Sbr45a \Sbxp55ha{{\tilde S}_\Q}
                    \Sbu34a                   \Sbu54a
\fStrup13aa{{\tilde R}_\Q}\Sbrsmall23a
\Sbxpsmall33aa{\scriptstyle{\tilde R}_\Q\oplus {\tilde A}_{NW}}\Sbrsmall43a 
               \Sbxpsmall53ha{\scriptstyle{\tilde S}_{\I 0}
\oplus\tilde A_{p\J}}
                                           \Sbu52a
                                         \Sbxp51ha{{\tilde A}_\Q}
\fmbx5511 \fmbx5577
}
\BDI(0,500)
\UB(0,250)\boxten
\EDI

\normalsize

\eject
\vsize=8.4truein
\nin
and is thus a complex $\cal S$ of sheaves on $\cal X.$ 
The inclusion ${\cal X}
\hookrightarrow {\cal Y}$ has a continuous 
retraction $\pi:{\cal Y}
\longrightarrow {\cal X}$, sending $p$ to $\cs_{iq}.$ 
Then $\pi^*{\cal S}$
is a complex of sheaves on $\cal Y,$ which lies in 
$D^b({\cal Y})^{\ge 0}$
for our choice of perversity. The object $\alpha$ in 
$D^b({\cal Y})_{[0,0]}$, which is given by the pair of diagrams

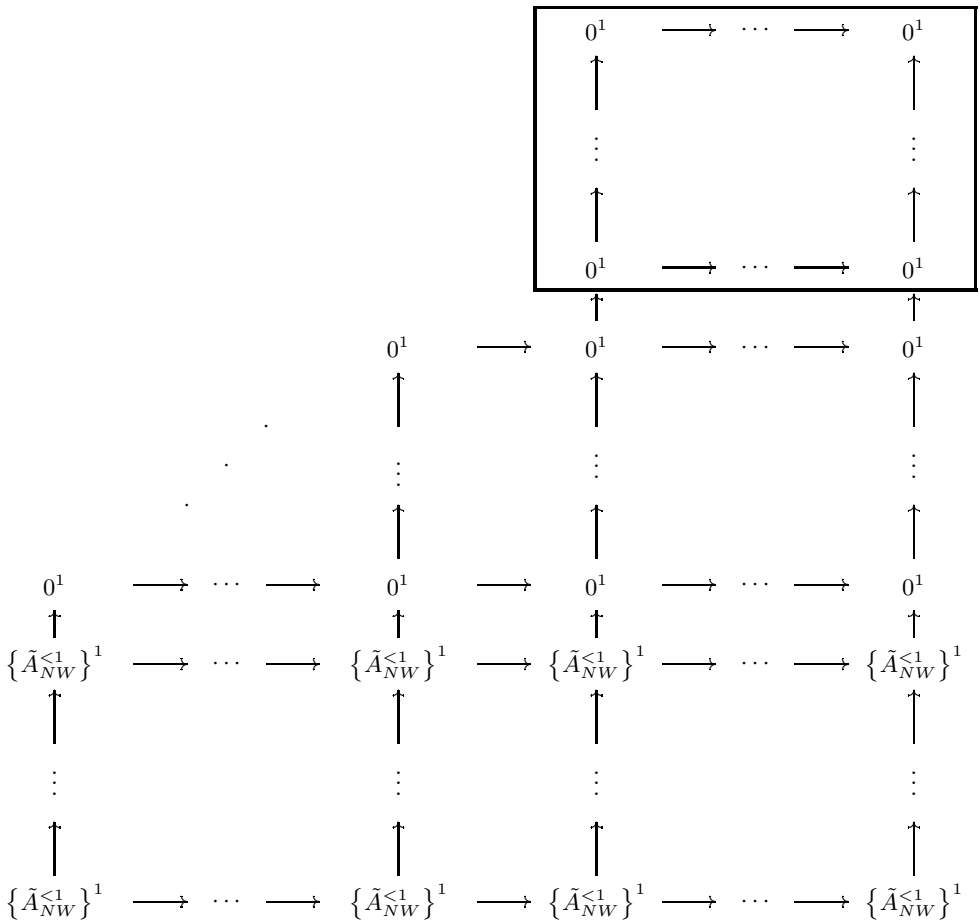
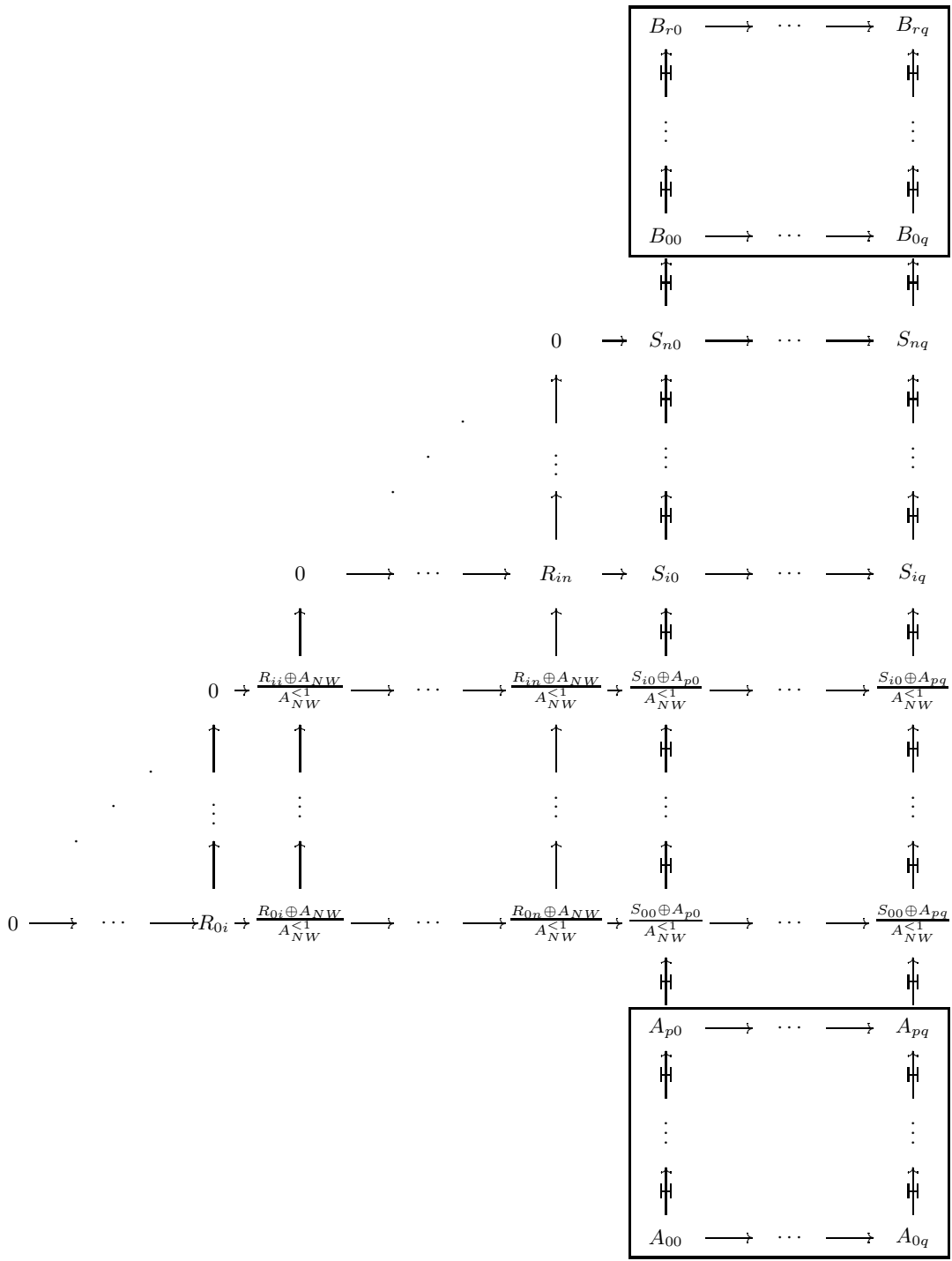
\small

\siv00{19}2{17}0000
\sjv00{11}1{11}1{11}00
\A10i
\A3in
\A50q
\B10p
\B30i
\B5in
\B70r

\SB\boxten{
                                      \Sbxp57aa{0^1}
                                      \Sbu56a
                    \fStrup35aa{0^1}\Sbr45a \Sbxp55aa{0^1}
                    \Sbu34a                   \Sbu54a
                   \Sbxp33aa{\left\{{\tilde A}_{NW}^{< 1}
                                               \right\}^1}\Sbr43a 
                    \Sbxp53aa{\left\{{\tilde A}_{NW}^{< 1}\right\}^1}
\fmbx5577
}
\BDI(0,350)
\UB(0,175)\boxten
\EDI

\normalsize

\nin
together with

\small

\siv{10}1{15}1{14}0000
\sjv{11}1{11}000000
\A10i
\A3in
\A50q
\B10p
\B30i
\B5in
\B70r

\SB\boxten{
\fStrup13aa{0^0}\Sbr23a\Sbxp33aa{0^0}\Sbr43a 
                                         \Sbxp53aa{0^0}
                                           \Sbu52a
                                         \Sbxp51aa{0^0}
\fmbx5511
}
\BDI(0,230)
\UB(0,115)\boxten
\EDI

\normalsize

\nin
injects into  $\pi^*{\cal S}.$  The quotient $\pi^*{\cal S}/\alpha$ 
is quasi--isomorphic to 
a complex of acyclic sheaves on $\cal Y,$ and in it we easily 
recognize a subset, giving a lifting of the 
diagram

\small

\siv{10}1{15}1{14}0000
\sjv{11}2{13}2{13}2{11}00
\A10i
\A3in
\A50q
\B10p
\B30i
\B5in
\B70r

\SB\boxten{
                                      \Sbxp57ha{B_\Q}
                                      \Sbu56h
                    \Strup35aa{R_\Q}\Sbr45a \Sbxp55ha{S_\Q}
                    \Sbu34a                   \Sbu54h
\Strup13aa{R_\Q}\Sbrsmall23a   
\Sbxpsmall33aa{{R_\Q\oplus A_{NW}}\over{A_{NW}^{<1}}}
\Sbrsmall43a 
                 \Sbxpsmall53ha{{S_{\I 0}\oplus A_{p\J}}\over{A_{NW}^{<1}}}
                                           \Sbu52h
                                         \Sbxp51ha{A_\Q}
\fmbx5511 \fmbx5577  
}
\BDI(0,540)
\UB(0,270)\boxten
\EDI

\normalsize

\nin
We deduce that the diagram has a lifting.

Every time we use \tst\ truncations in the 
proof, the Group 2 reader needs to go
through an argument similar to the one we have just given. Fortunately, we use 
\tst\ truncations only four times; once in the 
proof of Theorem~\ref{safe kernels}, as we 
have just seen. The next time we will see a homotopy using the truncation is in
Lemma~\ref{begin}. Then it occurs in
Lemma~\ref{kernels on side}. The fourth and last 
time we will see such a homotopy
will be in Lemma~\ref{a factoring of alpha}. 
Hopefully the reader will not object 
too vigorously, if I leave him to check for 
himself the other three occurrences of 
such an argument.
\ermk

{\it Continuation of the Proof of Theorem~\ref{safe kernels}.}
\label{page where} 
We are now
finished proving that the homotopy

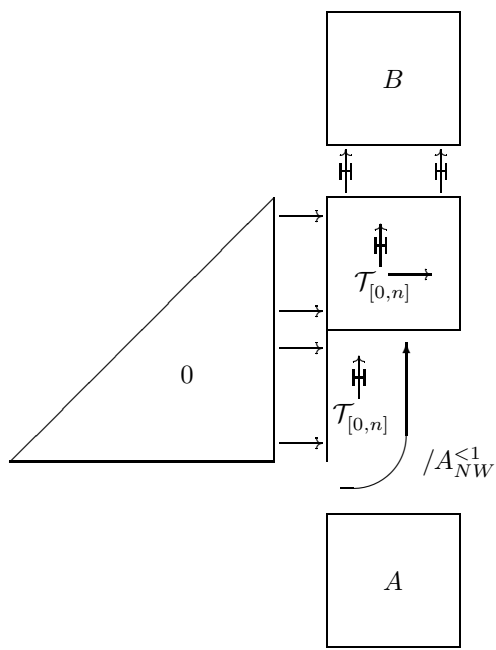
\sjv525525000
\siv{10}25000000

\SB\boxone{
                     \Fbx36B
                     \bu35h
\trui130  \br24a     \bx34ha0n
          \br23a     \htp33h0nu   \postr33{/A_{NW}^{<1}}e
                     \Fbx31A
}

\BC
\usebox\boxone
\EC

\nin
is well defined, even for the simplicial set that the Group 2 reader has been
working with. But this is not yet enough. The well--defined homotopy
above shows the identity is homotopic to some other map. 
We have to show that this other 
map is, in turn,
homotopic to the null map.

\nin
{\bf Warning for the Group 2 Reader.} As far as 
you are concerned, the next two 
paragraphs contain brazen, unmitigated lies. Remark~\ref{lies} will clarify
what the lies are, and how to fix them.

We have shown that the identity on

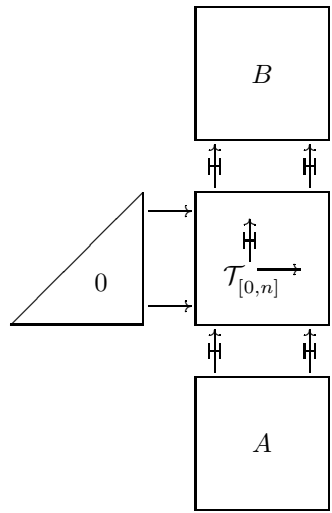
\siv525000000
\sjv525250000

\BP
                     \Fbx35B
                     \bu34h
\trui130  \br23a     \bx33ha0n
                     \bu32h
                     \Fbx31A
\EP

\nin
is homotopic to the simplicial map

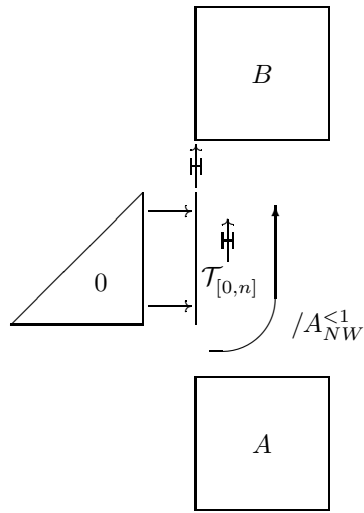
\siv525000000
\sjv525250000

\BP
                     \Fbx35B
                     \posar34hl
\trui130  \br23a     \htp33h0nu \postr33{/A_{NW}^{<1}}e
                     \Fbx31A
\EP

\nin
and this map factors through
the simplicial set

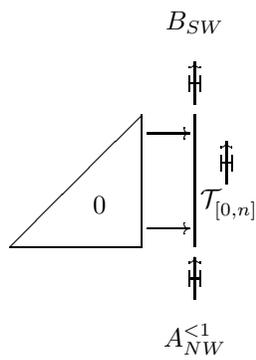
\siv522000000
\sjv325230000

\BP
                     \postr35{B_{SW}}l
                     \posar24he
\trui130   \br23a    \li33h0n{}l
                     \posar32hw
                     \postr31{A_{NW}^{<1}}l
\EP

\nin
But because kernels and cokernels agree,
the above simplicial set is equal to

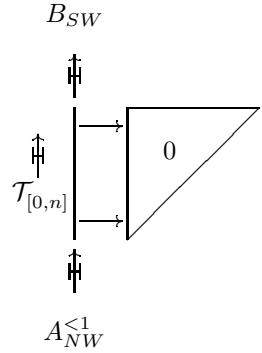
\siv002250000
\sjv325230000

\BP
                     \postr35{B_{SW}}r
                     \posar34he
                     \li33h0n{}r   \br43a \trdi530
                     \posar42hw
                     \postr31{A_{NW}^{<1}}r
\EP

\nin
and now, by the functoriality of the truncation, this set agrees with

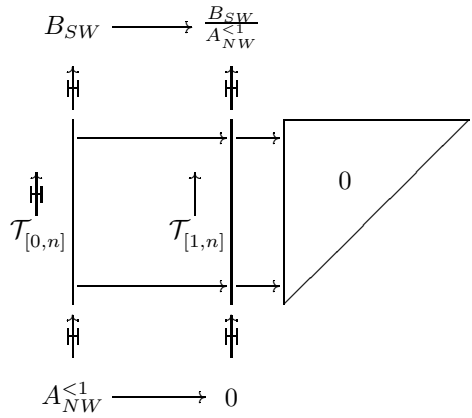
\siv242270000
\sjv327230000

\BP
\postr15{B_{SW}}r     \setra{35}{155}{30}a                \postr35{{B_{SW}}\over{A_{NW}^{<1}}}r
\posar14he                     \posar34he
\li13h0n{}r     \fbr143a                \li33a1n{}r   \br43a \trdi530
\posar22hw                     \posar42hw
\postr11{A_{NW}^{<1}}r      \setra{35}{15}{35}a               \postr31{0}r
\EP

\nin
and this makes it clear that the contraction to the initial 
object contracts the set.  \hfill{$\Box$}

\rmk{lies} As I said, as far as a Group 2 reader goes, this was a brazen lie.
What was wrong about it, and how can it be fixed?

As we correctly showed, the identity on

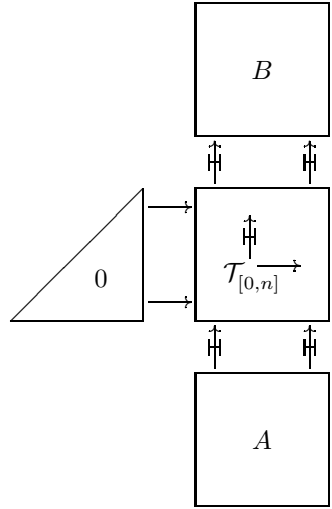
\siv525000000
\sjv525250000

\BP
                     \Fbx35B
                     \bu34h
\trui130  \br23a     \bx33ha0n
                     \bu32h
                     \Fbx31A
\EP

\nin
is homotopic to the simplicial map

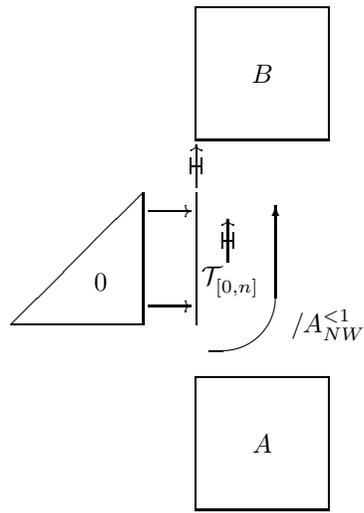
\siv525000000
\sjv525250000

\BP
                     \Fbx35B
                     \posar34hl
\trui130  \br23a     \htp33h0nu \postr33{/A_{NW}^{<1}}e
                     \Fbx31A
\EP

\nin
But when the author cavalierly suggested that the map factors through
the simplicial set

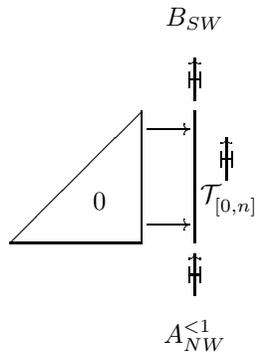
\siv522000000
\sjv325230000

\BP
                     \postr35{B_{SW}}l
                     \posar24he
\trui130   \br23a    \li33h0n{}l
                     \posar32hw
                     \postr31{A_{NW}^{<1}}l
\EP

\nin
he was blatantly lying, at least from the point 
of view of a Group 2 reader. The
simplicial map

\siv525000000
\sjv525250000

\BP
                     \Fbx35B
                     \posar34hl
\trui130  \br23a     \htp33h0nu \postr33{/A_{NW}^{<1}}e
                     \Fbx31A
\EP

\nin 
takes a simplex

\siv{13}2{13}000000
\sjv{13}2{13}2{13}0000

\A10n

\A30q
\B10p
\B30n

\B50r

\SB\boxone{
                           \Sbxp35ha{B_\Q}
                           \Sbu34h
\Strup13aa{R_\Q} \Sbr23a   \Sbxp33ha{S_\Q}
                           \Sbu32h
                           \Sbxp31ha{A_\Q}
\fmbx3355 \fmbx3311
}

\barr
s_n=\mtrx\boxone
\earr

\nin
to the simplex

\siv{13}2{17}000000
\sjv{13}2{13}2{13}0000

\A10n

\A30q
\B10p
\B30n

\B50r

\BP
                           \Sbxp35ha{B_\Q}
                           \Sbu34h
\Strup13aa{R_\Q} \Sbr23a   \Sbxp33ha{{S_{\I 0}\oplus A_{p\J}}\over{A_{NW}^{<1}}}
                           \Sbu32h
                           \Sbxp31ha{A_\Q}
\fmbx3355 \fmbx3311
\EP

\nin
Of course, this simplex is obtained from

\siv{13}2{17}000000
\sjv{13}2{13}2{13}0000

\A10n

\A30q
\B10p
\B30n

\B50r

\BP
                           \Sbxp35ha{B_\Q\oplus B_{00}}
                           \Sbu34h
\Strup13aa{R_\Q} \Sbr23a   \Sbxp33ha{{S_{\I 0}\oplus A_{p\J}}}
                           \Sbu32a
                           \Sbxp31ha{A_\Q}
\fmbx3355 \fmbx3311
\EP

\nin
by dividing every object in the top rectangle by the (diagonal)
inclusion of $B_{00}$, and everything in the middle rectangle
by the diagonal inclusion of $A_{NW}^{<1}$.
In other words, we obtain it from the direct sum 
of the fixed simplex

\siv{13}2{13}000000
\sjv{13}2{13}2{13}0000

\A10n

\A30q
\B10p
\B30n

\B50r

\BP
                           \Sbxp35ha{B_\Q}
                           \Sbu34h
\Strup13aa{0} \Sbr23a   \Sbxp33ha{{A_{p\J}}}
                           \Sbu32h
                           \Sbxp31ha{A_\Q}
\fmbx3355 \fmbx3311
\EP

\nin
and the simplex

\siv{13}2{17}000000
\sjv{13}2{13}2{13}0000

\A10n

\A30q
\B10p
\B30n

\B50r

\BP
                           \Sbxp35ha{B_{00}}
                           \Sbu34h
\Strup13aa{R_\Q} \Sbr23a   \Sbxp33ha{{S_{\I 0}}}
                           \Sbu32a
                           \Sbxp31ha{0}
\fmbx3355 \fmbx3311
\EP

\nin
This latest simplex, which contains all the information that varies as
we vary the integer $n,$ does indeed lie in the simplicial set

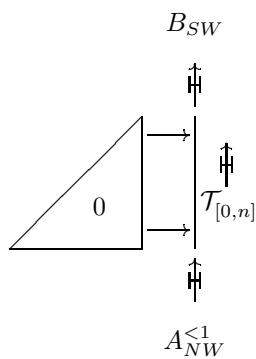
\siv522000000
\sjv325230000

\BP
                     \postr35{B_{SW}}l
                     \posar24he
\trui130   \br23a    \li33h0n{}l
                     \posar32hw
                     \postr31{A_{NW}^{<1}}l
\EP

\nin
So the reader may well wonder what all the fuss is about. 
Surely the map does factor
up to homotopy. All the variable information is contained in the smaller
simplicial set.

But there is a problem. Precisely, the problem is with the lifting to model
categories. The simplices

\siv{13}2{13}000000
\sjv{13}2{13}2{13}0000

\A10n

\A30q
\B10p
\B30n

\B50r

\BP
                           \Sbxp35ha{B_\Q}
                           \Sbu34h
\Strup13aa{0} \Sbr23a   \Sbxp33ha{{A_{p\J}}}
                           \Sbu32h
                           \Sbxp31ha{A_\Q}
\fmbx3355 \fmbx3311
\EP

\nin
and 

\siv{13}2{17}000000
\sjv{13}2{13}2{13}0000

\A10n

\A30q
\B10p
\B30n

\B50r

\BP
                           \Sbxp35ha{B_{0 0}}
                           \Sbu34h
\Strup13aa{R_\Q} \Sbr23a   \Sbxp33ha{{S_{\I 0}}}
                           \Sbu32a
                           \Sbxp31ha{0}
\fmbx3355 \fmbx3311
\EP

\nin
split as direct sums of diagrams having liftings to model categories. But
in fact, if we assume that

\siv{13}2{17}000000
\sjv{13}2{13}2{13}0000

\A10n

\A30q
\B10p
\B30n

\B50r

\BP
                           \Sbxp35ha{B_{0 0}}
                           \Sbu34h
\Strup13aa{R_\Q} \Sbr23a   \Sbxp33ha{{S_{\I 0}}}
                           \Sbu32a
                           \Sbxp31ha{0}
\fmbx3355 \fmbx3311
\EP

\nin
is just an arbitrary simplex in

\siv522000000
\sjv325230000

\BP
                     \postr35{B_{SW}}l
                     \posar24he
\trui130   \br23a    \li33h0n{}l
                     \posar32hw
                     \postr31{A_{NW}^{<1}}l
\EP

\nin
then these liftings will in general be to unrelated model categories, and the
careful reader will notice that in constructing the simplex

\siv{13}2{17}000000
\sjv{13}2{13}2{13}0000

\A10n

\A30q
\B10p
\B30n

\B50r

\BP
                           \Sbxp35ha{B_\Q}
                           \Sbu34h
\Strup13aa{R_\Q} \Sbr23a   \Sbxp33ha{{S_{\I 0}\oplus A_{p\J}}\over{A_{NW}^{<1}}}
                           \Sbu32h
                           \Sbxp31ha{A_\Q}
\fmbx3355 \fmbx3311
\EP

\nin
we had to assume a compatibility of liftings. Precisely, we needed to assume
that

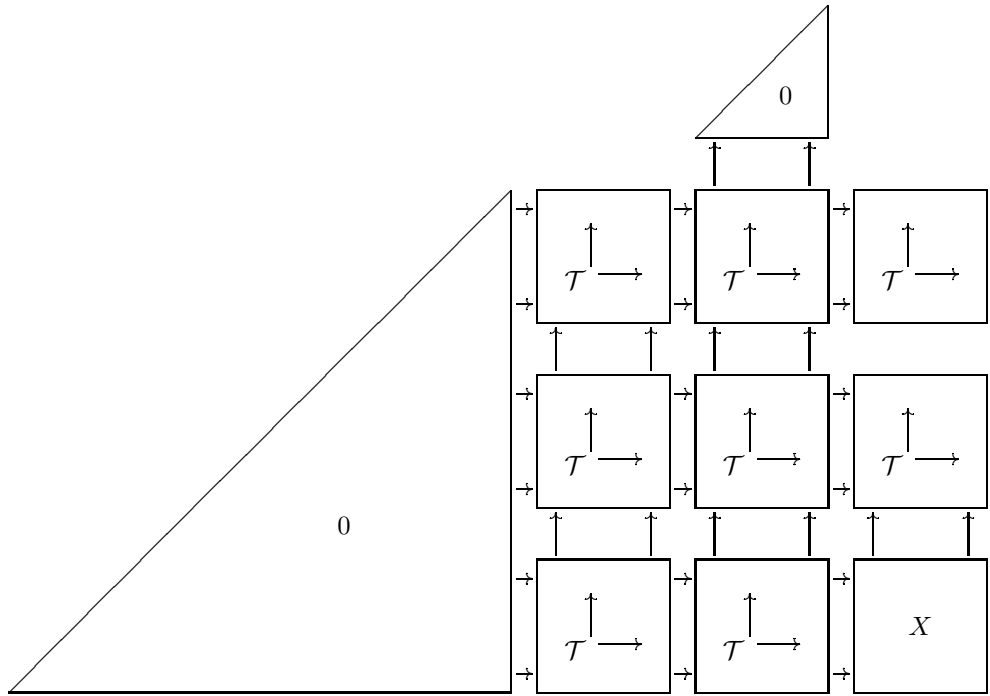
\siv{13}2{13}000000
\sjv{13}2{13}2{13}0000

\A10n

\A30q
\B10p
\B30n

\B50r

\BP
                           \Sbxp35ha{B_\Q}
                           \Sbu34h
\Strup13aa{0} \Sbr23a   \Sbxp33ha{{A_{p\J}}}
                           \Sbu32h
                           \Sbxp31ha{A_\Q}
\fmbx3355 \fmbx3311
\EP

\nin
and 

\siv{13}2{17}000000
\sjv{13}2{13}2{13}0000

\A10n

\A30q
\B10p
\B30n

\B50r

\BP
                           \Sbxp35ha{B_{0 0}}
                           \Sbu34h
\Strup13aa{R_\Q} \Sbr23a   \Sbxp33ha{{S_{\I 0}}}
                           \Sbu32a
                           \Sbxp31ha{0}
\fmbx3355 \fmbx3311
\EP

\nin
both split compatibly, that is one can choose a splitting of each so that the
summands correspond in pairs, and each pair lifts to the same model category.
One furthermore needs that the homology $H^0(S^r_{i 0})$ agree with the homology
$H^0(A^r_{p 0})$ for the $r^{th}$ summands $A^r_{p 0}$ of $A_{p 0}$ and
$S^r_{i 0}$ of $S_{i 0}.$

Having said what the difficulty is, it is also clear how to fix it. The point 
is that the contraction to the initial object preserves this direct sum
decomposition into pairs of simplices with liftings to the same model
category.

There is only one place in the remainder of the article where we will allow 
ourselves a similar lie. The reader is given notice that the same problem,
and the same solution, apply to the proof of Lemma~\ref{kernels on side}.

\ermk

\siv{19}15151500
\sjv525252500

\SB\boxeig{
                                     \trui570
                                      \bu56a
         \br25a  \fbx35\ct aa  \br45a    \fbx55\ct aa    \br65a   \fbx75\ct aa
                 \bu34a                \bu54a               
         \br23a   \fbx33\ct aa  \br43a    \fbx53\ct aa    \br63a   \fbx73\ct aa
                  \bu32a               \bu52a                   \bu72a
\trui110 \br21a   \fbx31\ct aa  \br41a    \fbx51\ct aa    \br61a   \Fbx71X
}

Now may be a good time to divulge a small secret. The blueprint simplicial set

\BC
\usebox\boxeig
\EC

\nin
that we have been considering until now, is unnecessarily large and clumsy.
Suppose we replace it by the smaller blueprint

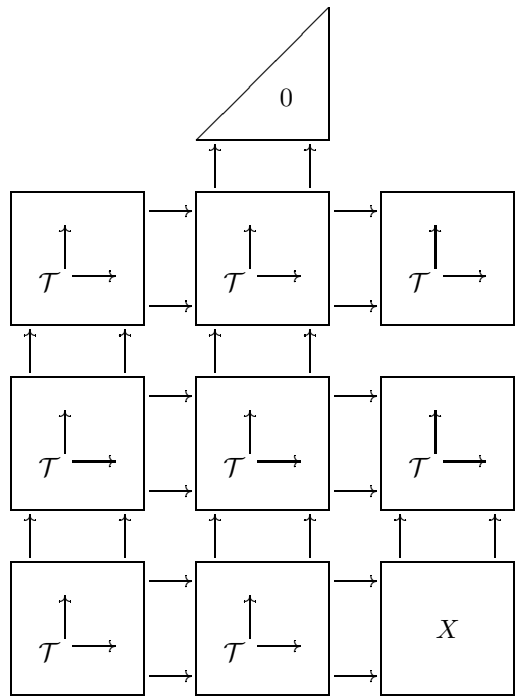
\siv005252500
\sjv525252500

\BP
                                     \trui570
                                      \bu56a
\fbx35\ct aa  \br45a    \fbx55\ct aa    \br65a   \fbx75\ct aa
                 \bu34a                \bu54a               
\fbx33\ct aa  \br43a    \fbx53\ct aa    \br63a   \fbx73\ct aa
                  \bu32a               \bu52a                   \bu72a
\fbx31\ct aa  \br41a    \fbx51\ct aa    \br61a   \Fbx71X
\EP

\nin
and we accept that on it the homotopy

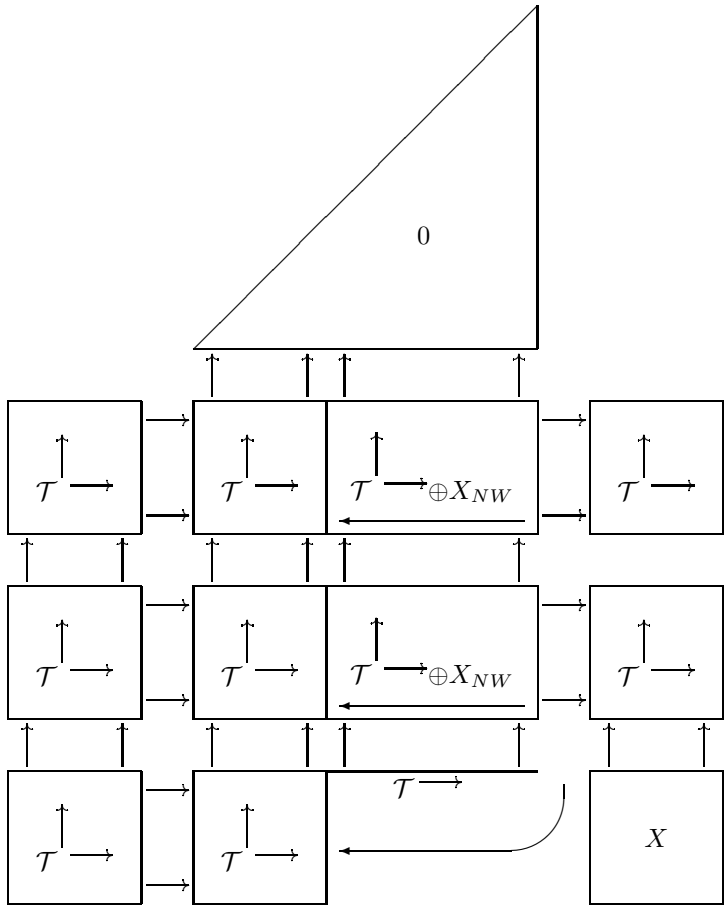
\siv005258250
\sjv525252{13}00

\BP
                                     \truj670
                                      \bu56a  \bu66a
\fbx35\ct aa  \br45a    \fbx55\ct aa \fhbx65\ct aa{ \oplus X_{NW}}l   \br75a   \fbx85\ct aa
                 \bu34a                \bu54a        \bu64a       
\fbx33\ct aa  \br43a    \fbx53\ct aa \fhbx63\ct aa{\oplus X_{NW}} l   \br73a   \fbx83\ct aa
                  \bu32a               \bu52a        \bu62a                \bu82a
\fbx31\ct aa  \br41a    \fbx51\ct aa  \fhtpl6611 \fli61\ct a{} t            \Fbx81X
\EP

\nin
is well--defined. Let us refer to it, for now, 
as the ``compact blueprint homotopy''. Then it is a formal 
consequence that the blueprint homotopy

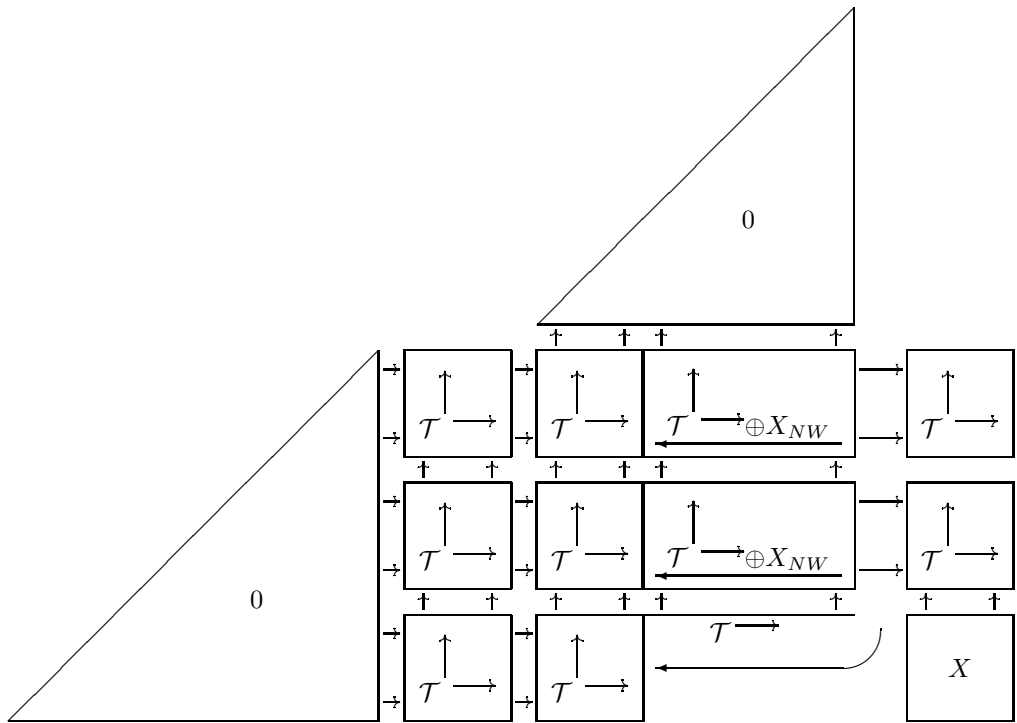
\siv{14}14148240
\sjv414141{12}00

\SB\boxten{
                                     \truj670
                                      \bu56a  \bu66a
         \br25a  \fbx35\ct aa  \br45a    \fbx55\ct aa \fhbx65\ct aa{ \oplus X_{NW}}l   \br75a   \fbx85\ct aa
                 \bu34a                \bu54a        \bu64a       
         \br23a   \fbx33\ct aa  \br43a    \fbx53\ct aa \fhbx63\ct aa{\oplus X_{NW}} l   \br73a   \fbx83\ct aa
                  \bu32a               \bu52a        \bu62a                \bu82a
\trui110 \br21a   \fbx31\ct aa  \br41a    \fbx51\ct aa  \fhtpl6611 \fli61\ct a{} t
            \Fbx81X
}
\BDI(0,270)
\UB(0,135)\boxten
\EDI

\nin
is also well defined.
How does one prove this fact? Simple. In the simplicial set

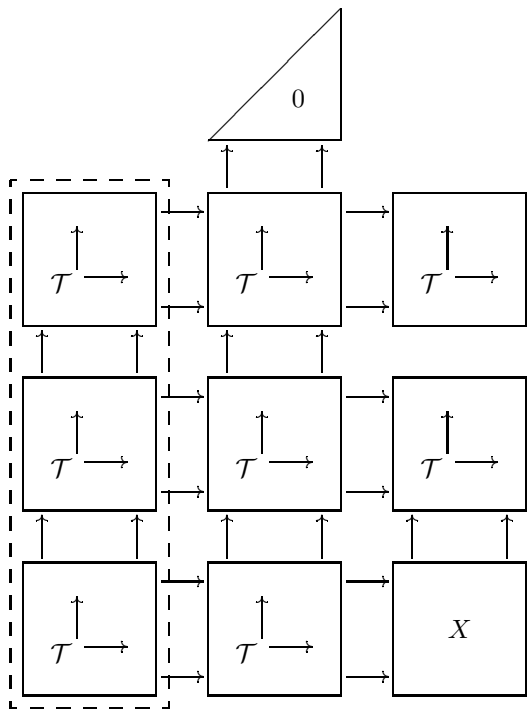
\siv005252500
\sjv525252500

\BP
                                     \trui570
                                      \bu56a
\fbx35\ct aa  \br45a    \fbx55\ct aa    \br65a   \fbx75\ct aa
                 \bu34a                \bu54a               
\fbx33\ct aa  \br43a    \fbx53\ct aa    \br63a   \fbx73\ct aa
                  \bu32a               \bu52a                   \bu72a
\fbx31\ct aa  \br41a    \fbx51\ct aa    \br61a   \Fbx71X
\dshbx3315
\EP

\nin
the part enclosed by a dashbox can be harmlessly subdivided. Precisely,  on the
simplicial set

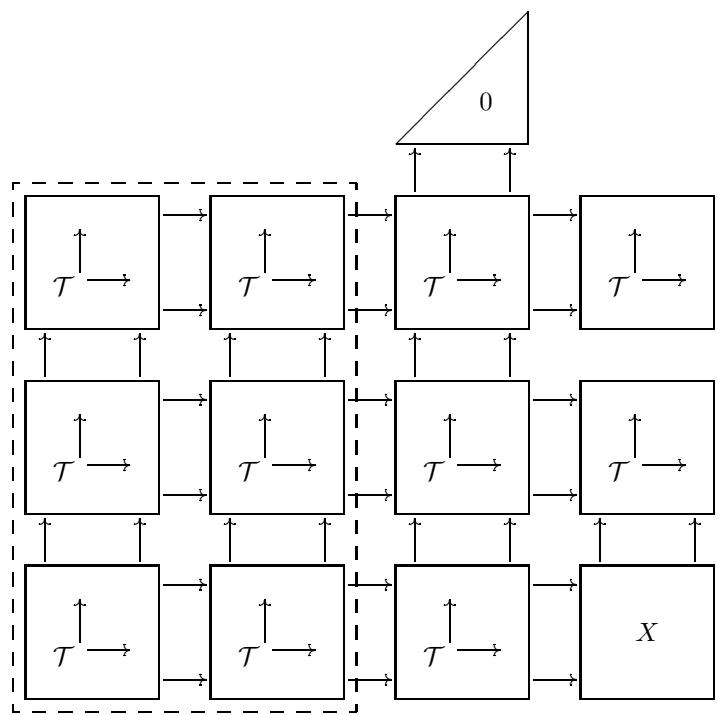
\siv{5}25252500
\sjv525252500

\BP
                                     \trui570
                                      \bu56a
\fbx15\ct aa    \br25a  \fbx35\ct aa  \br45a    \fbx55\ct aa    \br65a   \fbx75\ct aa
 \bu14a                \bu34a                \bu54a               
\fbx13\ct aa  \br23a   \fbx33\ct aa  \br43a    \fbx53\ct aa    \br63a   \fbx73\ct aa
\bu12a                  \bu32a               \bu52a                   \bu72a
\fbx11\ct aa \br21a   \fbx31\ct aa  \br41a    \fbx51\ct aa    \br61a   \Fbx71X
\dshbx1315
\EP

\nin 
there is a homotopy whose shorthand is simply

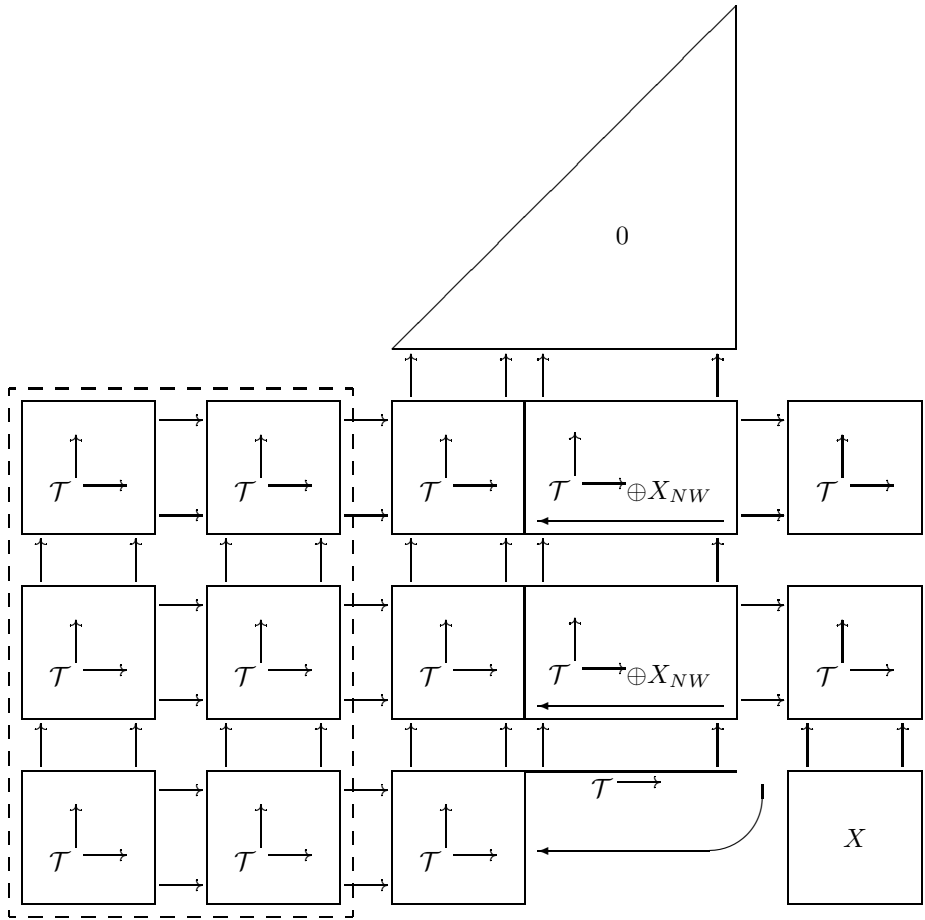
\siv{5}25258250
\sjv525252{13}00

\BP
                                     \truj670
                                      \bu56a  \bu66a
\fbx15\ct aa         \br25a  \fbx35\ct aa  \br45a    \fbx55\ct aa \fhbx65\ct aa{ \oplus X_{NW}}l   \br75a   \fbx85\ct aa
\bu14a                 \bu34a                \bu54a        \bu64a       
\fbx13\ct aa         \br23a   \fbx33\ct aa  \br43a    \fbx53\ct aa \fhbx63\ct aa{\oplus X_{NW}} l   \br73a   \fbx83\ct aa
\bu12a                  \bu32a               \bu52a        \bu62a                \bu82a
\fbx11\ct aa \br21a   \fbx31\ct aa  \br41a    \fbx51\ct aa  \fhtpl6611 \fli61\ct a{} t            \Fbx81X
\dshbx1315
\EP

\nin
since this is just a subdivided version of the compact blueprint homotopy.
But now the ordinary blueprint homotopy

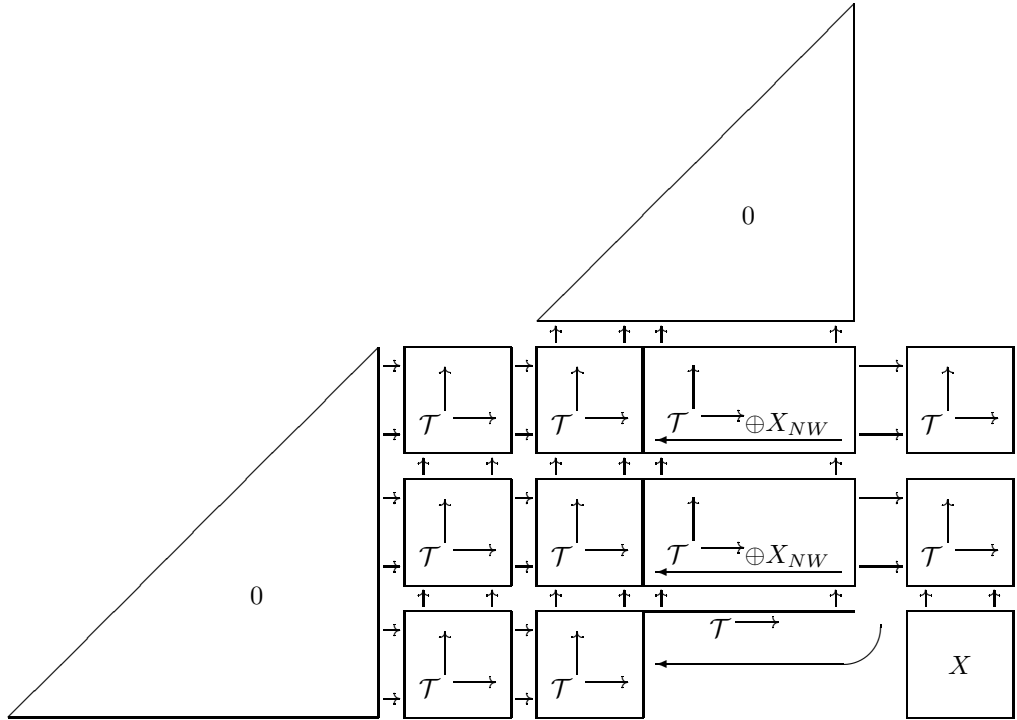
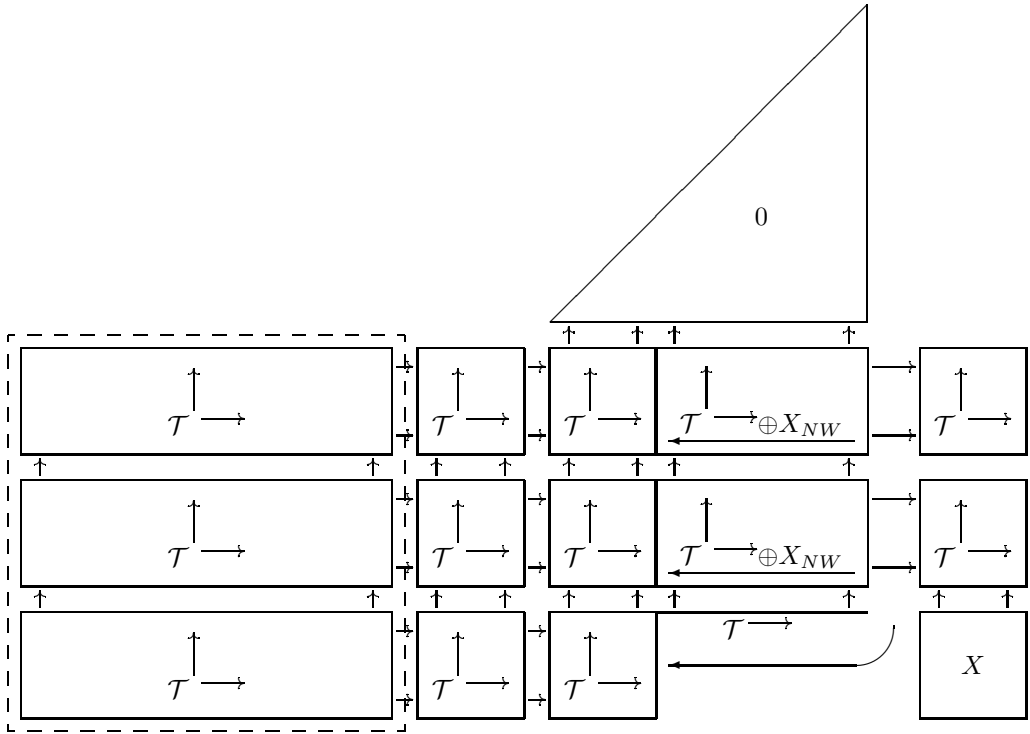
\siv{14}14148240
\sjv414141{12}00

\SB\boxten{
                                     \truj670
                                      \bu56a  \bu66a
         \br25a  \fbx35\ct aa  \br45a    \fbx55\ct aa \fhbx65\ct aa{ \oplus X_{NW}}l   \br75a   \fbx85\ct aa
                 \bu34a                \bu54a        \bu64a       
         \br23a   \fbx33\ct aa  \br43a    \fbx53\ct aa \fhbx63\ct aa{\oplus X_{NW}} l   \br73a   \fbx83\ct aa
                  \bu32a               \bu52a        \bu62a                \bu82a
\trui110 \br21a   \fbx31\ct aa  \br41a    \fbx51\ct aa  \fhtpl6611 \fli61\ct a{} t
            \Fbx81X
}
\BDI(0,290)
\UB(0,145)\boxten
\EDI

\nin
is obtained from

\siv{14}14148240
\sjv414141{12}00

\SB\boxten{
                                     \truj670
                                      \bu56a  \bu66a
\fbx15\ct aa \br25a  \fbx35\ct aa  \br45a    
\fbx55\ct aa \fhbx65\ct aa{ \oplus X_{NW}}l   \br75a   \fbx85\ct aa
 \bu14a                   \bu34a                \bu54a        \bu64a       
\fbx13\ct aa         \br23a   \fbx33\ct aa  
\br43a    \fbx53\ct aa \fhbx63\ct aa{\oplus X_{NW}} l   \br73a   \fbx83\ct aa
   \bu12a                 \bu32a               \bu52a        \bu62a   
             \bu82a
\fbx11\ct aa \br21a   \fbx31\ct aa  \br41a    
\fbx51\ct aa  \fhtpl6611 \fli61\ct a{} t
  \dshbx1115          \Fbx81X
}
\BDI(0,290)
\UB(0,145)\boxten
\EDI

\nin
by deleting some of the structure inside the dashbox. If we reflect back to the 
proof of Theorem~\ref{well defined homotopies}, it was based on the fact that
all squares are naturally Mayer--Vietoris. There are fewer squares in

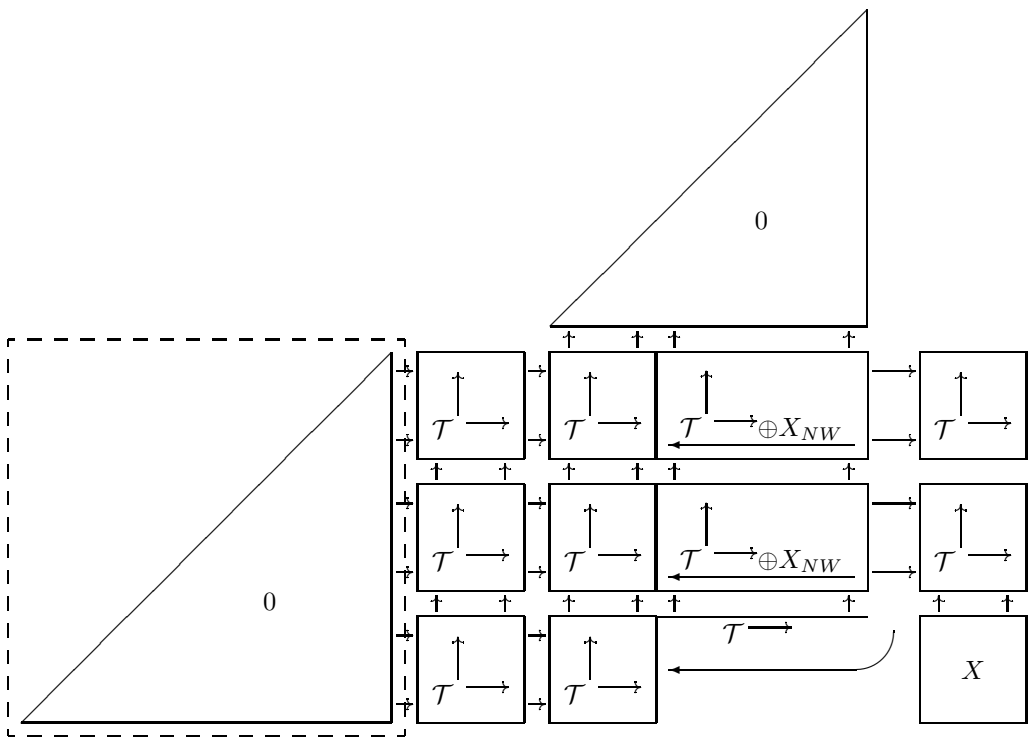
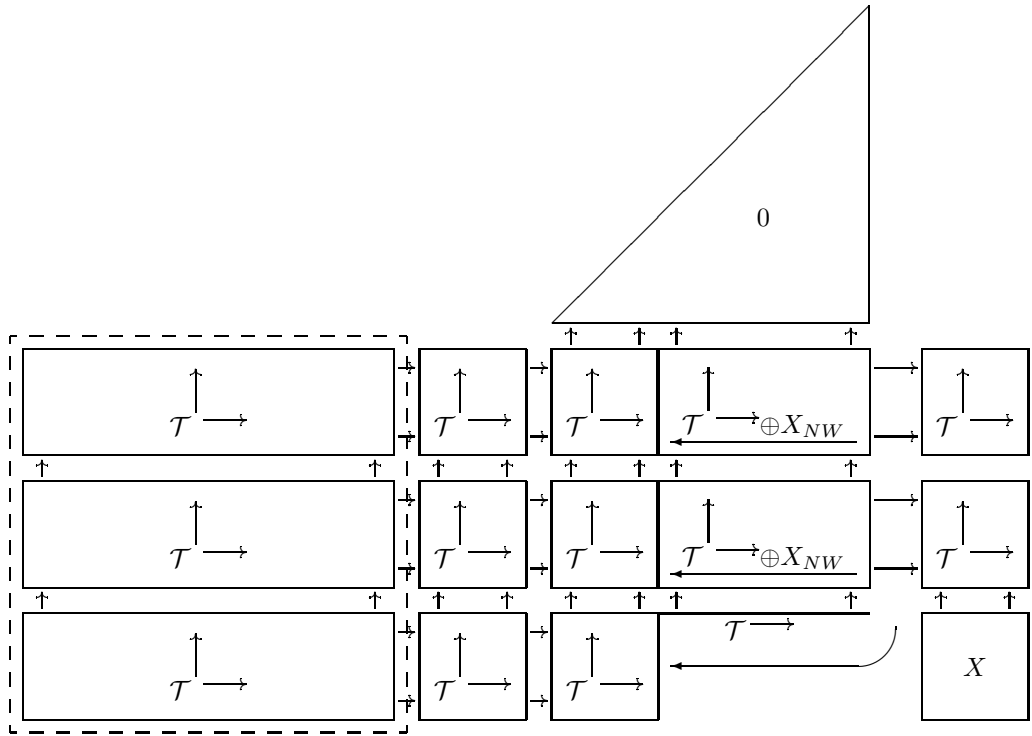
\siv{14}14148240
\sjv414141{12}00

\SB\boxten{
                                     \truj670
                                      \bu56a  \bu66a
         \br25a  \fbx35\ct aa  \br45a    \fbx55\ct aa \fhbx65\ct aa{ \oplus X_{NW}}l   \br75a   \fbx85\ct aa
                 \bu34a                \bu54a        \bu64a       
         \br23a   \fbx33\ct aa  \br43a    \fbx53\ct aa \fhbx63\ct aa{\oplus X_{NW}} l   \br73a   \fbx83\ct aa
                  \bu32a               \bu52a        \bu62a                \bu82a
\trui110 \br21a   \fbx31\ct aa  \br41a    \fbx51\ct aa  \fhtpl6611 \fli61\ct a{} t
   \dshbx1115            \Fbx81X
}
\BDI(0,290)
\UB(0,145)\boxten
\EDI

\nin
than in 

\siv{14}14148240
\sjv414141{12}00

\SB\boxten{
                                     \truj670
                                      \bu56a  \bu66a
\fbx15\ct aa \br25a  \fbx35\ct aa  \br45a    
\fbx55\ct aa \fhbx65\ct aa{ \oplus X_{NW}}l   \br75a   \fbx85\ct aa
 \bu14a                   \bu34a                \bu54a        \bu64a       
\fbx13\ct aa         \br23a   \fbx33\ct aa  
\br43a    \fbx53\ct aa \fhbx63\ct aa{\oplus X_{NW}} l   \br73a   \fbx83\ct aa
   \bu12a                 \bu32a               \bu52a        \bu62a   
             \bu82a
\fbx11\ct aa \br21a   \fbx31\ct aa  \br41a    
\fbx51\ct aa  \fhtpl6611 \fli61\ct a{} t
  \dshbx1115          \Fbx81X
}
\BDI(0,290)
\UB(0,145)\boxten
\EDI

\nin
because a triangle of squares is embedded in a rectangle. The fact that some
objects are restricted to be 0 in ordinary blueprint homotopy, but are 
free in the compact blueprint, only shows that the ordinary blueprint is even
more a special case of the compact blueprint than we might otherwise think.

There is no particularly good reason why I chose the blueprint homotopy to be the
one I gave. The compact blueprint homotopy does the job just as well, and it
can be made even more compact. The main point of this section is to convince the
reader, that the manipulations involved, in reducing a non--trivial homotopy to
a deletion of a subdivision of the blueprint, are essentially trivial. 
 From now on, we will feel free to leave this reduction to the reader.


\newpage

\mbox{ }


\begin{thebibliography}{99} 



\bibitem{BBD} {\sc A. A. Beilinson, J. Bernstein and P. Deligne}.  Analyse et 
topologie sur les {\'e}spaces singuliers.  Ast{\'e}risque 100, Soc. Math. France 
(1982).

\bibitem{B-S} {\sc A. Borel and J. P. Serre}.  Le th\'eor\`eme de RiemannRoch, Bull. 
Soc. Math. France 86 (1958), 97-136.
 
\bibitem{GN} {\sc C. Giffen,  and A. Neeman}.  $K$--theory for triangulated categories.  
Preprint.
 
\bibitem{H} {\sc R. Hartshorne}.  Residues and duality, SLN 20 (1966).
 
\bibitem{HS} {\sc V. A. Hinich and V. V. Schechtman}. Geometry of a category of complexes 
and algebraic $K$--theory.  Duke Math. J. 52 (1985), 399-430.
 
\bibitem{N1}  {\sc A. Neeman}.   Some new axioms for triangulated categories.  
J. 
of Algebra 139 (1992) 221--255.
 
\bibitem{N2}  {\sc A. Neeman}.  The Brown Representability Theorem and phantomless
triangulated categories. J. of Algebra 151 (1992) 118--155.
 
\bibitem{Q}  {\sc D. Quillen}.   Higher algebraic $K$--theory I, SLN 341 (1973), 85-147.

\bibitem {TT} {\sc R. Thomason and T. Trobaugh}. Higher algebraic K-theory of
schemes and of derived categories.
In: The Grothendieck Festschrift ( a collection 
of papers to honor Grothendieck's
60'th birthday) Volume 3 pp. 247--435, Birkh{\"a}user 1990.
 
\bibitem{V}  {\sc J. L. Verdier}.   Cat{\'e}gories d{\'e}riv{\'e}es, {\'e}tat 0.  SGA 
$4{1\over2}$, 262-308 (SLN 569, 1977).


 
\bibitem{W1}  {\sc F. Waldhausen}.  Algebraic $K$--theory of generalized free products I, 
II.  Ann. of Math. 108 (1978), 135-256.

\bibitem{W2} {\sc F. Waldhausen}. Algebraic K-theory of spaces.
SLN 1126 (1985), 318-419.

\end{thebibliography}
\end{document}